\documentclass[11pt]{amsart}
\usepackage[margin=1in]{geometry}
\usepackage{amsmath}
\usepackage{amsxtra}
\usepackage{amscd}
\usepackage{amsthm}
\usepackage{amssymb}
\usepackage[normalem]{ulem}
\usepackage{comment}
\usepackage{color}
\usepackage{tikz}
\usetikzlibrary{matrix,arrows,decorations.pathmorphing}
\usetikzlibrary{shapes.multipart}
\usepackage{hyperref}

\numberwithin{equation}{section}

\theoremstyle{definition}
\newtheorem{theorem}[equation]{Theorem}
\newtheorem{lemma}[equation]{Lemma}
\newtheorem{proposition}[equation]{Proposition}
\newtheorem{corollary}[equation]{Corollary}

\newtheorem{notation}[equation]{Notation}
\newtheorem{example}[equation]{Example}

\newtheorem{remark}[equation]{Remark}

\newtheorem{claim}[equation]{Claim}

\newcommand{\tworow}[2]{\genfrac{}{}{0pt}{}{#1}{#2}}

\newcommand{\C}{\mathbb{C}}

\newcommand{\Z}{\mathbb{Z}}
\newcommand{\N}{\mathbb{N}}

\newcommand{\del}{\partial}
\newcommand{\gdelg}{g\del_g}

\newcommand{\Hom}{\text{Hom}}

\newcommand{\Ext}{\text{Ext}}

\newcommand{\la}{\langle}
\newcommand{\ra}{\rangle}

\newcommand{\SL}{\mathsf{SL}}

\newcommand{\gr}{\operatorname{\mathsf{gr}}}

\newcommand{\onto}{\twoheadrightarrow}
\newcommand{\into}{\hookrightarrow}

\newcommand{\lra}{\longrightarrow}
\newcommand{\lla}{\longleftarrow}
\newcommand{\Mod}{\operatorname{Mod}}
\newcommand{\Spec}{\operatorname{Spec}}

\newcommand{\Rad}{\operatorname{Rad}}

\newcommand{\bN}{\mathbb{N}}
\newcommand{\mbf}{\mathbf}
\newcommand{\bb}{\mathbb}

\newcommand{\Hlam}{\mathsf{H_{\lambda}}}

\newcommand{\Hg}{\mathsf{H_g}}

\newcommand{\Hzero}{\mathsf{H_0}}

\newcommand{\Hab}{\bar{H}_{\lambda;\alpha,\beta}}
\newcommand{\Hzerob}{\bar{H}_{\lambda;0,\beta}}

\newcommand{\Mzerob}{S_{0,\beta}}
\newcommand{\Mzerobprime}{S_{0,\beta'}}
\newcommand{\Da}{\Delta_{\alpha}}

\newcommand{\Zgee}{\mathsf{Z_g}}

\newcommand{\Sab}{S_{\alpha,\beta}}
\newcommand{\Sabprime}{S_{\alpha,\beta'}}

\newcommand{\Weyl}{\mathsf{Weyl}}

\newcommand{\msf}{\mathsf}
\newcommand{\R}{\mathsf{R}}
\newcommand{\Zg}{\msf{Z_{g_1}}}
\newcommand{\Zlam}{\msf{Z_{\lambda}}}
\newcommand{\Zzero}{\msf{Z_0}}

\newcommand{\GL}{\mathsf{GL}}
\newcommand{\g}{ \begin{pmatrix} 1&1\\ 0&1\end{pmatrix}}
\newcommand{\E}{(\bb{Z}/p)^r}

\title{Symplectic reflection algebras in positive characteristic as Ore extensions}
\author{Emily Norton}
\begin{document}

\begin{abstract}
We investigate PBW deformations $\Hlam$ of $k[x,y]\rtimes G$ where $G$ is the cyclic group of order $p$ and $k$ also has characteristic $p$; in these deformations, $[x,y]$ takes a value in $kG$. The algebras $\Hlam$ are a version of symplectic reflection algebras that only exist in positive characteristic. They also happen to possess a presentation as an Ore extension over a commutative subring $\R$, and via the derivation defining the extension, have interesting connections to certain polynomials appearing in combinatorics and related to alternating permutations (the Andr\'{e} polynomials). We find the center of these algebras, their Verma modules, their simple modules, and the Ext groups between simples.  The Verma modules coincide with the fibers of $\Mod\Hlam\onto\Spec\R$, and $\Mod\Hlam$ turns out to be the disjoint union, in the sense of Smith and Zhang, of the Verma modules. As in non-defining characteristic, there are some distinctions between the $t=0$ and $t=1$ cases. We also obtain results when $G$ is replaced with $G^r$: in particular, we find the center and the simple modules when $[x,y]\in k[G^r]$.
\end{abstract}

\maketitle

\tableofcontents
\section{Introduction}
The goal of this paper is to study symplectic reflection algebras associated to a cyclic group of order $p$, and more generally an elementary abelian $p$-group, when the characteristic of the ground field is also $p$. The family of algebras $\Hlam$ as $\lambda$ ranges through $kG$ can be expressed as Ore extensions over their commutative subring $\R=k[x]\otimes kG$; this construction makes it possible to find the center of the algebras as well as to draw conclusions about the structure of the algebras and their modules.

For any $\lambda$, the algebra $\Hlam$ is Noetherian because $\R$ is Noetherian, and an Ore extension of a Noetherian ring is Noetherian.  The powers of $x^p$ (central elements) generate an infinite descending chain of two-sided ideals and thus $\Hlam$ is neither left nor right Artinian. It is not a domain since $kG$, which embeds into $\Hlam$, is not a domain.  It follows from the PBW theorem that $\Hlam$ has GK-dimension equal to that of $k[x,y]\rtimes G$, and hence has GK-dimension $2$.  $\Hlam$ has infinite global dimension.  The Jacobson radical of $\Hlam$ is $0$, so $\Hlam$ is semiprimitive (but not semisimple, since it is not Artinian). The only idempotents $\Hlam$ possesses are $0$ and $1$. The last two facts together imply that no simple $\Hlam$-module has a projective cover.

The algebra $\Hlam$ for a fixed $\lambda$ is always a finite module over its center. This is not obvious and takes some work to establish. Contrast this with the story when $kG$ is semisimple (say $k=\bb{C}$): the associated symplectic reflection algebra $\msf{H}_{t,c}$ has an important subalgebra called the spherical subalgebra, and the Satake isomorphism states that the center of $\msf{H}_{t,c}$ coincides with the center of the spherical subalgebra \cite{Gwyn}. The spherical subalgebra is defined as $e\msf{H}_{t,c}e$ where $e$ is the idempotent $\frac{1}{|G|}\sum_{g\in G}g$. When char $k=p$ divides the order of $G$, such an idempotent does not exist and there is no spherical subalgebra of $\Hlam$. Nonetheless the center $\Zlam$ of $\Hlam$ is large. The degrees of the generators of $\Zlam$ depend on whether $\lambda$ has a constant term or not. These cases correspond, respectively, to the $t=1$ and $t=0$ cases for symplectic reflection algebras in characteristic $0$. In characteristic $0$, the center of $\msf{H}_{t,c}$ is $\bb{C}$ if $t=1$ and is large if $t=0$. What we see in the algebras $\Hlam$ is that in the $t=1$ setting the center $\Zlam$ is isomorphic to a polynomial ring in two variables, generated by elements of degrees $p$ and $p^2$; in the $t=0$ case, $\Zlam$ is generated by elements of degrees $2,\;p,$ and $p$ with a singular relation between the first two generators.

Simple representations of $\Hlam$ are all finite-dimensional and are parametrized by $k\times k$. The simples are obtained as quotients of Verma modules; however, each Verma has a one-parameter family of maximal submodules, and so a one-parameter family of simple quotients; there is a one-parameter family of the Verma modules. The Verma modules are fiber modules of the morphism $\Spec\Hlam\onto\Spec\R$ induced by the inclusion $\R\into\Hlam$ ($\Spec\Hlam$ can be defined as a module category $\Mod\Hlam$ and is a noncommutative space \cite{SZ}). Computations of Ext groups between simple modules together with a theorem of Smith and Zhang \cite{SZ} imply that  $\Mod\Hlam$ is the disjoint union of its fibers.

Here is an outline of the contents of the sections of this paper. Section $2$ defines the algebras $\Hlam$ in imitation of the definition of a symplectic reflection algebra in dimension $2$, and observes that $\Hlam$ satisfies a PBW theorem. Section $3$ shows that $\Hlam$ is an Ore extension of its commutative subring $\R:=k[x]\otimes kG$ and establishes a general result about the combinatorics of $\R[y;\delta]$ for $\delta=f(x,g)\del_x+xg\del_g$ which will be necessary in proving that the center of $\Hlam$ is large. Section $4$ constructs explicit generators for the center of $\Hlam$ in both the $t=0$ and $t=1$ cases, showing that $\Hlam$ is always a finite module over its center. Section $5$ describes simple modules and Verma modules over $\Hlam$ and shows that $\Mod \Hlam$ is the disjoint union of its fibers. In contrast to what happens in characteristic $0$, the geometry of the center does not completely determine the dimensions of the simple modules, and the deformation parameter $\lambda\in kG$ has an effect on the outcome: when $t=0$ the Azumaya and smooth loci coincide only when $\lambda\in\Rad kG$. When $t=1$ the center is smooth but the dimension of simple modules drops over the $0$-fiber of the projection $\Mod\Hlam\onto\Spec\R$. Section $6$ generalizes the construction of $\Hlam$ from the cyclic group $G$ to an elementary $p$-group $E=\E$. For $r>1$, we first consider the case when $[x,y]\in E$, which reduces to the case $\lambda=g=\g$. We prove that the center of $\Hg$ in this case is generated by a quadratic element $x^2-2D^{p^r-2}(g),\;x^p,$ and $$y^{p^r}-y\frac{\delta^{p^r}(g)}{\delta(g)}=\prod_{\zeta\in\Xi}\prod_{a=1}^p(y-(\zeta+a)x)$$
where $D=\xi_i\sum g_i\del_i$, $\xi_i$ an $\bb{F}_p$ basis of $\bb{F}_q$, and where $\Xi$ consists of all elements of $\bb{F}_q$ orthogonal to $\bb{F}_p\subset\bb{F}_q$ thought of as the $\bb{F}_p$-linear subspace spanned by $1$. That is, the big central element is the product over all $e\in E$ of $e\cdot y$, taken in a certain order. We then look at simple $\Hg$-representations and show that for all $r>1$, the Azumaya locus and the smooth locus of the center coincide. As for the center of $\Hlam$ in the case of an arbitrary parameter $\lambda\in kE$, we show that in the $t=0$ case that $x^2-2D^{p^r-2}(\lambda)$, $x^p$, and $y^{p^r}-y\frac{\delta^{p^r}(g_1)}{\delta(g_1)}$ generate $\Zlam$, while in the $t=1$ case, $x^p$ and $\prod\limits_{a\in\bb{F}_p}\left(y^{p^r}-y\frac{\delta^{p^r}(g_1)}{\delta(g_1)}-ax^{p^r}\right)$ generate $\Zlam.$ We then look at the structure of Verma modules $\Da$ and show that simple modules obtained as quotients of $\Delta_0$ either have dimension $p$ or $1$, depending on whether $\lambda$ is invertible or a zero-divisor; simple quotients of $\Da$ when $\alpha\neq0$ have dimension equal to the degree of the big central element.

\textbf{Acknowledgements.} I would like to thank Victor Ginzburg for suggesting the problem, and for useful conversations. 

\section{$\Hlam$ as a symplectic reflection algebra in positive characteristic}
The name ``symplectic reflection algebra" derives from a phenomenon that is not necessarily true in characteristic $p$: a finite subgroup $G$ of $\SL(2,\bb{C})$ acting on $\bb{C}^2$ which preserves the symplectic form on $\bb{C}^2$ must be generated by ``symplectic reflections:" elements $s\in G$ such that the rank of $Id-s$ is $2$ \cite{EG}. In characteristic $p$, however, the cyclic subgroup $G$ of $\SL(2,\bb{F}_p)$ consisting of upper-triangular matrices with $1$'s on the diagonal is not generated by symplectic reflections but nonetheless is a finite subgroup which preserves the symplectic form on $\bb{F}_p^2$. The latter makes it possible to define the symplectic reflection algebra associated to $G$ exactly as it is defined in the $2$-dimensional case in characteristic $0$. Starting from such observations, Balagovic and Chen studied rational Cherednik algebras in characteristic $p$ for $\GL(n,p)$, obtained character tables for those algebras, and described irreducible representations in Category O for $\GL(2,p)$ \cite{Martina}. Our object of study is closely related to their work but the problem and methods are different: our construction attaches only a single copy of the underlying vector space and so our algebras are of the wreath product type rather than the rational Cherednik algebra type which have extra symmetry built into them (the two major styles of symplectic reflection algebras). And while a version of Category O can be defined for $\Hlam$, it's not as natural to put $\Hlam$-modules in that box, and its utility is unclear -- for example, it contains no projective objects (the argument in Remark 4.9 of \cite{Serga} shows this for $\Hlam$).

Let $k=\bar{k}$ be an algebraically closed field of characteristic $p$, $p$ odd. Let $G$ be the cyclic group of order $p$ realized as the unipotent subgroup of $\SL(2,p)$ generated by $g=\begin{pmatrix}
1&1\\
0&1 
\end{pmatrix}$. There is a natural action of $G$ on $k\la x,y\ra$, the algebra of noncommuting polynomials in $x$ and $y$, i.e. the tensor algebra of the two-dimensional vector space $V$ with coordinates $x$ and $y$; the action is given by $g^a\cdot f(x,y)=f(x,ax+y)$; and we can form the smash product algebra $k\la x,y\ra\rtimes G$. As a $k$-vector space $k\la x,y\ra\rtimes G\cong k\la x,y\ra\otimes kG$ but the multiplication in $k\la x,y\ra\rtimes G$ twists polynomials in $x$ and $y$ by the action of $G$ when an element of $kG$ moves past a polynomial: 
$$(f(x,y)\otimes \sum c_ag^a)(h(x,y)\otimes g^b)=\sum c_af(x,y)h(x,ax+y)\otimes g^{a+b}.$$

 Take $\lambda=c_0+c_1g+c_2g^2+...+c_{p-1}g^{p-1}\in kG$, a polynomial in $g$ of degree at most $p-1$. $G$ is abelian, so the center of $kG$ is $kG$. We define $\Hlam$ to be the quotient of $k\la x,y\ra\rtimes G$ by one quadratic relation which makes the commutator of $x$ and $y$ take the value $\lambda\in kG$: 
 $$\Hlam:=\frac{k\la x,y\ra\rtimes G}{(xy-yx-\lambda)}.$$
 \remark $x=gyg^{-1}-y$, so $\Hlam$ is generated by $y$ and $g$.
 
  For any $\lambda\in kG$, $\Hlam$ is a filtered algebra with the filtration given by assigning $x$ and $y$ degree $1$ and $G$ degree $0$. When $\lambda=0$ the algebra $\Hzero=k[x,y]\rtimes G$ is not only filtered but graded. A PBW theorem holds for $\Hlam$:
 
 \theorem (PBW theorem). For any $\lambda\in kG$, $\gr\Hlam=\Hzero$.\\
 
 The PBW theorem follows by checking the conditions given in the theorem by Griffeth which determines when a Drinfeld-Hecke algebra over a field of arbitrary characteristic satisfies a PBW theorem \cite{Griffeth}, and these conditions are easy to check for $\Hlam$. One may also prove it directly for $\Hlam$ (and other algebras similar to $\Hlam$) by adapting Bergman's Diamond Lemma \cite{Bergman}.The PBW also follows immediately from the fact that $\Hlam$ admits a description as an Ore extension -- see the next section. The PBW theorem allows us to write any element of the algebra $\Hlam$ in normal form. We will adopt the convention that $y$'s precede $x$'s precede $g$'s. Thus any element of $\Hlam$ can be written uniquely as a $k$-linear combination of monomials of the form $y^ix^jg^k$.
\notation Rescaling $\lambda$ by $t\in k$ has no effect on $\Hlam$: $$\Hlam\cong\msf{H}_{t\lambda}$$ and so, as in the study of symplectic reflection algebras more generally, it is useful to consider two cases: when $\lambda$ contains a constant term, and when $\lambda$ does not. In the former situation, we rescale $\lambda=\sum_{i=0}^{p-1}c_ig^i$ so that the constant term $c_0=1$, and call this the $t=1$ case. The latter situation where $c_0=0$ is called the $t=0$ case. 

It turns out that $\Hlam$ behaves like a symplectic reflection algebra in this respect: the precise shapes of the center and modules depend entirely on whether $t=0$ or $t=1$. However, we will see that in both cases $\Hlam$ is finite over its center, and $\Mod\Hlam$ is the disjoint union of Verma modules.

\section{$\Hlam$ as an Ore extension}
Another way to view $\Hlam$ is as an Ore extension of its commutative subring $$\R:=k[x]\otimes k[G]\cong kG[x].$$ An Ore extension is like a polynomial ring but with noncommutative relations between the added-on variable $y$ and the ground ring $\R$; a derivation and an automorphism of $\R$ control the relation between $\R$ and $y$. Formally, given an injective ring homomorphism $\sigma:\R\to \R$ and a $\sigma$-twisted derivation $\delta:\R\to \R$, that is, $\delta$ satisfying $$\delta(r_1r_2)=r_1\delta(r_2)+\delta(r_1)\sigma(r_2),$$ the (left) Ore extension of $\R$ with respect to $\sigma$ and $\delta$ is defined to be $\R[y]$ as a free right $\R$-module but with the multiplication rule $$ry=y\sigma(r)+\delta(r).$$ 
The notation for the Ore extension is $\R[y;\sigma, \delta]$. If $\sigma=Id$ then $\delta$ may be any derivation of $\R$ and the resulting Ore extension is simply denoted $\R[y;\delta]$. Ore extensions $\R[y;\delta]$ are called differential operator extensions, since the variable $y$ acts on $\R$ via the derivation $\delta$ and thus $\R[y;\delta]$ is an algebra of differential operators on $\R$ \cite{Good}.
\proposition $\Hlam\cong\R[y;\delta]$, the differential operator Ore extension over $\R=k[x]\otimes k[G]$ with $\delta$ the derivation defined on the generators $x$ and $g$ of $\R$ by $\delta(x)=\lambda$, $\delta(g)=xg$.
\begin{proof}
First, $\delta$ respects the relations defining $\R$: \begin{align*}&\delta(xg)=x^2g+\lambda g=\delta(gx)\\
 &\delta(g^p)=pxg^p=0=\delta(1)\end{align*}

The added-on variable $y$ in $\R[y;\delta]$ satisfies the relations \begin{align*}&xy=yx+\delta(x)=yx+\lambda\\ & gy=yg+\delta(g)=yg+xg\end{align*}
and these are the relations of $y$ with $x$ and $g$ in $\Hlam$.
\end{proof}
Since any Ore extension $\R[y;\sigma,\delta]$ over a $k$-algebra $\R$ is isomorphic to $\R[y]$ as a $k$-vector space, the PBW theorem for $\Hlam$ is an immediate consequence of this construction:
\corollary The PBW theorem holds for $\Hlam$: $$\gr\Hlam=\Hzero=k[x,y]\rtimes kG$$

\subsection{Powers of $\delta=f(x,g)\del_x+xg\del_g$ on $g$ in differential operator extensions of $kG[x]$}

The symplectic-reflection-type algebras $\Hlam$, viewed instead as Ore extensions over their commutative subring $\R=k[x]\otimes kG$, fit into a larger family of differential operator Ore extensions of $k[x]\otimes kG$. Namely, we may consider what happens more generally when $k=\bar{\bb{F}}_p$, $G$ is the cyclic group of order $p$, and the derivation $\delta$ takes the form $$\delta=f(x,g)\del_x+xg\del_g$$ with $f(x,g)$ a polynomial in $kG[x]$. Specifying $xg\del_g$ in the formula for $\delta$ means that these algebras are all quotients of $k\la x,y\ra \rtimes G$, so that in all the noncommutative algebras produced by such $\delta$, the group action remains present. 

\lemma A general formula: 

Let $\delta=f\del_x+xg\del_g$ where $f\in k[x,g]/(g^p-1)$. Then for any $n\geq1,$ $g$ divides $\delta^n(g)$ and \begin{align*}
\delta^n(g)=\delta\left(\frac{\delta^{n-1}(g)}{g}g\right)&=x\delta^{n-1}(g)+\delta\left(\frac{\delta^{n-1}(g)}{g}\right)g\\&=x\left( x\delta^{n-2}(g)+\delta\left(\frac{\delta^{n-2}(g)}{g}\right)g\right)+\delta\left(\frac{\delta^{n-1}(g)}{g}\right)g\\
&=x\left(x\left(   x\delta^{n-3}(g)+\delta\left(   \frac{\delta^{n-3}(g)}{g}  \right)g     \right)+\delta\left(\frac{\delta^{n-2}(g)}{g}\right)g\right)+\delta\left(\frac{\delta^{n-1}(g)}{g}\right)g\\&=...\\&=x^ng+x^{n-2}\delta\left(\frac{\delta(g)}{g}\right)g+x^{n-3}\delta\left(\frac{\delta^2(g)}{g}\right)g+...+x\delta\left(\frac{\delta^{n-2}(g)}{g}\right)g+\delta\left(\frac{\delta^{n-1}(g)}{g}\right)g
\end{align*}

 Here is the beginning of the sequence $\delta^n(g)$ written in terms of $\delta$-powers of $f$ and ordinary powers of $x$:

\begin{align*}
&\delta(g)=xg,\qquad\delta^2(g)=(f+x^2)g,\qquad\delta^3(g)=\left(\delta(f)+3fx+x^3\right)g,\\
&\delta^4(g)=\left(\delta^2(f)+3f^2+4\delta(f)x+6fx^2+x^4\right)g\\
&\delta^5(g)=\left( \delta^3(f)+10f\delta(f)+\left(5\delta^2(f)+15f^2\right)x+10\delta(f)x^2+10fx^3+x^5 \right)g\\
&\delta^6(g)=\big(\delta^4(f)+10\delta(f)^2+15f\delta^2(f) +15f^3 +\left(6\delta^3(f)+60f\delta(f)\right)x+\left(15\delta^2(f)+45f^2\right)x^2\\&\qquad\qquad\qquad+20\delta(f)x^3+15fx^4+x^6 \big)  g
\end{align*}

We visualize the powers of delta in a recurrence diagram. The $n$th row corresponds to $\delta^n$ and the $m$th column corresponds to the coefficient of $x^mg$ as a polynomial in $f,\delta(f),...$.

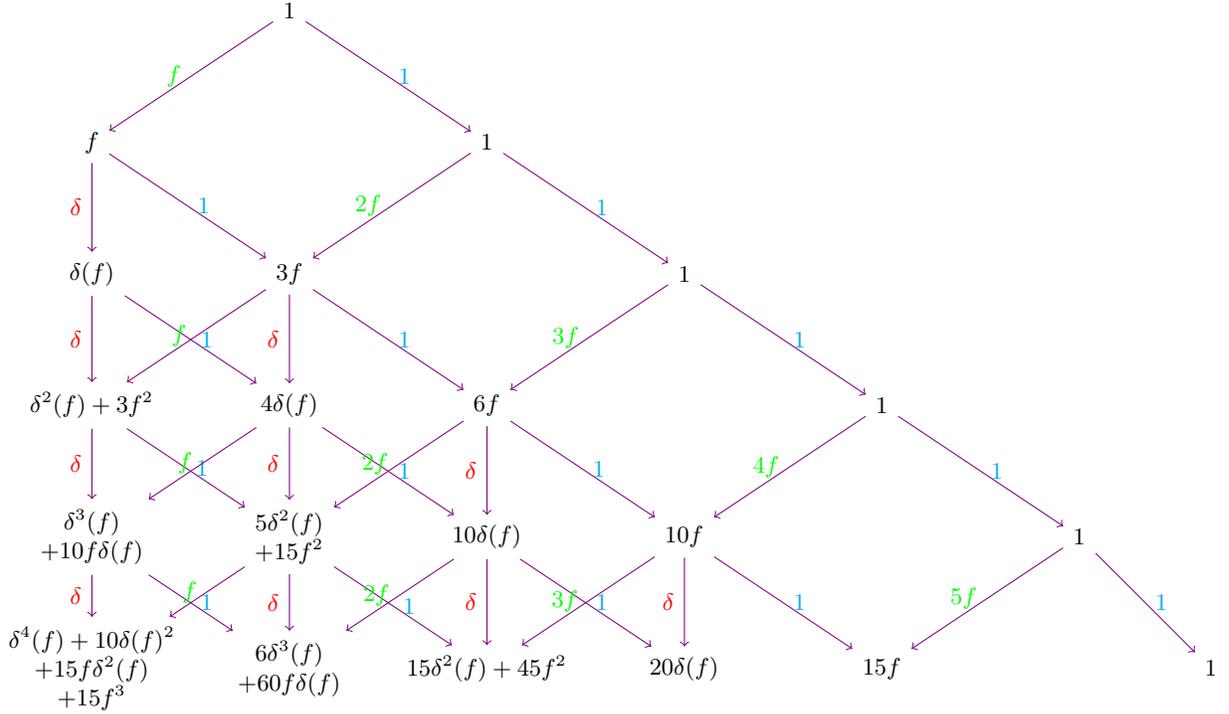
\begin{tikzpicture}[scale=1.75, font=\footnotesize,every text node part/.style={align=center}]
\node (A11) at (1,0) {$1$};
\node (A20) at (-.5,-1) {$f$};
\node (A22) at (2.5,-1) {$1$};
\node (A30) at (-.5,-2) {$\delta(f)$};
\node (A31) at (1,-2) {$3f$};
\node (A33) at (4,-2) {$1$};
\node (A40) at (-.5,-3) {$\delta^2(f)+3f^2$};
\node (A41) at (1,-3) {$4\delta(f)$};
\node (A42) at (2.5,-3) {$6f$};
\node (A44) at (5.5,-3) {$1$};
\node (A50) at (-.5,-4) {$\delta^3(f)$\\$+10f\delta(f)$}; 
\node (A51) at (1,-4) {$5\delta^2(f)$\\$+15f^2$};
\node (A52) at (2.5,-4) {$10\delta(f)$};
\node (A53) at (4,-4) {$10f$};
\node (A55) at (7,-4) {$1$};
\node (A60) at (-.5,-5) {$\delta^4(f)+10\delta(f)^2$\\$+15f\delta^2(f)$\\$ +15f^3$};
\node (A61) at (1,-5) {$6\delta^3(f)$\\$+60f\delta(f)$};
\node (A62) at (2.5,-5) {$15\delta^2(f)+45f^2$};
\node (A63) at (4,-5) {$20\delta(f)$};
\node (A64) at (5.5,-5) {$15f$};
\node (A66) at (8,-5) {$1$};
\path[->,violet,font=\footnotesize,]
(A11) edge node [green] [left] {$f$} (A20)
(A11) edge node [cyan] [right] {$1$} (A22)
(A20) edge node [red] [left] {$\delta$} (A30)
(A20) edge node [cyan] [right] {$1$} (A31)
(A22) edge node [green] [left] {$2f$} (A31)
(A22) edge node [cyan] [right] {$1$} (A33)
(A30) edge node [red] [left] {$\delta$} (A40)
(A30) edge node [cyan] [right] {$1$} (A41)
(A31) edge node [green] [left] {$f$} (A40)
(A31) edge node [red] [left] {$\delta$} (A41)
(A31) edge node [cyan] [right] {$1$} (A42)
(A33) edge node [green] [left] {$3f$} (A42)
(A33) edge node [cyan] [right] {$1$} (A44)
(A40) edge node [red] [left] {$\delta$} (A50)
(A40) edge node [cyan] [right] {$1$} (A51)
(A41) edge node [red] [left] {$\delta$} (A51)
(A41) edge node [cyan] [right] {$1$} (A52)
(A41) edge node [green] [left] {$f$} (A50)
(A42) edge node [red] [left] {$\delta$} (A52)
(A42) edge node [cyan] [right] {$1$} (A53)
(A42) edge node [green] [left] {$2f$} (A51)
(A44) edge node [cyan] [right] {$1$} (A55)
(A44) edge node [green] [left] {$4f$} (A53)
(A50) edge node [red] [left] {$\delta$} (A60)
(A50) edge node [cyan] [right] {$1$} (A61)
(A51) edge node [green] [left] {$f$} (A60)
(A51) edge node [red] [left] {$\delta$} (A61)
(A51) edge node [cyan] [right] {$1$} (A62)
(A52) edge node [green] [left] {$2f$} (A61)
(A52) edge node [red] [left] {$\delta$} (A62)
(A52) edge node [cyan] [right] {$1$} (A63)
(A53) edge node [green] [left] {$3f$} (A62)
(A53) edge node [red] [left] {$\delta$} (A63)
(A53) edge node [cyan] [right] {$1$} (A64)
(A55) edge node [green] [left] {$5f$} (A64)
(A55) edge node [cyan] [right] {$1$} (A66);
\end{tikzpicture}

\proposition $\delta^p(g)=\delta^{p-2}(f)g+x^pg=\delta^{p-1}(x)g+x^pg$.

\begin{proof}
The diagram together with the formula in the lemma lead to the following expressions for the first few powers of $\delta$, from which the general pattern is clear, and which explains the familiar numbers increasing down diagonals:  

\begin{align*}
\delta(g)&=xg\\
\delta^2(g)&=x^2g+fg\\
\delta^3(g)&=x^3g+{\color{violet}(1+2)}fxg+\delta(f)g\\
\delta^4(g)&=x^4g+{\color{violet}(1+2+3)}fx^2g+{\color{magenta}(1+3)}\delta(f)xg+\left(\delta^2(f)+3f^2\right)g\\
\delta^5(g)&=x^5g+{\color{violet}(1+2+3+4)}fx^3g+{\color{magenta}(1+3+6)}\delta(f)x^2g+{\color{red}(1+4)}(\delta^2(f)+3f^2)xg+(\delta^3+10f\delta(f))g\\
\delta^6(g)&=x^6g+{\color{violet}(1+2+3+4+5)}fx^4g+{\color{magenta}(1+3+6+10)}\delta(f)x^3g+{\color{red}(1+4+10)}\left(\delta^2(f)+3f^2\right)x^2g\\
&\qquad+{\color{orange}(1+5)}\left(\delta^3(f)+10f\delta(f)\right)xg+\left(\delta^4(f)+10\delta(f)^2+15f\delta^2(f)+15f^3\right)g\\
...\\
\delta^n(g)&=x^ng+{\color{violet}\sum_{k=1}^{n-1}{k\choose 1}}fx^{n-2}g+{\color{magenta}\sum_{k=2}^{n-1}{k\choose 2}}\delta(f)x^{n-3}g\\&\qquad\qquad\qquad\qquad\qquad\qquad\qquad+{\color{red}\sum_{k=3}^{n-1}{k\choose 3}}\left(\delta^2f+3f^2\right)x^{n-4}g+...+{\color{blue}\sum_{k=n-2}^{n-1}{k\choose n-2}}F_{n-1}xg+F_ng\\
&=x^ng+{\color{violet}{n\choose2}}fx^{n-2}g+{\color{magenta}{n\choose3}}\delta(f)x^{n-3}g+{\color{red}{n\choose4}}\left(\delta^2(f)+3f^2\right)x^{n-4}g+...+{\color{blue}{n\choose n-1}}F_{n-1}xg+F_ng\\
\end{align*}

We have let $F_j$ be the $kG[x]$-sequence of polynomials in $f,\delta(f), \delta^2(f),...$ that appears along the left vertical edge of the diagram above, i.e. the coefficient of $x^0$ in this latest formula. Since the nontrivial binomial coefficients for $n=p$ vanish mod $p$, it follows that:

$$\delta^p(g)=x^pg+F_pg$$

All that remains is to understand the sequence $F_j$, which starts from $j=2$. A look at the recursion diagram for $\delta^n(g)$ indicates the fact that the sequence $F_n$ satisfies the recursion
$$F_n=\delta(F_{n-1})+(n-1)fF_{n-2}$$

\notation As a shorthand, the partition $\mbf{\alpha}:=(2^{a_2},3^{a_3},4^{a_4},...,(j+2)^{a_{j+2}})$ denotes the monomial $f^{a_2}\delta(f)^{a_3}\delta^2(f)^{a_4}\cdot\cdot\cdot\delta^j(f)^{a_{j+2}}.$  We denote the integer being partitioned by $|\mbf{\alpha}|$: $|\mbf{\alpha}|=2a_2+...+(j+2)a_{j+2}$, and write $\alpha\vdash|\alpha|$ to say $\alpha$ is a partition of $|\alpha|$. The exponent $a_j$ in a partition denotes that part $j$ being repeated $a_j$ times; for example, $(2^3,3,4^4)$ denotes the partition $(2,2,2,3,4,4,4,4).$

The partitions that appear in $F_n$ are all the partitions of $n$ with $k-1$ parts, for all $k>0$ such that $n-2k\geq0$. Indeed, this follows by induction on $n$ together with the recursion formula for $F_n$: $\delta$ raises $|\mbf{\alpha}|$ by $1$ while multiplication by $f$ corresponds to padding the partition with $2$. Moreover, the sums $$S(n):=\mbox{ sum of the coefficients of }F_n$$ form the sequence \href{http://oeis.org/A000296}{A000296}: $S_n$ counts the number of ways $n$ people can arrange themselves into cliques, where a clique contains at least two people.

By abuse of notation, write $F_n$ as the linear combination of the partitions associated to the monomials appearing in $F_n$. Then:

\begin{align*}
F_2&=(2)\\
F_3&=(3)\\
F_4&=(4)+3(2^2)\\
F_5&=(5)+10(2,3)\\
F_6&=(6)+15(2,4)+10(3,3)+15(2^3)\\
F_7&=(7)+21(2,5)+35(3,4)+105(2^2,3)\\
F_8&=(8)+28(2,6)+56(3,5)+35(4^2)+280(2,3^2)+210(2^2,4)+105(2^4)\\
F_9&=(9)+36(2,7)+84(3,6)+126(4,5)+280(3^3)+1260(2,3,4)+378(2^2,5)+1260(2^3,3)\\
...
\end{align*}

Now we notice that the coefficient of $\mbf{\alpha}=(2^{a_2},...,n^{a_n})\vdash n$ counts the number of ways $n$ people can be divided up into $a_2$ cliques of size $2$, $a_3$ cliques of size $3$,... The formula for the coefficients of the partition summands of $F_n$ is therefore: 

$$F_n=\sum_{\tworow{\mbf{\alpha}\vdash n}{a_i>0\implies i>1}} \frac{n!}{(2!^{a_2}3!^{a_3}\cdot\cdot\cdot n!^{a_n})(a_2!a_3!\cdot\cdot\cdot a_n!)}(2^{a_2},3^{a_3},...,n^{a_n})$$

If $n=p$ is prime it follows that all of these coefficients except $1$ vanish, and we conclude: $$F_n\equiv\delta^{p-2}(\lambda)\equiv\delta^{p-1}(x)\mbox{ mod }p$$
\end{proof}

Perhaps more general theorems that can be proved about Ore extensions of $kG[x]$ as above, about their centers, structure, representation theory -- we have some ideas -- but this paper will remain focused on the family of algebras $\Hlam$. We owe the idea to use partitions as a shorthand for polynomials in $\delta$-powers of $f$ (or $\lambda$, where $\Hlam$ is concerned), to the second Benkart-Lopes-Ondrus paper on Ore extensions of $k[x]$ \cite{Benk2}, and there are similarities between their formulas and ours. It is interesting to compare the combinatorics of $\delta=h(x)\del_x$ for Ore extensions $k[x;\delta]$ with the combinatorics of $\delta=f(x,g)\del_x+xg\del_g$ for Ore extensions of $kG[x][y;\delta]$ as above. For instance, in the table of coefficients following Corollary 9.5 in \cite{Benk2}, we see some numbers that are familiar from the case of $\lambda=g$ (see Section 4.3 below of our paper):  the sequence composed of the second entry from the right in each row of their table is the sequence of Eulerian numbers \href{http://oeis.org/A000295}{A000295}, which is the second row of the diagram of Andr\'{e} numbers (cf. 4.3). The sequence composed of the third entry from the right in each row of their table is the sequence of reduced tangent numbers\href{http://oeis.org/A002105}{A002105}, which is the leftmost column of the diagram of Andr\'{e} numbers (cf. 4.3). Moreover, the sum of the entries in the $n$th row of Benkart et al's table is $(n-1)!$ (\cite{Benk2}, Proposition 9.7); when $\lambda=g^2$, the sum of the coefficients of $\delta^n(g)$ for $\delta=\lambda\del_x+xg\del_g$ is $n!$: see Section 4.3 below.

\section{The Center of $\Hlam$}
\subsection{Powers of $\delta$ on $g$ and on $x$: partition recurrence diagrams}
Consider the effect of $$\delta=\lambda\del_x+xg\del_g$$ applied successively to the element $g$. It is easy to check that $\delta^n(g)$ is divisible by $g$ for all $n\geq0$ and that, considered as a polynomial in $kG[x]$, it has degree $n$; moreover, only the coefficients of $x^{n-2k}$, $n-2k\geq0$, are nonzero.  For a given monomial $C_{n,m}x^mg$ in $\delta^n(g)$ (it is convenient to factor out the omnipresent extra $g$ from the $kG$-coefficient), the $kG$-coefficient $C_{n,m}$ obeys the recurrence relation:
\begin{align*}
C_{0,0}&=1,\qquad C_{n,m}=0\mbox{ if } n<m\mbox{ or } m<0\\
C_{n,m}&=(D+1)(C_{n-1,m-1})+(m+1)\lambda C_{n-1,m+1}\mbox{ for }n>0
\end{align*}
where $$D:=\gdelg.$$ The $C_{n.m}$ are polynomials in $\lambda,\;D(\lambda),\;D^2(\lambda),...,\;D^{p-2}(\lambda)$. Such polynomials may be represented by linear combinations of partitions. To this end, we assign $\lambda$ the value $1$, $D(\lambda)$ the value $2$, and so forth, so that $\lambda^{a_1}D(\lambda)^{a_2}\cdot\cdot\cdot D^{p-2}(\lambda)^{a_{p-2}}$ corresponds to the partition $((p-2)^{a_{p-2}},...,2^{a_2},1^{a_1})$. The monomials produced by the recurrence equation, and thus the powers of $\delta$ on $g$, may be represented by a directed graph whose vertices are the monomials $C_{n,m}$ and whose labeled edges signify the recurrence relation feeding one vertex into the ones below it to its right and left. A given vertex $C_{n,m}$ is the sum over the vertices $C_{i,j}$ at the tails of its incoming arrows with the operations labeling the incoming arrows applied to the $C_{i,j}$. \cite{nicebook}

The following diagram represents $\delta^n(g)$, $n\geq0$. The $n$th row down denotes $\delta^n(g)$, while the $m$th column to the right picks out the $kG$-coefficient of $x^mg$. Both $n$ and $m$ start from $0$. Multiplication by $\lambda$ corresponds to padding the partition by $1$: 
$$\lambda\cdot((p-2)^{a_{p-2}},...,2^{a_2},1^{a_1})=((p-2)^{a_{p-2}},...,2^{a_2},1^{a_1+1})$$
Southeast arrows in the diagram always correspond to applying $D+1$ so we have omitted their labels. Application of $D$ to a partition corresponds to differentiating that partition part by part using the Leibnitz rule as if a partition were a product of its parts and where $D(i)=i+1$. That is:  \begin{align*}D((p-2)^{a_{p-2}},...,2^{a_2},1^{a_1})=&a_{p-2}((p-2)^{a_{p-2}-1},...,2^{a_2},1^{a_1+1})+...+a_2((p-2)^{a_{p-2}},...,3^{a_3+1},2^{a_2-1},1^{a_1})\\&\qquad\qquad+a_1((p-2)^{a_{p-2}},...,2^{a_2+1},1^{a_1-1})
\end{align*}
(For simplicity, we have assumed $D^{p-1}(\lambda)=\lambda$ (the $t=0$ case) to illustrate the rule; otherwise an additional scalar term gets produced if $a_{p-2}$ is nonzero.)\\

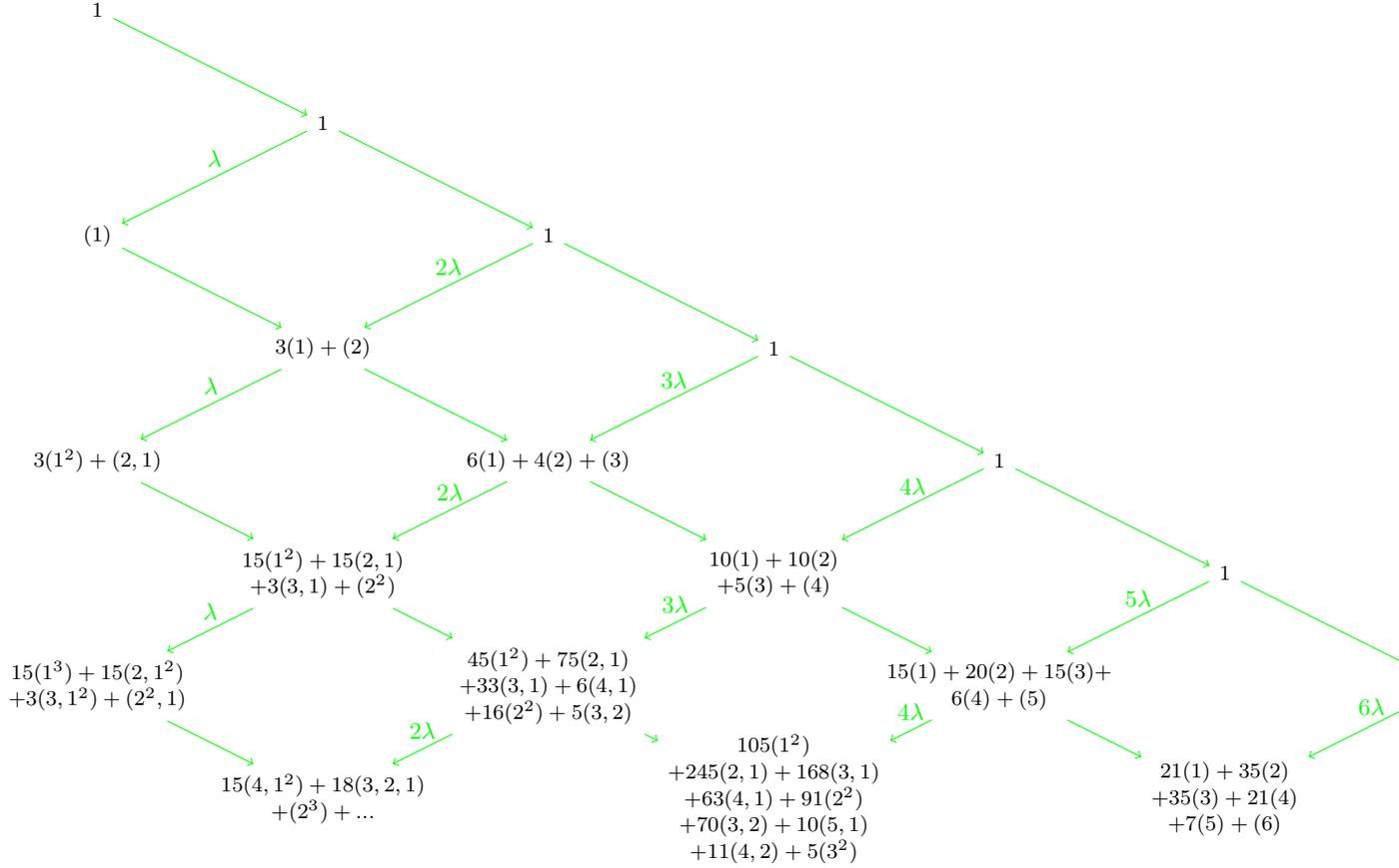
\begin{tikzpicture}[scale=3, font=\scriptsize,every text node part/.style={align=center}]
\node (A00) at (0,0) {$1$};
\node (A11) at (1, -.5) {$1$};
\node (A20) at (0,-1) {$(1)$};
\node (A22) at (2, -1) {$1$};
\node (A31) at (1,-1.5) {$3(1)+(2)$};
\node (A33) at (3,-1.5) {$1$};
\node (A40) at (0,-2) {$3(1^2)+(2,1)$};
\node (A42) at (2,-2) {$6(1)+4(2)+(3)$};
\node (A44) at (4,-2) {$1$};
\node (A51) at (1,-2.5) {$15(1^2)+15(2,1)$\\$+3(3,1)+(2^2)$};
\node (A53) at (3,-2.5) {$10(1)+10(2)$\\$+5(3)+(4)$};
\node (A55) at (5,-2.5) {$1$};
\node (A60) at (0,-3) {$15(1^3)+15(2,1^2)$\\$+3(3,1^2)+(2^2,1)$};
\node (A62) at (2,-3) {$45(1^2)+75(2,1)$\\$+33(3,1)+6(4,1)$\\$+16(2^2)+5(3,2)$};
\node (A64) at (4,-3) {$15(1)+20(2)+15(3)+$\\$6(4)+(5)$};
\node (A66) at (6,-3) {$1$};
\node (A71) at (1,-3.5) {$15(4,1^2)+18(3,2,1)$\\$+(2^3)+...$};
\node (A73) at (3,-3.5) {$105(1^2)$\\$+245(2,1)+168(3,1)$\\$+63(4,1)+91(2^2)$\\$+70(3,2)+10(5,1)$\\$+11(4,2)+5(3^2)$};
\node (A75) at (5,-3.5) {$21(1)+35(2)$\\$+35(3)+21(4)$\\$+7(5)+(6)$};
\node (A77) at (7,-3.5) {$1$};
\path[->,green,font=\footnotesize]
(A00) edge node {}(A11)
(A11) edge node [green][above]{$\lambda$} (A20)
(A11) edge node {}(A22)
(A20)edge node {}(A31)
(A22) edge node [green][above]{$2\lambda$}(A31)
(A22) edge node {}(A33)
(A31) edge node [green][above]{$\lambda$}(A40)
(A31) edge node {} (A42)
(A33) edge node [green][above]{$3\lambda$}(A42)
(A33) edge node {}(A44)
(A40) edge node {} (A51)
(A42) edge node [green][above]{$2\lambda$} (A51)
(A42) edge node {} (A53)
(A44) edge node [green][above]{$4\lambda$} (A53)
(A44) edge node {}(A55)
(A51) edge node [green][above]{$\lambda$} (A60)
(A51) edge node {} (A62)
(A53) edge node [green][above]{$3\lambda$} (A62)
(A53) edge node {} (A64)
(A55) edge node [green][above] {$5\lambda$} (A64)
(A55) edge node {}(A66)
(A60) edge node {}(A71)
(A62) edge node [green][above] {$2\lambda$} (A71)
(A62) edge node {}(A73)
(A64) edge node [green][above] {$4\lambda$} (A73)
(A64) edge node {} (A75)
(A66) edge node [green][above] {$6\lambda$} (A75)
(A66) edge node {}(A77);
\end{tikzpicture}\\

Compare this with the diagram for $\delta^n(\lambda)=\delta^{n+1}(x)$, starting from $\delta^0(\lambda)=\delta(x)=\lambda$. The $kG$-coefficient $\tilde{C}_{n,m}$ of $x^m$ in $\delta^n(\lambda)$ obeys the recurrence relation:
\begin{align*}
\tilde{C}_{0,0}&=\lambda,\qquad \tilde{C}_{n,m}=0\mbox{ if } n<m\mbox{ or } m<0\\
\tilde{C}_{n,m}&=D(\tilde{C}_{n-1,m-1})+(m+1)\lambda \tilde{C}_{n-1,m+1}\mbox{ for }n>0
\end{align*} Actually, we could tack on $\tilde{C}_{-1,1}=1$ since $\lambda=\tilde{C}_{0,0}=\delta(x)$. But it is not necessary. Southeast arrows always correspond to applying $D$. Notice that the diagram replicates the subdiagram of $\delta^n(g)$ where all partitions whose sum isn't the maximum possible for their node are erased, and the top diagonal is omitted. Thus the statement that $\delta^p(g)=\delta^{p-1}(x)g+x^pg=\delta^{p-2}(\lambda)g+x^pg$ is equivalent to the statement that the only partitions that survive in $\delta^p(g)$ mod $p$ are those whose sum is the maximum possible at those nodes, namely $p-1$. Here now is the diagram for $\delta^n(\lambda)=\delta^{n+1}(x)$ up to $n=6$:

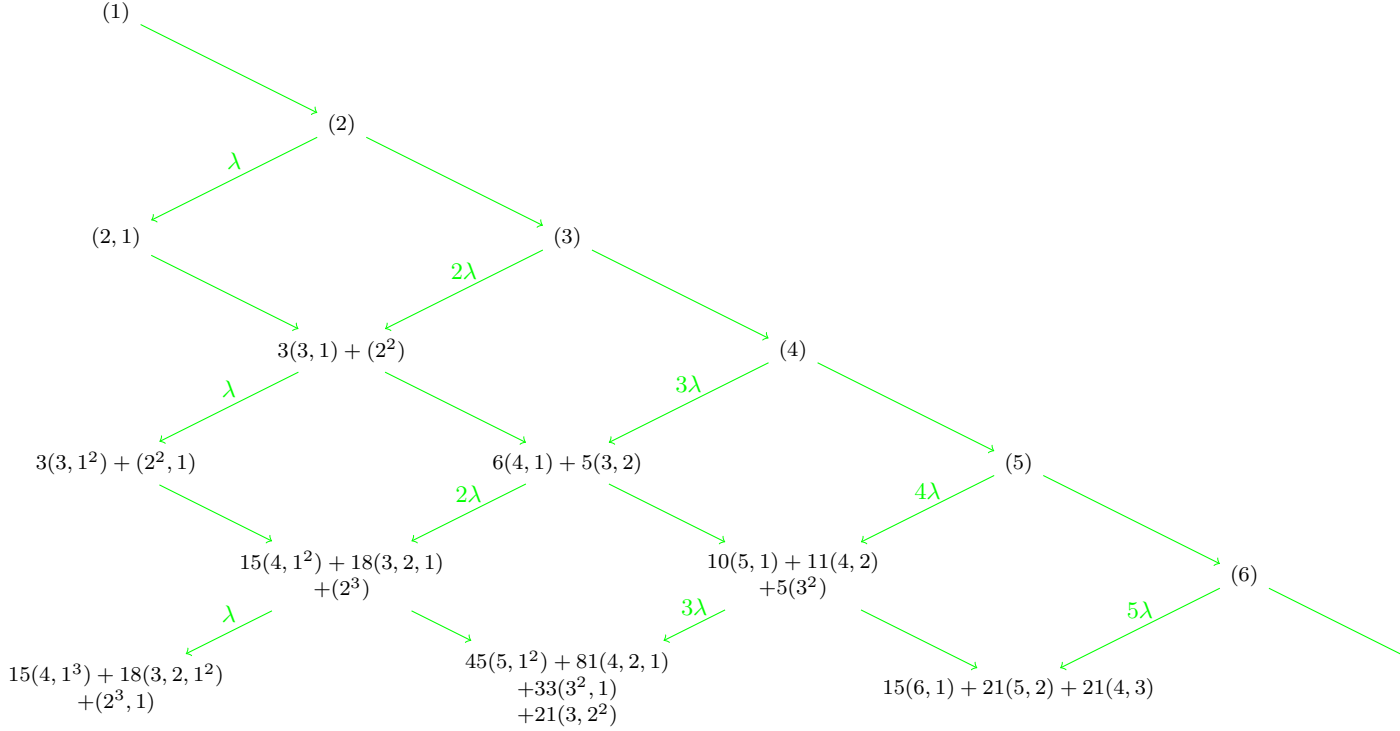
\begin{tikzpicture}[scale=3, font=\scriptsize,every text node part/.style={align=center}]
\node (A00) at (0,0) {$(1)$};
\node (A11) at (1, -.5) {$(2)$};
\node (A20) at (0,-1) {$(2,1)$};
\node (A22) at (2, -1) {$(3)$};
\node (A31) at (1,-1.5) {$3(3,1)+(2^2)$};
\node (A33) at (3,-1.5) {$(4)$};
\node (A40) at (0,-2) {$3(3,1^2)+(2^2,1)$};
\node (A42) at (2,-2) {$6(4,1)+5(3,2)$};
\node (A44) at (4,-2) {$(5)$};
\node (A51) at (1,-2.5) {$15(4,1^2)+18(3,2,1)$\\$+(2^3)$};
\node (A53) at (3,-2.5) {$10(5,1)+11(4,2)$\\$+5(3^2)$};
\node (A55) at (5,-2.5) {$(6)$};
\node (A60) at (0,-3) {$15(4,1^3)+18(3,2,1^2)$\\$+(2^3,1)$};
\node (A62) at (2,-3) {$45(5,1^2)+81(4,2,1)$\\$+33(3^2,1)$\\$+21(3,2^2)$};
\node (A64) at (4,-3) {$15(6,1)+21(5,2)+21(4,3)$};
\node (A66) at (6,-3) {$(7)$};
\path[->,green,font=\footnotesize]
(A00) edge node {}(A11)
(A11) edge node [above]{$\lambda$} (A20)
(A11) edge node {}(A22)
(A20)edge node {}(A31)
(A22) edge node [above]{$2\lambda$}(A31)
(A22) edge node {}(A33)
(A31) edge node [above]{$\lambda$}(A40)
(A31) edge node {} (A42)
(A33) edge node [above]{$3\lambda$}(A42)
(A33) edge node {}(A44)
(A40) edge node {} (A51)
(A42) edge node [above]{$2\lambda$} (A51)
(A42) edge node {} (A53)
(A44) edge node [above]{$4\lambda$} (A53)
(A44) edge node {}(A55)
(A51) edge node [above]{$\lambda$} (A60)
(A51) edge node {} (A62)
(A53) edge node [above]{$3\lambda$} (A62)
(A53) edge node {} (A64)
(A55) edge node [above] {$5\lambda$} (A64)
(A55) edge node {}(A66);
\end{tikzpicture}\\

Moreover, notice how the coefficients in the triangles for $\delta^n(g)$ and $\delta^n(\lambda)$ start off following the same pattern but then begin to diverge. In fact, between the two diagrams there's a simple and precise relationship: the triangle for $\delta^n(\lambda)$ is obtained from that for $\delta^n(g)$ by ``homogenizing" the partitions in each row. Namely, take a node in the $n$th row, $m$th column of the $\delta^n(g)$ triangle. Add a part to each partition at the node, so that every partition at that node has the sum of its parts equal to $n+1$. The result will be the $(n,m)$th node in the triangle for $\delta^n(\lambda)$ triangle. The following lemma proves this observation:
\lemma The triangle of $kG$-coefficients for $\delta^n(\lambda),\;n\geq0$, is the homogenization of the triangle of $kG$-coefficients for $\delta^n(g),\;n\geq0$.
\begin{proof}
We will find a recurrence relation for the ``homogenization" of $\delta^n(g)$ and show it's the same as the recurrence relation for $\delta^n(\lambda)$ with the same initial condition.

Write
$$C_{n,k}=\sum_{j\leq n}\sum_{\tworow{P\vdash j}{\# parts(P)=\frac{n-k}{2}}}c_PP$$
where $C_{n,k}$ is the coefficient of $x^mg$ in $\delta^n(g)$ and $c_P\in k$. For a partition $P=(P_1,...,P_{\frac{n-k}{2}})$ appearing in $C_{n,k}$ define $$\hat{P}:=(P_1,...,P_{\frac{n-k}{2}},n+1-j),$$ so that $\hat{P}$ is $P$ padded by $n+1-j$ where $j$ is the sum of the parts in $P$. Then define the homogenization of $C_{n,k}$ to be
$$\hat{C}_{n,k}=\sum_{j\leq n}\sum_{\tworow{P\vdash j}{\# parts(P)=\frac{n-k}{2}}}c_P\hat{P}$$
Then $$D(\hat{C}_{n,k-1})=\widehat{D(C_{n,k-1})}+\hat{C}_{n,k-1}$$ because differentiating all the parts of a partition $P$ except for the padded part $n+1-j$ ends up taking $D(P)$ which is a partition of $j+1$ and filling it out with $n+1-j=n+2-(j+1)$ to give a partition of $n+2$ -- this explains the first term on the right-hand-side. The second term on the right-hand-side comes from differentiating the part $n+1-j$ in each $P$ to get $n+2-j$. Moreover, it is clear from the definition that $$(k+1)\lambda\hat{C}_{n,k+1}=\widehat{(k+1)\lambda C_{n,k+1}}.$$
Since $C_{n,k}$ obeys the recurrence relation $$C_{n+1,k}=D(C_{n,k-1})+C_{n,k-1}+(k+1)\lambda C_{n,k+1}$$ it follows that $\hat{C}_{n,k}$ obeys the recurrence relation
\begin{align*}
\hat{C}_{n+1,k}&=\widehat{D(C_{n,k-1})}+\hat{C}_{n,k-1}+\widehat{(k+1)\lambda C_{n,k+1}}\\
&=D(\hat{C}_{n,k-1})+(k+1)\lambda\hat{C}_{n,k+1}\\
\end{align*}
which is the same recurrence relation satisfied by $\tilde{C}_{n,k},$ the coefficients in $\delta^n(\lambda)$. Moreover, ${C}_{0,0}$ is the empty partition, being a partition of $0$ by numbers greater than $0$, and thus its homogenization is a partition of $0+1=1$ given by padding the empty partition with $1$. Therefore $$\hat{C}_{0,0}=(1)=\tilde{C}_{0,0}$$ so that the two sequences share the same initial condition. It follows that $\hat{C}_{n,m}=\tilde{C}_{n,m}$ for all $n,m$, that is, $\delta^n(\lambda)$ is the homogenization of $\delta^n(g)$.
\end{proof}

\subsection{Generators and relations for the center of $\Hlam$}
Let $\Zlam$ denote the center of $\Hlam$. The main result of this section shows that $\Zlam$ is a finite module over its center.
\theorem \begin{enumerate}
\item ($t=0$ case). Let $\lambda\in kG$, $\lambda$ without constant term. Then the center of $\Hlam$ is generated by \begin{align*}
&x^2-2D^{p-2}(\lambda),\\
&x^p,\\
&y^p-y\frac{\delta^p(g)}{\delta(g)}
\end{align*}
and so $$\Zlam\cong \frac{k[A,B,C] }{(A^p-B^2)}$$
\item ($t=1$ case). Let $\lambda\in kG$, $\lambda$ with nonzero constant term. Then the center of $\Hlam$ is generated by \begin{align*}
&x^p,\\
&y^{p^2}-y^p\frac{\delta^{p^2}(g)}{\delta^p(g)}=\prod_{a=1}^p(Y-ax^p)\\
\end{align*}
where $Y:=y^p-y\frac{\delta^p(g)}{\delta(g)}$; and so: $$\Zlam\cong k[B,D]$$
\end{enumerate}
\begin{proof}
First, we recall a couple of basic facts about derivations, particularly in prime characteristic.
\lemma $ry^n=\sum_{i=0}^n{n\choose i}y^{n-i}\delta^i(r)$ \cite{Good}.
\lemma If $pR=0$ then $ry^p=y^pr+\delta^p(r)$.

Lemma 4.3 is easily proved by induction and holds over any field. Lemma 4.4 is an immediate consequence of the formula in Lemma 4.3. Lemma 4.4 implies 
\corollary $\delta^p$ is a derivation.
\begin{proof} For any $r,s\in R,$ $$\delta^p(rs)=[rs,y^p]=r[s,y^p]+[r,y^p]s=r\delta^p(s)+\delta^p(r)s$$\end{proof}

\remark From the perspective of $\Hlam$ as an Ore extension, saying that $B=x^p$ and $A=x^2-2(g\del_g)^{p-2}(\lambda)$ are in $\Zlam$ is the same as saying that they are killed by $\delta.$ 
\corollary $x^p\in\Zlam$.
\begin{proof}
$\delta(x^p)=px^{p-1}=0$.
\end{proof}
\lemma If $\lambda$ has no constant term then $x^2-2(\gdelg)^{p-2}(\lambda)\in\Zlam$. If $\lambda$ contains a constant term then $\Zlam\cap R=k[x^p]$.
\begin{proof}
If $\lambda$ has no constant term, $(\gdelg)^{p-1}(\lambda)=\lambda$, so that $\delta\left((g\del_g)^{p-2}(\lambda)\right)=x\lambda$, and therefore
$$\delta\left(x^2-2(g\del_g)^{p-2}(\lambda)\right)=2x\lambda-2x\lambda=0$$

If $\lambda$ has a constant term then $B=x^p\in \Zlam$ by Corollary 4.7. Suppose $f\in \R$ is of minimal degree in $x$ such that $f\in \Zlam$ but $f\notin k[x^p]$. Let $\phi:\Hlam\to \msf{H_1}$ be the map sending $g$ to $1$, $x$ to $x$ and $y$ to $y$. Then $\phi(\Zlam)\subset \msf{Z_1}=k[x^p, y^{p^2}-y^px^{p(p-1)}]$ so the powers of $x$ appearing in $f$ whose coefficients are not zero-divisors in $kG$ must be multiples of $p$.  On the other hand, the associated graded map $\gr:\Hlam\to\Hzero$ also maps $\Zlam$ into the center of $k[x,y]\rtimes G$: $\msf{Z_0}=k[x,y^p-yx^{p-1}]$. Then $\gr(f)$, which is equal to the top degree term of $f$, belongs to $\msf{Z_0}$ and thus the top degree term of $f$ has its coefficient in $k$. Therefore $f=x^{pn}+f'$ for some $n\in\mathbb{N}$ and some $f'\in R$, $f'\notin k[x^p]$, with $f'$ of strictly smaller degree than $f$ in the variable $x$. But $x^p$ is central in $\Hlam$, so $f'=f-x^{pn}$ is central, contradicting the minimality of $f$.
\end{proof}

The PBW theorem imposes restrictions on the generators of $\Zlam$. If $Z\in\Zlam$ then $$\gr(Z)\in\msf{Z}_0=k[x,y^p-yx^{p-1}]$$ and so a central element $Z$ of degree $p$ involving $y$ must be of the form $y^p-yx^{p-1}+l.o.t.'s$. Moreover, for $C=y^p-yf(x,g)$ to belong to the center would have to mean that $$f(x,g)=\frac{\delta^p(g)}{\delta(g)},$$ because $gy^p=y^pg+\delta^p(g)$ and $gy=yg+\delta(g)$. On the other hand, if $\lambda$ has a nonzero constant term then $\Zlam$ is like a deformation of $\msf{H}_1$, and the center of $\msf{H}_1=\msf{A}_1\rtimes G$ consists of the $G$-invariants of the center of the Weyl algebra: $$\msf{Z}_1=Z(\msf{A}_1)^G=k[x^p,y^p]^G=k[x^p,y^{p^2}-y^px^{p^2-p}].$$ Should $\Hlam$ lack a central element of the form $C$ as above when $t=1$, it will hopefully have one of the form $y^{p^2}-y^px^{p^2-p}+...$ so that when the deformation parameter $\lambda-1$ goes to $0$, so to speak, the resulting expression will be the second generator of $\msf{Z}_1$. Anyway, by associated graded considerations, $D\in\Zlam$ of degree $p^2$ and containing the term $y^{p^2}$ must have those leading terms. Moreover, for $D=y^{p^2}-y^ph(x,g)$ of degree $p^2$ to belong to $\Zlam$ would have to mean that $$h(x,g)=\frac{\delta^{p^2}(g)}{\delta^p(g)}$$ in order for $g$ and $D$ to commute. So we know what the central elements have to look like, if they exist.

\lemma $\delta(g)$ divides $\delta^p(g)$.
\begin{proof}We have $$\delta(g)=xg.$$Suppose by induction that $g$ divides $\delta^n(g)$, and write $\delta^n(g)=rg$ for some $r\in R$. Then $\delta^{n+1}(g)=\delta(r)g+r\delta(g)=\delta(r)g+rxg$ is also divisible by $g$. 

As for being divisible by $x$, $\delta(g)=xg$ and $\delta^2(g)=x^2g+\lambda g$. Suppose by induction that $$\delta^n(g)=x^n+x^{n-2}f_{n-2}(g)+...+x^{\epsilon}f_{\epsilon}(g)$$ for some polynomials $f_i(g)\in kG$, and where $\epsilon=0$ or $1$ depending on whether $n$ is even or odd; i.e. suppose the only powers of $x$ appearing in $\delta^n(g)$ are congruent to $n$ mod $2$. Applying $\delta$ to a monomial $x^kf_k(g)$ gives $$\delta(x^kf_k(g))=kx^{k-1}\delta(x)f_k(g)+x^{k+1}g\del_g(f_k(g))$$
and thus $$\delta^{n+1}(g)=x^{n+1}+x^{n-1}(f_{n-2}(g)+n\lambda f_n(g))+x^{n-3}(f_{n-4}+(n-2)\lambda f_{n-2}(g))+...$$
It follows that $x$ divides $\delta^n(g)$ if and only if $n$ is odd.

Consequently $\delta(g)=xg$ divides $\delta^p(g)$, and so $\frac{\delta^p(g)}{\delta(g)}$ is actually an element of $R$.
\end{proof}

Let $C_{n,k}$ denote the coefficient of $x^kg$ in $\delta^n(g)$, and let $\tilde{C}_{n,k}$ denote the coefficient of $x^k$ in $\delta^n(\lambda)$. Then $$\tilde{C}_{p-2,k}\equiv C_{p,k}\mbox{ mod }p$$ for $k=1,3,...,p-2$ since $$\delta^p(g)=\delta^{p-2}(\lambda)g+x^pg$$ by Proposition 3.4. Moreover, $$\tilde{C}_{p,k}=D(\lambda)C_{p,k}$$ for $k=1,3,...,p-2,p$ by Lemma 4.1 together with the observation that only the partitions in $C_{p,k}$ the sum of whose parts is maximal have nonzero coefficients mod $p$. 

Suppose $t=0.$ We claim that
$$\delta\left(\frac{\delta^p(g)}{\delta(g)}\right)=0.$$ $\delta(g)=xg$ and $\frac{\delta^p(g)}{xg}=C_{p,1}+C_{p,3}x^2+C_{p,5}x^4+...+C_{p,p-2}x^{p-3}+x^{p-1}$. Then 
\begin{align*}\delta\left(\frac{\delta^p(g)}{\delta(g)}\right)&=\left(D(C_{p,1})+2\lambda C_{p,3}\right)x+\left(D(C_{p,3})+4\lambda C_{p,5}\right)x^3+...+\left(D(C_{p,p-2})+(p-1)\lambda\right)x^{p-2}\\
&=\left(D(\tilde{C}_{p-2,1})+2\lambda \tilde{C}_{p-2,3}\right)x+\left(D(\tilde{C}_{p-2,3})+4\lambda \tilde{C}_{p-2,5}\right)x^3+...+\left(D(\tilde{C}_{p-2,p-2})+(p-1)\lambda\right)x^{p-2}
\end{align*}
What we must show, then, is that $$D(\tilde{C}_{p-2,k})=-(k+1)\lambda \tilde{C}_{p-2,k+2}$$ for $k=1,3,5,...,p-2$. This is true for $k=1$: the linear recursion for $\tilde{C}_{n,k}$ gives that $\tilde{C}_{p-1,0}=\lambda \tilde{C}_{p-2,1}$. Applying the result of Lemma 4.1 on the one hand, and the linear recursion relation on the other, \begin{align*}
D(\lambda)\tilde{C}_{p-2,1}=\tilde{C}_{p,1}&=D(\tilde{C}_{p-1,0})+2\lambda \tilde{C}_{p-1,2}\\&=D(\lambda)\tilde{C}_{p-2,1}+\lambda D(\tilde{C}_{p-2,1})+2\lambda\left(D(\tilde{C}_{p-2,1})+3\lambda \tilde{C}_{p-2,3} \right)\end{align*} and so $$\lambda D(\tilde{C}_{p-2,1})=-2\lambda^2 \tilde{C}_{p-2,3}$$
Now pretend that $\lambda$ is invertible so that we can cancel a $\lambda$ on both sides to get $$D(\tilde{C}_{p-2,1})=-2\lambda \tilde{C}_{p-2,3}$$
But now observe that the case that $\lambda$ is invertible proves it, in fact, for $\lambda$ a zero-divisor as well, since the coefficients of the partitions are the same no matter what $\lambda$ is.

From $D(\tilde{C}_{p-2,1})=-2\lambda \tilde{C}_{p-2,3}$ it follows that $\tilde{C}_{p-1,2}=\lambda \tilde{C}_{p-2,3}$; indeed, if $D(\tilde{C}_{p-2,k})=-(k+1)\tilde{C}_{p-2,k+2}$ then $\tilde{C}_{p-1,k+1}=\lambda \tilde{C}_{p-2,k+2}$ as an immediate consequence of the linear recursion which says $\tilde{C}_{p-1,k+1}=D(\tilde{C}_{p-2,k})+(k+2)\lambda \tilde{C}_{p-2,k+2}$. A diagram illustrates the situation:\\

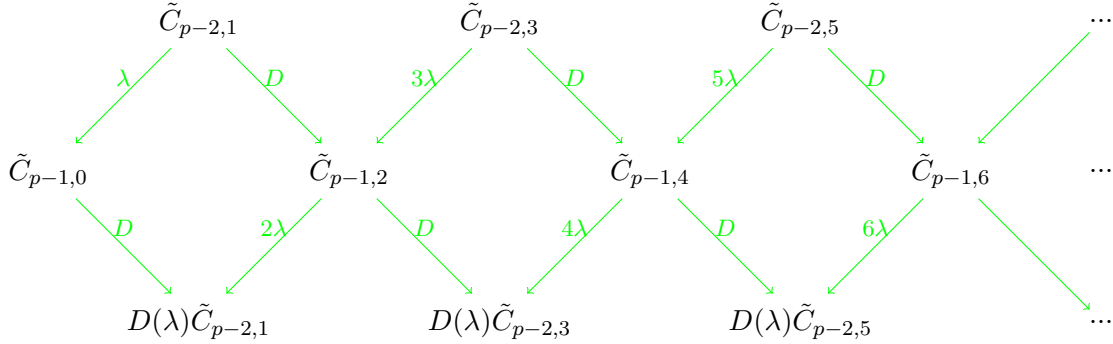
\begin{tikzpicture}
\node (C21) at (1,0) {$\tilde{C}_{p-2,1}$};
\node (C23) at (5,0) {$\tilde{C}_{p-2,3}$};
\node (C25) at (9,0) {$\tilde{C}_{p-2,5}$};
\node (C27) at (13,0) {...};
\node (C10) at (-1,-2) {$\tilde{C}_{p-1,0}$};
\node (C12) at (3,-2) {$\tilde{C}_{p-1,2}$};
\node (C14) at (7,-2) {$\tilde{C}_{p-1,4}$};
\node (C16) at (11,-2) {$\tilde{C}_{p-1,6}$};
\node (C01) at (1,-4) {$D(\lambda)\tilde{C}_{p-2,1}$};
\node (C03) at (5,-4) {$D(\lambda)\tilde{C}_{p-2,3}$};
\node (C05) at (9,-4) {$D(\lambda)\tilde{C}_{p-2,5}$};
\node (C07) at (13,-4) {...};
\node (C17) at (13,-2) {...};
\path[->,green,font=\footnotesize]
(C21) edge node [above] {$\lambda$} (C10)
(C21) edge node [above] {$D$} (C12)
(C10) edge node [above] {$D$} (C01)
(C12) edge node [above] {$2\lambda$} (C01)
(C23) edge node [above] {$3\lambda$} (C12)
(C23) edge node [above] {$D$} (C14)
(C12) edge node [above] {$D$} (C03)
(C14) edge node [above] {$4\lambda$} (C03)
(C14) edge node [above] {$D$} (C05)
(C25) edge node [above] {$5\lambda$} (C14)
(C25) edge node [above] {$D$} (C16)
(C16) edge node [above] {$6\lambda$} (C05)
(C16) edge node {} (C07)
(C27) edge node {} (C16);
\end{tikzpicture}

Induction takes care of $k=3,5,...,p-4$. If $D(\tilde{C}_{p-2,k})=-(k+1)\tilde{C}_{p-2,k+2}$ then $\tilde{C}_{p-1,k+1}=\lambda \tilde{C}_{p-2,k+2}$, from which it follows in turn, by the same calculation just done, that $D(\tilde{C}_{p-2,k+2})=-(k+3)\tilde{C}_{p-2,k+4}$. 

The only subtlety depending on the assumption $t=0$ or $t=1$ arises with the term $$D(\tilde{C}_{p-2,p-2})+(p-1)\lambda.$$ We have $\tilde{C}_{p-2,p-2}=D^{p-2}(\lambda)$ and thus $D(\tilde{C}_{p-2,p-2})=D^{p-1}(\lambda)$. The assumption that $\lambda$ does not contain a constant term, so $t=0$, implies that $D^{p-1}(\lambda)=\lambda$. In this case, then, \linebreak $D(\tilde{C}_{p-2,p-2})-\lambda=0$ and it follows that $\delta\left(\frac{\delta^p(g)}{\delta(g)}\right)=0$. This implies that $$\left[y,y^p-y\frac{\delta^p(g)}{\delta(g)}\right]=0.$$ By construction, $$\left[g,y^p-y\frac{\delta^p(g)}{\delta(g)}\right]=0$$ and we are done.

In the case $t=1$ where $\lambda$ contains a constant term, assumed after rescaling to be $1$, we then have $$D(\tilde{C}_{p-2,p-2})=D^{p-1}(\lambda)=\lambda-1$$
from which it follows that $$\delta\left(\frac{\delta^p(g)}{\delta(g)}\right)=-x^p.$$
Then $$\left[y,y^p-y\frac{\delta^p(g)}{\delta(g)}\right]=-yx^p.$$
Set $Y:=y^p-y\frac{\delta^p(g)}{\delta(g)}$.
\claim $$\prod_{a=1}^p\left(Y-ax^p\right)\in\Zlam$$
\begin{proof}
The element clearly commutes with $g$ and we must show it commutes with $y$. Let $\tilde{y}:=Y$, $\tilde{x}:=x^p$, $\tilde{g}:=y$. We know $$[\tilde{g},\tilde{y}]=-\tilde{x}\tilde{g},\qquad[\tilde{g},\tilde{x}]=[\tilde{x},\tilde{y}]=0$$ and we want to show that $\tilde{g}$ commutes with $\prod_{a=1}^p\left(\tilde{y}-a\tilde{x}\right)$.  But the relations between $\tilde{x}$, $\tilde{y},$ and $\tilde{g}$ are exactly those defining $k[\tilde{x},\tilde{y}]\rtimes \la \tilde{g}\ra$ where $\tilde{g}=g^{-1}$ instead of $g$ is taken as the generator. The center of $k[\tilde{x},\tilde{y}]\rtimes \la \tilde{g}\ra$, equal to the $\tilde{g}$-invariants of $k[\tilde{x},\tilde{y}]$, contains $$\prod_{a=1}^p\left(\tilde{y}-a\tilde{x}\right)=\prod_{a=1}^p\left(\tilde{y}-(a-1)\tilde{x}\right)=\tilde{g}\cdot\prod_{a=1}^p\left( \tilde{y}-a\tilde{x}\right).$$

\end{proof}
This concludes the proof of the theorem about the center of $\Hlam$.
\end{proof}

\subsection{The special case of $\lambda=g$ and the Andr\'{e} polynomials}
When $\lambda=g$ there are other nice explicit formulas for the central generator involving $y^p$:
$$C=y^p-y\frac{\delta^p(g)}{\delta(g)}=y^p-y\left(x^2-2g\right)^{\frac{p-1}{2}}=y^p-\sum_{k=0}^{\frac{p-1}{2}} \frac{(2k-1)!!}{k!}yx^{p-2k-1}g^{k}=\prod_{a=1}^p(y-ax)$$

Moreover, in the case $\lambda=g$, $\delta^n(g)$ has a significance in combinatorics: $$\delta^n(g)=A_n(g,x)$$ is the $n$th Andr\'{e} polynomial. The Andr\'{e} polynomials are a sequence of bivariate polynomials which, when evaluated at $(1,1)$, produce the Euler numbers $1,1,2,5,16,61,...$ (\href{http://oeis.org/A000111}{OEIS A000111}). The Euler numbers enumerate alternating permutations in the symmetric group $S_n$, that is, permutations $\sigma$ such that $\sigma(1)<\sigma(2)>\sigma(3)<\sigma(4)>...$. There is a monograph about the Andr\'{e} polynomials and some variations on them by Dominique Foata \cite{Foata}. The Andr\'{e} in question is D\'{e}sir\'{e} Andr\'{e}, a French mathematician of the late nineteenth century who discovered that the generating function for the number of alternating permutations in $S_n$ is $\sec x+\tan x$. No connection has previously been made between the Andr\'{e} polynomials and the algebra $\Hg$, or more broadly $\Hlam$.

Here are the first few Andr\'{e} polynomials:

\begin{align*}
&\del^0(g)=g\\
&\del^1(g)=xg\\
&\del^2(g)=x^2g+g^2\\
&\del^3(g)=x^3g+4xg^2\\
&\del^4(g)=x^4g+11x^2g^2+4g^3\\
&\del^5(g)=x^5g+26x^3g^2+34xg^3\\
&\del^6(g)=x^6g+57x^4g^2+180x^2g^3+34g^4\\
\end{align*}

We are interested in the integer sequence formed by the coefficients of the Andr\'{e} polynomials \href{http://oeis.org/A094503}{A094503}. After reindexing, we visualize the triangles of numbers as forming square arrays. The reason for visualizing the sequence in the square array is that the $n$th row keeps track of the coefficients of monomials of total degree $N-2n$ in $x$ and $y$ in the expansion of $(x+y)^N$ in $\Hg$. Precisely: the coefficient of $y^{N-2n-k}x^kg^{n}$ is then given by ${N\choose k}A_{n,k}$. Everything can be done over $\Z$, thinking of $g$ as infinite-order, generating a cyclic semigroup (isomorphic to $\N$) acting on $k\la x,y\ra$. The $n$'th knight's-move antidiagonal, starting from the top row and walking two spaces left and one down, two spaces left and one down,... until it reaches the leftmost column, reads off the coefficients of the $n$th Andr\'{e} polynomial as above.

 Each array is a directed graph with vertices at integer lattice points in the lower right quadrant of the plane and with edges given by the vectors $(1,0)$ and $(-1,-1)$ between the vertices, with multiplicities. Note that the square array arising from $\delta^m(x)$ coincides, for $\lambda=g$, with that for $\delta^m(g)$, but starts with the coefficient of $\delta^1(x)$ in the upper left corner (at the origin). The labeled arrows in the array show the linear recursion relation for the sequence as follows:  if vertex $v$ is at the head of $k$ arrows with tails at $w_1,...,w_k$, and those arrows are labeled with $a_1,...,a_k$ respectively, then $v=a_1w_1+...+a_kw_k$. The numbers $A_{n,k}$ in the sequence are placed at integer lattice points $(k,-n)$ in the lower right quadrant of the plane, starting with the initial value of $A_{0,0}=1$ at $(0,0)$. The numbers above the arrows can be interpreted as multiplicities of arrows, so that the value of an integer in the sequence placed at $(k,-n)$ equals the number of paths from $(0,0)$ to $(k,-n)$ in the graph \cite{nicebook}.\\

\begin{tikzpicture}[scale=2]
\node (T00) at (0,0) {$1$};
\node (T01) at (1,0) {$1$};
\node (T02) at (2,0) {$1$};
\node (T03) at (3,0) {$1$};
\node (T04) at (4,0) {$1$};
\node (T05) at (5,0) {$1$};
\node (T06) at (6,0) {$1$};
\node (T07) at (7,0) {$1$};
\node (T08) at (8,0) {$1...$};
\node (T10) at (0,-1) {$1$};
\node (T11) at (1,-1) {$4$};
\node (T12) at (2,-1) {$11$};
\node (T13) at (3,-1) {$26$};
\node (T14) at (4,-1) {$57$};
\node (T15) at (5,-1) {$120$};
\node (T16) at (6,-1) {$247$};
\node (T17) at (7,-1) {$502$};
\node (T18) at (8,-1){...};
\node (T20) at (0,-2) {$4$};
\node (T21) at (1,-2) {$34$};
\node (T22) at (2,-2) {$180$};
\node (T23) at (3,-2) {$768$};
\node (T24) at (4,-2) {$2904$};
\node (T25) at (5,-2) {$10194$};
\node (T26) at (6,-2) {...};
\node (T27) at (7,-2) {...};
\node (T28) at (8,-2) {...};
\node (T30) at (0,-3) {$34$};
\node (T31) at (1,-3) {$496$};
\node (T32) at (2,-3) {$4288$};
\node (T33) at (3,-3) {$28768$};
\node (T34) at (4,-3) {...};
\node (T35) at (5,-3) {...};
\node (T36) at (6,-3) {...};
\node (T37) at (7,-3) {...};
\node (T38) at (8,-3) {...};
\node (T40) at (0,-4) {496};
\node (T41) at (1,-4) {...};
\path[->,violet,font=\footnotesize,=angle 90]
(T00) edge node [above]{$1$} (T01)
(T01) edge node [above]{$1$} (T02)
(T02) edge node [above]{$1$} (T03)
(T03) edge node [above]{$1$} (T04)
(T04) edge node [above]{$1$} (T05)
(T05) edge node [above]{$1$} (T06)
(T06) edge node [above]{$1$} (T07)
(T07) edge node [above]{$1$} (T08)
(T01) edge node [left]{$1$} (T10)
(T02) edge node [left]{$2$} (T11)
(T03) edge node [left]{$3$} (T12)
(T04) edge node [left]{$4$} (T13)
(T05) edge node [left]{$5$} (T14)
(T06) edge node [left]{$6$} (T15)
(T07) edge node [left]{$7$} (T16)
(T08) edge node [left]{8} (T17)
(T10) edge node [above]{$2$} (T11)
(T11) edge node [above]{$2$} (T12)
(T12) edge node [above]{$2$} (T13)
(T13) edge node [above]{$2$} (T14)
(T14) edge node [above]{$2$} (T15)
(T15) edge node [above]{$2$} (T16)
(T16) edge node [above]{$2$} (T17)
(T17) edge node [above]{$2$} (T18)
(T11) edge node [left]{$1$} (T20)
(T12) edge node [left]{$2$} (T21)
(T13) edge node [left]{$3$} (T22)
(T14) edge node [left]{$4$} (T23)
(T15) edge node [left]{$5$} (T24)
(T16) edge node [left]{$6$} (T25)
(T17) edge node [left]{$7$} (T26)
(T20) edge node [above]{$3$} (T21)
(T21) edge node [above]{$3$} (T22)
(T22) edge node [above]{$3$} (T23)
(T23) edge node [above]{$3$} (T24)
(T24) edge node [above]{$3$} (T25)
(T25) edge node [above]{$3$} (T26)
(T21) edge node [left]{$1$} (T30)
(T22) edge node [left]{$2$} (T31)
(T23) edge node [left]{$3$} (T32)
(T24) edge node [left]{$4$} (T33)
(T25) edge node [left]{$5$} (T34)
(T30) edge node [above]{$4$} (T31)
(T31) edge node [above]{$4$} (T32)
(T32) edge node [above]{$4$} (T33)
(T33) edge node [above]{$4$} (T34)
(T31) edge node [left] {$1$} (T40);
\end{tikzpicture}

\proposition The integers $A_{n,k}$ satisfy the linear and quadratic recursion relations \begin{align}
A_{n,k}&=(n+1)A_{n,k-1}+(k+1)A_{n-1,k+1}\\
A_{n,k}&=A_{n,k-1}+\sum_{l=0}^{n-1}\sum_{m=0}^{k}{2n+k-1\choose2l+m}A_{l,m}A_{n-1-l,k-m}
\end{align}
with initial condition $$A_{0,0}=1,\qquad A_{n,k}=0\mbox{ if }n<0\mbox{ or if }k<0\mbox{ and }n>0.$$\\
While the proposition follows from facts about $\delta$ already established earlier, and available also in Foata's monograph \cite{Foata}, it can be proved directly by finding the coefficients of $(x+y)^{N+1}=(x+y)(x+y)^N=(x+y)^N(x+y)$ by identifying the coefficient of $y^{N-2k-n}x^ng^k$ in $(x+y)^N$ as ${N\choose 2k+n}A_{n,k}$ by induction and then computing which terms combine to give each coefficient in $(x+y)^{N+1}$.

\subsubsection{Integer sequences from $\delta^n$ when $\delta=g^a\del_x+xg\del_g$}
The polynomials $\delta^n(g)$ in $\Hlam$ for arbitrary $\lambda$ can be thought of as generalizations of the Andr\'{e} polynomials. When $\delta^n(g)$ is viewed as a polynomial in $x$ and $\gdelg$-powers of $\lambda$, its coefficients sum to the $n$th Euler number; when $\lambda=g$, this polynomial gives back the $n$th Andr\'{e} polynomial, since $(\gdelg)^r(g)=g$. These polynomials are already interesting when $\lambda=g^a$ for $a>1$, $a\in\bb{Z}$. This section discusses positive integer sequences related to the Andr\'{e} numbers which result from replacing the differential operator $\delta=g\del_x+xg\del_g$ with $\delta=g^a\del_x+xg\del_g$.

For the first example, consider $\delta=g^2\del_x+xg\del_g$. This should be something special, since if we thought of $g$ as being the shadow of a generator of a polynomial ring where $g$, like $x$ and $y$, had degree $1$, then $\delta$ would be homogeneous and the resulting algebra would be graded.

Here are the first few powers of $\delta$ on $g$ and on $x$. Rather than computing them as polynomials in $g^2,\;\delta(g^2),...$ etc, we compute them as polynomials in $g$ and in $x$, and what we get is:

\begin{align*}
&\del^0(g)=g\\
&\del^1(g)=xg\\
&\del^2(g)=x^2g+g^3\\
&\del^3(g)=x^3g+5xg^3\\
&\del^4(g)=x^4g+18x^2g^3+5g^5\\
&\del^5(g)=x^5g+58x^3g^3+61xg^5\\
&\del^6(g)=x^6g+179x^4g^3+479x^3g^5+61g^7\\
&\del^7(g)=x^7g+543x^5g^3+3111x^3g^5+1385xg^7\\
&\del^8(g)=x^8g+1636x^6g^3+18270x^4g^5+19028x^3g^7+1385g^9\\
&\\
&\\
&\del^0(x)=x\\
&\del^1(x)=g^2\\
&\del^2(x)=2xg^2\\
&\del^3(x)=4x^2g^2+2g^4\\
&\del^4(x)=8x^3g^2+16xg^4\\
&\del^5(x)=16x^4g^2+88x^2g^4+16g^6\\
&\del^6(x)=32x^5g^2+416x^2g^4+272xg^6\\
&\del^7(x)=64x^6g^2+1824x^4g^4+2880x^3g^6+272g^8\\
&\del^8(x)=128x^7g^2+7680x^5g^4+24576x^2g^6+7936xg^8\\
\end{align*}

Indeed, something special happens with the coefficients...:

\begin{align*}
&\del^0(g)(1,1)={\color{green}1}={\color{blue}0!}\\
&\del^1(g)(1,1)={\color{green}1}={\color{blue}1!}\\
&\del^2(g)(1,1)={\color{green}1+1}=2={\color{blue}2!}\\
&\del^3(g)(1,1)={\color{green}1+5}=6={\color{blue}3!}\\
&\del^4(g)(1,1)={\color{green}1+18+5}=34={\color{blue}4!}\\
&\del^5(g)(1,1)={\color{green}1+58+61}=120={\color{blue}5!}\\
&\del^6(g)(1,1)={\color{green}1+179+479+61}=730={\color{blue}6!}\\
&\del^7(g)(1,1)={\color{green}1+543+3111+1385}=5040={\color{blue}7!}\\
&\del^8(g)(1,1)={\color{green}1+1636+18270+19028+1385}=40320={\color{blue}8!}\\
&\\
&\\
&\del^0(x)(1,1)={\color{green}1}={\color{blue}0!}\\
&\del^1(x)(1,1)={\color{green}1}={\color{blue}1!}\\
&\del^2(x)(1,1)={\color{green}2}={\color{blue}2!}\\
&\del^3(x)(1,1)={\color{green}4+2}=6={\color{blue}3!}\\
&\del^4(x)(1,1)={\color{green}8+16}=24={\color{blue}4!}\\
&\del^5(x)(1,1)={\color{green}16+88+16}=120={\color{blue}5!}\\
&\del^6(x)(1,1)={\color{green}32+416+272}=720={\color{blue}6!}\\
&\del^7(x)(1,1)={\color{green}64+1824+2880+272}=5040={\color{blue}7!}\\
&\del^8(x)(1,1)={\color{green}128+7680+24576+7936}=40320={\color{blue}8!}\\
\end{align*}

In both examples, the coefficients of the polynomials $\delta^m(g)$ and $\delta^m(x)$ (the numbers in green) form a triangle of positive integers whose rows partition the factorials (in blue). It is very easy to show that the $n$th row sums do indeed add up to $n!$, using quadratic recurrence relations for $\delta^n(g)$ and $\delta^n(x)$. For example: 
$$\delta^n(g)=\sum_{i=0}^{n-1}{n-1\choose i}\delta^i(x)\delta^{n-1-i}(g)$$
and applying induction and evaluating both sides at $(1,1)$, we get 
$$\delta^n(g)(1,1)=\sum_{i=0}^{n-1}{n-1\choose i} i!(n-1-i)!=\sum_{i=0}^{n-1}(n-1)!=n!$$
Likewise the same argument applies to the evaluation of $\delta^n(x)$ at $(1,1)$, but via the quadratic recurrence relation $$\delta^n(x)=2\sum_{i=1}^{n-1}{n-2\choose i-1}\delta^i(x)\delta^{n-1-i}(x).$$ This proves:
\proposition When $\delta=g^2\del_x+xg\del_g$, the sum of the integer coefficients of $\delta^n(x)$ and the sum of the integer coefficients of $\delta^n(g)$ coincide and equal $n!$.

Compare this result to a similar result in \cite{Benk2}, where a different partitioning of factorials is also obtained from certain coefficients appearing in powers of $\delta$ in the setting of an Ore extension of $k[x]$. It's an interesting coincidence.

Here are the first few rows of the array of coefficients arising from $\delta^m(g)$ in its square array form:\\

\begin{tikzpicture}[scale=2]
\node (T00) at (0,0) {$1$};
\node (T01) at (1,0) {$1$};
\node (T02) at (2,0) {$1$};
\node (T03) at (3,0) {$1$};
\node (T04) at (4,0) {$1$};
\node (T05) at (5,0) {$1$};
\node (T06) at (6,0) {$1$};
\node (T07) at (7,0) {$1$};
\node (T08) at (8,0) {$1...$};
\node (T10) at (0,-1) {$1$};
\node (T11) at (1,-1) {$5$};
\node (T12) at (2,-1) {$18$};
\node (T13) at (3,-1) {$58$};
\node (T14) at (4,-1) {$179$};
\node (T15) at (5,-1) {$543$};
\node (T16) at (6,-1) {$1636$};
\node (T17) at (7,-1) {$4916$};
\node (T18) at (8,-1){...};
\node (T20) at (0,-2) {$5$};
\node (T21) at (1,-2) {$61$};
\node (T22) at (2,-2) {$479$};
\node (T23) at (3,-2) {$3111$};
\node (T24) at (4,-2) {$18270$};
\node (T25) at (5,-2) {$101166$};
\node (T26) at (6,-2) {...};
\node (T27) at (7,-2) {...};
\node (T28) at (8,-2) {...};
\node (T30) at (0,-3) {$61$};
\node (T31) at (1,-3) {$1385$};
\node (T32) at (2,-3) {$19028$};
\node (T33) at (3,-3) {$206276$};
\node (T34) at (4,-3) {...};
\node (T35) at (5,-3) {...};
\node (T36) at (6,-3) {...};
\node (T37) at (7,-3) {...};
\node (T38) at (8,-3) {...};
\path[->,magenta,font=\footnotesize,=angle 90]
(T00) edge node [above]{$1$} (T01)
(T01) edge node [above]{$1$} (T02)
(T02) edge node [above]{$1$} (T03)
(T03) edge node [above]{$1$} (T04)
(T04) edge node [above]{$1$} (T05)
(T05) edge node [above]{$1$} (T06)
(T06) edge node [above]{$1$} (T07)
(T07) edge node [above]{$1$} (T08)
(T01) edge node [left]{$1$} (T10)
(T02) edge node [left]{$2$} (T11)
(T03) edge node [left]{$3$} (T12)
(T04) edge node [left]{$4$} (T13)
(T05) edge node [left]{$5$} (T14)
(T06) edge node [left]{$6$} (T15)
(T07) edge node [left]{$7$} (T16)
(T08) edge node [left]{8} (T17)
(T10) edge node [above]{$3$} (T11)
(T11) edge node [above]{$3$} (T12)
(T12) edge node [above]{$3$} (T13)
(T13) edge node [above]{$3$} (T14)
(T14) edge node [above]{$3$} (T15)
(T15) edge node [above]{$3$} (T16)
(T16) edge node [above]{$3$} (T17)
(T17) edge node [above]{$3$} (T18)
(T11) edge node [left]{$1$} (T20)
(T12) edge node [left]{$2$} (T21)
(T13) edge node [left]{$3$} (T22)
(T14) edge node [left]{$4$} (T23)
(T15) edge node [left]{$5$} (T24)
(T16) edge node [left]{$6$} (T25)
(T17) edge node [left]{$7$} (T26)
(T20) edge node [above]{$5$} (T21)
(T21) edge node [above]{$5$} (T22)
(T22) edge node [above]{$5$} (T23)
(T23) edge node [above]{$5$} (T24)
(T24) edge node [above]{$5$} (T25)
(T25) edge node [above]{$5$} (T26)
(T21) edge node [left]{$1$} (T30)
(T22) edge node [left]{$2$} (T31)
(T23) edge node [left]{$3$} (T32)
(T24) edge node [left]{$4$} (T33)
(T25) edge node [left]{$5$} (T34)
(T30) edge node [above]{$7$} (T31)
(T31) edge node [above]{$7$} (T32)
(T32) edge node [above]{$7$} (T33)
(T33) edge node [above]{$7$} (T34);
\end{tikzpicture}

Here is the array of coefficients arising from $\delta^m(x)$ starting from the coefficient of $\delta(x)$ in the upper right hand corner:\\

\begin{tikzpicture}[scale=2]
\node (T00) at (0,0) {$1$};
\node (T01) at (1,0) {$2$};
\node (T02) at (2,0) {$4$};
\node (T03) at (3,0) {$8$};
\node (T04) at (4,0) {$16$};
\node (T05) at (5,0) {$32$};
\node (T06) at (6,0) {$64$};
\node (T07) at (7,0) {$128$};
\node (T08) at (8,0) {$256...$};
\node (T10) at (0,-1) {$2$};
\node (T11) at (1,-1) {$16$};
\node (T12) at (2,-1) {$88$};
\node (T13) at (3,-1) {$416$};
\node (T14) at (4,-1) {$1824$};
\node (T15) at (5,-1) {$7680$};
\node (T16) at (6,-1) {$31616$};
\node (T17) at (7,-1) {$128512$};
\node (T18) at (8,-1){...};
\node (T20) at (0,-2) {$16$};
\node (T21) at (1,-2) {$272$};
\node (T22) at (2,-2) {$2880$};
\node (T23) at (3,-2) {$24576$};
\node (T24) at (4,-2) {$185856$};
\node (T25) at (5,-2) {$1304832$};
\node (T26) at (6,-2) {...};
\node (T27) at (7,-2) {...};
\node (T28) at (8,-2) {...};
\node (T30) at (0,-3) {$272$};
\node (T31) at (1,-3) {$7936$};
\node (T32) at (2,-3) {$137216$};
\node (T33) at (3,-3) {$1841152$};
\node (T34) at (4,-3) {...};
\node (T35) at (5,-3) {...};
\node (T36) at (6,-3) {...};
\node (T37) at (7,-3) {...};
\node (T38) at (8,-3) {...};
\path[->,cyan,font=\footnotesize,=angle 90]
(T00) edge node [above]{$2$} (T01)
(T01) edge node [above]{$2$} (T02)
(T02) edge node [above]{$2$} (T03)
(T03) edge node [above]{$2$} (T04)
(T04) edge node [above]{$2$} (T05)
(T05) edge node [above]{$2$} (T06)
(T06) edge node [above]{$2$} (T07)
(T07) edge node [above]{$2$} (T08)
(T01) edge node [left]{$1$} (T10)
(T02) edge node [left]{$2$} (T11)
(T03) edge node [left]{$3$} (T12)
(T04) edge node [left]{$4$} (T13)
(T05) edge node [left]{$5$} (T14)
(T06) edge node [left]{$6$} (T15)
(T07) edge node [left]{$7$} (T16)
(T08) edge node [left]{8} (T17)
(T10) edge node [above]{$4$} (T11)
(T11) edge node [above]{$4$} (T12)
(T12) edge node [above]{$4$} (T13)
(T13) edge node [above]{$4$} (T14)
(T14) edge node [above]{$4$} (T15)
(T15) edge node [above]{$4$} (T16)
(T16) edge node [above]{$4$} (T17)
(T17) edge node [above]{$4$} (T18)
(T11) edge node [left]{$1$} (T20)
(T12) edge node [left]{$2$} (T21)
(T13) edge node [left]{$3$} (T22)
(T14) edge node [left]{$4$} (T23)
(T15) edge node [left]{$5$} (T24)
(T16) edge node [left]{$6$} (T25)
(T17) edge node [left]{$7$} (T26)
(T20) edge node [above]{$6$} (T21)
(T21) edge node [above]{$6$} (T22)
(T22) edge node [above]{$6$} (T23)
(T23) edge node [above]{$6$} (T24)
(T24) edge node [above]{$6$} (T25)
(T25) edge node [above]{$6$} (T26)
(T21) edge node [left]{$1$} (T30)
(T22) edge node [left]{$2$} (T31)
(T23) edge node [left]{$3$} (T32)
(T24) edge node [left]{$4$} (T33)
(T25) edge node [left]{$5$} (T34)
(T30) edge node [above]{$8$} (T31)
(T31) edge node [above]{$8$} (T32)
(T32) edge node [above]{$8$} (T33)
(T33) edge node [above]{$8$} (T34);
\end{tikzpicture}

Now replace $2$ with any $a\in\bb{N}$ and consider $\delta=g^a\del_x+xg\del_g$. 
The linear recursion relation for the array arising from $\delta^n(g)$ as in the magenta sequence above is $$T_{n,k}=(an+1)T_{n,k-1}+(k+1)T_{n-1,k+1}$$ while the linear recursion relation for the array arising from $\delta^n(x)$ as in the cyan sequence above is $$U_{n,k}=a(n+1)U_{n,k-1}+(k+1)U_{n-1,k+1}.$$ The sequence $T_{n,k}$ obeys a quadratic recurrence relation: recall that $\delta^m(g)$ obeys a quadratic recurrence because $\delta(g)=xg$: 
$$\delta^m(g)=\sum_{i=0}^{m-1}{m-1\choose i}\delta^i(x)\delta^{m-1-i}(g).$$ After a change of variables, and looking at each monomial in $\delta^m(g)$, it follows that
$$T_{n,k}=T_{n,k-1}+\sum_{i=0}^{n-1}\sum_{j=0}^k{2n+k-1\choose2i+j}T_{i,j}U_{n-1-i,k-j},$$ and thus the sequences $T_{n,k}$ and $U_{n,k}$ are interlinked. On the other hand, the sequence $U_{n,k}$ obtained from $\delta^m(x)$ also satisfies a quadratic recurrence relation but one that does not involve $T_{n,k}$:
$$U_{n,k}=aU_{n,k-1}+\sum_{i=0}^{n-1}\sum_{j=0}^k{2n+k-1\choose2i+j}aU_{i,j}U_{n-1-i,k-j}.$$  We recognize this equation, absent the $a$'s, from the quadratic recurrence equation defining the Andr\'{e} numbers $\{A_{n,k}\}$. As the starting values $U_{0,0}$ and $T_{0,0}$ for $\{U_{n,k}\}$ and $\{A_{n,k}\}$ are both $1$, it follows by induction that \proposition $$U_{n,k}=a^{n+k}A_{n,k}.$$ So the sequence $U_{n,k}$ arising from $\delta^m(x)$ simply looks like:\\
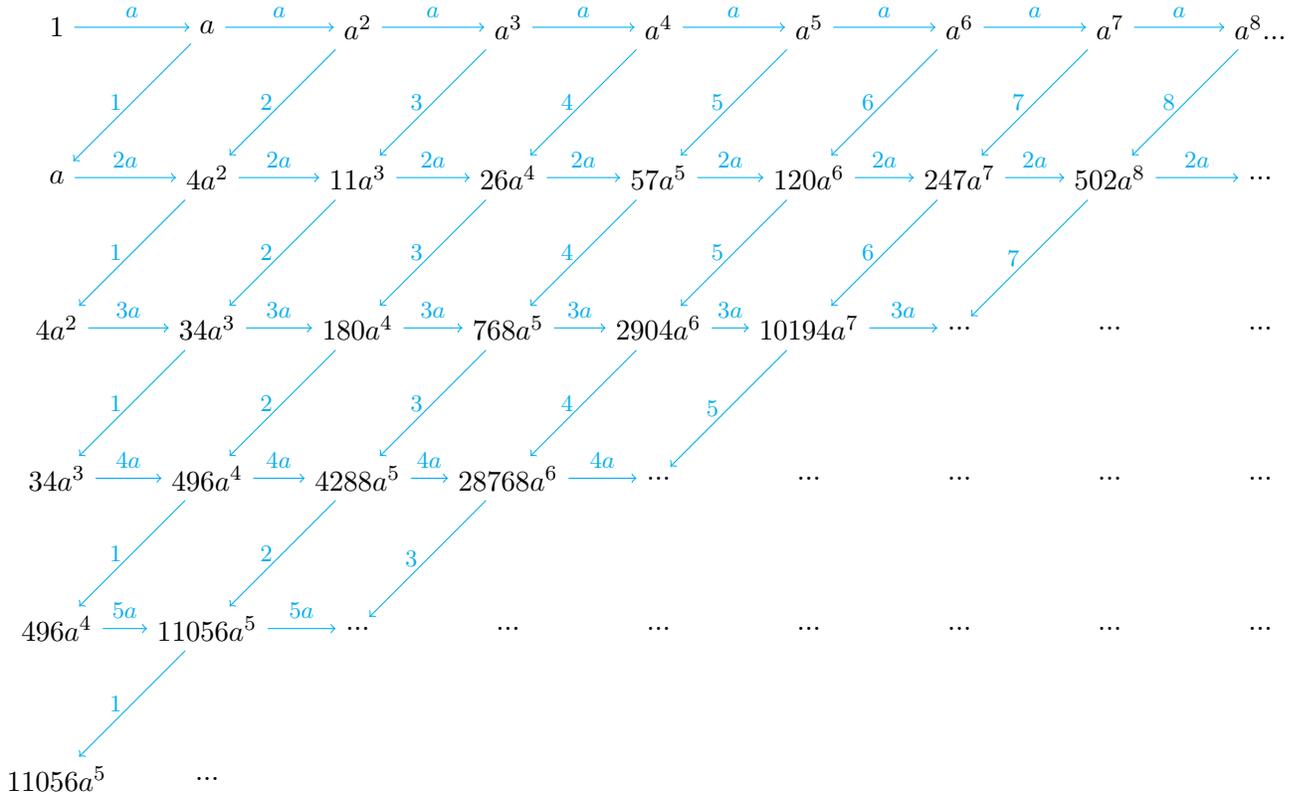
\begin{tikzpicture}[scale=2]
\node (T00) at (0,0) {$1$};
\node (T01) at (1,0) {$a$};
\node (T02) at (2,0) {$a^2$};
\node (T03) at (3,0) {$a^3$};
\node (T04) at (4,0) {$a^4$};
\node (T05) at (5,0) {$a^5$};
\node (T06) at (6,0) {$a^6$};
\node (T07) at (7,0) {$a^7$};
\node (T08) at (8,0) {$a^8...$};
\node (T10) at (0,-1) {$a$};
\node (T11) at (1,-1) {$4a^2$};
\node (T12) at (2,-1) {$11a^3$};
\node (T13) at (3,-1) {$26a^4$};
\node (T14) at (4,-1) {$57a^5$};
\node (T15) at (5,-1) {$120a^6$};
\node (T16) at (6,-1) {$247a^7$};
\node (T17) at (7,-1) {$502a^8$};
\node (T18) at (8,-1){...};
\node (T20) at (0,-2) {$4a^2$};
\node (T21) at (1,-2) {$34a^3$};
\node (T22) at (2,-2) {$180a^4$};
\node (T23) at (3,-2) {$768a^5$};
\node (T24) at (4,-2) {$2904a^6$};
\node (T25) at (5,-2) {$10194a^7$};
\node (T26) at (6,-2) {...};
\node (T27) at (7,-2) {...};
\node (T28) at (8,-2) {...};
\node (T30) at (0,-3) {$34a^3$};
\node (T31) at (1,-3) {$496a^4$};
\node (T32) at (2,-3) {$4288a^5$};
\node (T33) at (3,-3) {$28768a^6$};
\node (T34) at (4,-3) {...};
\node (T35) at (5,-3) {...};
\node (T36) at (6,-3) {...};
\node (T37) at (7,-3) {...};
\node (T38) at (8,-3) {...};
\node (T40) at (0,-4) {$496a^4$};
\node (T41) at (1,-4) {$11056a^5$};
\node (T42) at (2,-4) {...};
\node (T43) at (3,-4) {...};
\node (T44) at (4,-4) {...};
\node (T45) at (5,-4) {...};
\node (T46) at (6,-4) {...};
\node (T47) at (7,-4) {...};
\node (T48) at (8,-4) {...};
\node (T50) at (0,-5) {$11056a^5$};
\node (T51) at (1,-5) {...};
\path[->,cyan,font=\footnotesize,=angle 90]
(T00) edge node [above]{$a$} (T01)
(T01) edge node [above]{$a$} (T02)
(T02) edge node [above]{$a$} (T03)
(T03) edge node [above]{$a$} (T04)
(T04) edge node [above]{$a$} (T05)
(T05) edge node [above]{$a$} (T06)
(T06) edge node [above]{$a$} (T07)
(T07) edge node [above]{$a$} (T08)
(T01) edge node [left]{$1$} (T10)
(T02) edge node [left]{$2$} (T11)
(T03) edge node [left]{$3$} (T12)
(T04) edge node [left]{$4$} (T13)
(T05) edge node [left]{$5$} (T14)
(T06) edge node [left]{$6$} (T15)
(T07) edge node [left]{$7$} (T16)
(T08) edge node [left]{8} (T17)
(T10) edge node [above]{$2a$} (T11)
(T11) edge node [above]{$2a$} (T12)
(T12) edge node [above]{$2a$} (T13)
(T13) edge node [above]{$2a$} (T14)
(T14) edge node [above]{$2a$} (T15)
(T15) edge node [above]{$2a$} (T16)
(T16) edge node [above]{$2a$} (T17)
(T17) edge node [above]{$2a$} (T18)
(T11) edge node [left]{$1$} (T20)
(T12) edge node [left]{$2$} (T21)
(T13) edge node [left]{$3$} (T22)
(T14) edge node [left]{$4$} (T23)
(T15) edge node [left]{$5$} (T24)
(T16) edge node [left]{$6$} (T25)
(T17) edge node [left]{$7$} (T26)
(T20) edge node [above]{$3a$} (T21)
(T21) edge node [above]{$3a$} (T22)
(T22) edge node [above]{$3a$} (T23)
(T23) edge node [above]{$3a$} (T24)
(T24) edge node [above]{$3a$} (T25)
(T25) edge node [above]{$3a$} (T26)
(T21) edge node [left]{$1$} (T30)
(T22) edge node [left]{$2$} (T31)
(T23) edge node [left]{$3$} (T32)
(T24) edge node [left]{$4$} (T33)
(T25) edge node [left]{$5$} (T34)
(T30) edge node [above]{$4a$} (T31)
(T31) edge node [above]{$4a$} (T32)
(T32) edge node [above]{$4a$} (T33)
(T33) edge node [above]{$4a$} (T34)
(T31) edge node [left]{$1$} (T40)
(T32) edge node [left]{$2$} (T41)
(T33) edge node [left]{$3$} (T42)
(T40) edge node [above]{$5a$} (T41)
(T41) edge node [above]{$5a$} (T42)
(T41) edge node [left]{$1$} (T50);
\end{tikzpicture}

\corollary Returning to the example of $xy=yx+g^2$ and the triangular array for $\delta^n(x)$, it follows that weighted row sums of Andr\'{e} numbers partition the factorials: 
$$n!=2^nA_{0,n}+2^{n-2}A_{1,n-2}+2^{n-4}A_{2,n-4}+...+2^{n-2\lfloor \frac{n}{2}\rfloor}A_{\lfloor \frac{n}{2}\rfloor,n-2\lfloor \frac{n}{2}\rfloor}$$

There's got to be a combinatorial explanation for this formula: $n!$ is the order of the symmetric group $S_n$, and $A_{k,n-2k}$ counts some type of permutation (``Andr\'{e} permutation") in $S_n$ \cite{Foata}.

\remark In the case $a=3$, the leftmost column of the square array corresponding to the sequence $\delta^n(g)$ is OEIS sequence \href{http://oeis.org/A126151}{A126151}: $1,6,96,2976,151416,...$.

\remark We also note that in the case $a=4$, the row sums of the triangle corresponding to the sequence $\delta^n(x)$ (knight's-move diagonal sums in the square version), i.e. the sequence $\delta^n(x)(1,1)$, count the number of minimax trees on $n$ vertices (\href{http://oeis.org/A080795}{A080795}). Minimax trees, a family of recursively generated binary trees, were defined and studied by Foata and Han and are closely related to Andr\'{e} numbers \cite{FH}.

\section{Representations of $\Hlam$}
\subsection{Simple $\Hlam$-modules in the $t=0$ case}
With an explicit description of the center, it is not difficult to find the simple modules over $\Hg$, and in general, $\Hlam$ when $t=0$. We follow the strategy of \cite{GS} and identify simple modules over $\Hlam$ on which the center $\Zlam$ acts by a fixed character $(\alpha,\beta)\in \Spec \Zlam$ with simple modules over a finite-dimensional algebra $\Hab$ obtained as the quotient of $\Hlam$ by the two-sided ideal where the central generators are set equal to corresponding scalars. This process is called central reduction, and the geometry of the center gives important information about simple modules and their dimensions. 
\remark We know a priori that all simple modules over $\Hlam$ are finite-dimensional because $\Hlam$ is finite over its center.

Points in the affine space $\Spec \Zlam$ correspond to central characters. In $\Spec \Zlam$, there are three loci of particular representation-theoretic interest, which Gordon and Smith explain as follows in their study of representations of wreath product symplectic reflection algebras \cite{GS}: 
\begin{itemize}
\item The smooth locus or non-singular locus: where the Jacobian of the relations among the generators defining $\Spec \Zlam$ is non-degenerate
\item The Azumaya locus: a point is in the Azumaya locus if the simple module on which $\Hlam$ acts with the corresponding central character is of the maximal possible dimension for a simple module over $\Hlam$.
\item The locus of points where each point corresponds to a unique simple module up to isomorphism, i.e. where $\Hlam$ only has one simple module on which it acts with that central character.
\end{itemize}
In characteristic $0$, Gordon and Smith proved that the smooth locus, the Azumaya locus, and the one-to-one \{simple modules\}-to-\{points in Spec of the center\} locus are all the same for Spec of the center of a wreath product symplectic reflection algebra (with $t=0$).

In characteristic $p$, the theorem of Gordon and Smith need not hold. We will see that there is a unique simple for every point in $\Spec \Zlam$ and thus the one-to-one locus always strictly contains the smooth locus, which is $\Spec \Zlam$ minus a line. The Azumaya locus includes the singular line (i.e. is the whole space $\Spec \Zlam$) when $\lambda$ is invertible, so that when $\lambda$ is invertible the last two out of the three loci coincide. The Azumaya locus coincides exactly with the smooth locus when $\lambda$ is a $0$-divisor, and so in this case the first two of the three loci coincide.

\subsubsection{Simple representations of $\Hg$.}
First we'll describe central reduction for $\Hg$ and the simple modules in this special case.  $\Zgee$ is generated by $X=x^2-2g$, $Y=x^p$, and $Z=y^p-yX^{\frac{p-1}{2}}$ with the relation $X^p+2=Y^2$. Take $(\alpha, \beta)\in\mathbb{\bar{F}}_p^2$. Then $\mathfrak{m}_x:=(X-\sqrt[p]{\alpha}+2,Y-\sqrt{\alpha},Z-\beta)$ is a maximal two-sided ideal of $Z(\Hg)$. Let $\Hab:=\Hg/\mathfrak{m}_x\Hg$. Thus $\Hab$ is the algebra $\Hg$ subject to the additional relations: $$x^p=\sqrt{\alpha},$$ $$x^2-2g=\sqrt[p]{\alpha}-2,$$ $$y^p-y(\sqrt[p]{\alpha}-2)^{\frac{p-1}{2}}=\beta$$ Simple modules over $\Hab$ correspond to simple modules over $\Hg$ with central character $(\alpha,\beta)$.

For any $(\alpha,\beta)$, $\{1,x,..,x^{p-1}\}\otimes\{ 1,y,...,y^{p-1}\}$ is a $k$-basis for $\Hab$. Form a left module $\Sab$ over $\Hab$ by quotienting $\Hab$ by the sum of the left ideals generated by $(x-\sqrt[2p]{\alpha})$ and $(g-1)$: $$\Sab:=\Hab/(\Hab(x-\sqrt[2p]{\alpha})+\Hab(g-1))$$
\claim $\Sab$ is simple of dimension $p$.
\begin{proof}
It is clear that $1,y,y^2,...,y^{p-1}$ span $\Sab$. Let $f(y)\in k[y]$ and suppose $f(y)\in\Sab$ generates a proper submodule. Let $d=\textrm{degree }(f(y))$, $d<p$, and write $f(y)=y^d+a_1y^{d-1}+...+a_d$. Without loss of generality assume deg $f$ is minimal among those generating proper submodules.

Suppose first that $(\alpha,\beta)$ does not define a singular point of the center, i.e. $\alpha\neq 0$. Now consider $gf(y)$:
\[
gf(y)=g\cdot f(y)g=g\cdot f(y)=(x+y)^d+a_1(x+y)^{d-1}+a_2(x+y)^{d-2}+...+a_d
\]
The coefficient of $y^{d-1}x$ in $(x+y)^d$ is given by ${d\choose1}=d$ so in expanding $(x+y)^d$, the following term appears: $$dy^{d-1}x=dy^{d-1}\sqrt[2p]{\alpha}\neq0$$
since $\alpha, d\neq 0$. Therefore $gf(y)=f(y)+h(y)$ where $h(y)\in k[y]$ has leading term $d\sqrt[2p]{\alpha}y^{d-1}$, and so $h(y)$ is an element of degree $d-1$ in the cyclic module generated $f(y)$; in particular, it must too generate a proper submodule of $\Sab$, but that contradicts the assumption that $d$ is minimal.

Now let's consider the situation when $\alpha=0$, corresponding to a singular point of Spec $ Z(\Hg)$. Thus in $\Hzerob$, $Y=0$, $X=-2$, and $Z=\beta$, and to form the module $\Mzerob$, quotient $\Hzerob$ by the left ideal $\Hzerob x+\Hzerob(g-1)$. If $d>1$, $gf(y)$ contains the term ${d\choose2}y^{d-2}$ as a lower-order term and after $f(y)$ is subtracted from $gf(y)$ there is nothing to cancel ${d\choose2}y^{d-2}$ (since degrees of lower order terms drop by increments of $2$, a term of degree $d-1$ in $f$ is irrelevant here); thus ${d\choose2}y^{d-2}$ plus whatever other lower-order terms are left over in $gf(y)-f(y)$ is an element of lesser degree than $d$ which also generates a proper submodule, contradiction. On the other hand, if $d=1$ and $f$ generates a proper submodule then $yf$ is of degree $2$ but then $yf$ generates everything, and thus so does $f$. So $\Mzerob$ is also simple.
\end{proof}

\corollary The Azumaya locus and the nonsingular locus of $\Hg$ do not coincide.
\begin{proof}
The Azumaya locus of the center is all of $\Spec(Z(\Hg))$, while the nonsingular locus is $\Spec(Z(\Hg))$ minus a line.
\end{proof}
For simple modules in the singular locus, here are explicit matrices that give the action of $g$, $y$, and $x$. Let $\delta=g\del_x+xg\del_g$ and let $A_0(x,g)=g,\;A_1(x,g)=xg,\; A_2(x,g)=x^g+g^2,...,A_n(x,g)=\delta(A_{n-1}(x,g)),...$ be the Andr\'{e} polynomials. Set $\alpha_n=A_n(0,1)$. Note that for $n$ odd, $\alpha_n=0$. Set $\alpha_{i,j}={j-1\choose i-1}\alpha_{i},\; i<j$. Then as $\beta$ ranges through $k$, there are distinct simple representations of $\Hg$ over the singular locus given by:
\begin{align*}
&g=\begin{pmatrix}1&0&\alpha_2&0&\alpha_4&0&\alpha_6&0&...\\
0&1&0&{3\choose2}\alpha_2&0&{5\choose4}\alpha_4&0&{7\choose6}\alpha_6&...\\
0&0&1&0&{4\choose2}\alpha_2&0&{6\choose4}\alpha_4&0&...\\
0&0&0&1&0&{5\choose2}\alpha_2&0&{7\choose4}\alpha_4&...\\
0&0&0&0&1&0&{6\choose2}\alpha_2&0&...\\
&&&&...\\&&&&...\\&&&&...
 \end{pmatrix}\\
&x=\begin{pmatrix}0&1&0&\alpha_2&0&\alpha_4&0&...\\ 0&0&2&0&{4\choose3}\alpha_2&0&{6\choose5}\alpha_4&...\\0&0&0&3&0&{5\choose3}\alpha_2&0&....\\0&0&0&0&4&0&{6\choose3}\alpha_2&...\\0&0&0&0&0&5&0&...\\&&&&...\\&&&&...\\&&&&...\end{pmatrix}
\\
&\qquad y=\begin{pmatrix}0&0&0&0&...&0&\beta\\1&0&0&0&...&0&T_{\frac{p+3}{2}}\\0&1&0&0&...&0&0\\0&0&1&0&...&0&0\\&&&&...\\0&0&0&0&...&1&0\end{pmatrix}
\end{align*}
$T_j$ denotes the $j$th reduced tangent number modulo $p$: see sequence \href{http://oeis.org/A002105}{A002105} in the OEIS which begins $T_1=1,\;T_2=1,\;T_3=4\;,T_4=34,\;T_5=496,\;...$, $T_n=2^n(2^{2n}-1)\frac{|B_{2n}|}{n}$ where $B_j$ is the $j$th Bernoulli number.

\subsubsection{Finite-dimensional quotients of $\Hlam$}
The simple modules over $\Hlam$ may be found by the same procedure of central reduction that was used for $\Hg$, though we'll shift gears and use Verma modules for $\Hlam$ to define the simple modules in the end. The description of simple modules over $\Hlam$ depends on whether $\lambda\in kG$ is invertible or a zero-divisor. Let $\lambda=\sum\limits_{i=1}^{p-1}c_ig^i$, then $\lambda^p=\sum\limits_{i=1}^{p-1}c_i^p$, so $\lambda$ is invertible if $p$ does not divide $c_{\lambda}:=\sum\limits_{i=1}^{p-1}c_i$, while $\lambda$ is zero-divisor, in fact nilpotent, if $\sum\limits_{i=1}^{p-1}c_i$ is divisible by $p$.
If $\lambda\in kG$ is invertible then the simple modules for $\Hlam$ behave the same as those for $\Hg$:  all simple representations of $\Hlam$ are $p$-dimensional and the smooth locus of the center of the algebra is strictly contained in its Azumaya locus which equals Spec $\Zlam$. When $\lambda$ is a zero-divisor, the simple modules behave differently at the singular locus of the center. 

This dichotomy results from two possible different structures for $\Hlam$ based on the nature of $\lambda$:
\proposition Let $\lambda=\sum\limits_{i=1}^{p-1}c_ig^i$ and suppose $p$ divides $c_{\lambda}:=\sum\limits_{i=1}^{p-1}c_i$, i.e. $\lambda$ is a zero-divisor. Then $\Hlam$ surjects onto $k[y]$ with kernel generated by the (2-sided) ideal $(g-1)$. If $p$ does not divide $c_{\lambda}$, i.e. if $\lambda$ is an invertible element of $kG$, then $\frac{\Hlam}{(g-1)}=0$.
\begin{proof} Let $\phi:\Hlam\onto\frac{\Hlam}{(g-1)}$. Recall that $g$ and $y$ generate $\Hlam$. We check what happens to the relations. Then $\phi(y)=\phi(gy)=\phi(yg+xg)=\phi(y)+\phi(x)$, so the image of $x$ is $0$ in the quotient. Thus $[x,y]=0$ in the quotient, and this is compatible with the relation $[x,y]=\lambda$ if $\lambda$ is a zero-divisor since then $\phi(\lambda)=c_{\lambda}=0$. If $\lambda$ is not a zero-divisor then $0=\phi(\lambda)=c_{\lambda}\neq0$ and thus $\frac{\Hlam}{(g-1)}=0$.

\end{proof}

\corollary A $1$-dimensional representation of $\Hlam$ exists if and only if $\lambda\in kG$ is nilpotent. In the latter case, $\Hlam$ has infinitely many (simple) $1$-dimensional representations via setting $g=1$ and $y=\beta$ for each $\beta\in k$. 

\begin{proof}As for the statement that if $\lambda$ is not a zero-divisor there cannot be any $1$-dimensional representation of $\Hlam$: $G$ must act by a character on any $1$-dimensional representation, but the cyclic group of order $p$ has only the trivial character over a field of the same characteristic. Therefore $\lambda$ must act by the scalar $c_\lambda\neq0$ if $\lambda$ is invertible, so $xy-yx=\lambda$ must act by $c_\lambda$, but $xy-yx$ acts by $0$ on any $1$-dimensional representation.\end{proof}
Note that $x$ must act by $0$ on any $1$-dimensional representation, because of the relation $gy=xg+yg$.

\vspace{5mm}

Fix $\lambda\in kG$. Given a finite-dimensional representation $V$ of $\Hlam$, the center of $\Hlam$ must act by scalars on it, and we call the scalars by which the generators of the center act the ``central character."  If $Y$ acts by $\alpha$ then $X$ must act by $\sqrt[p]{\alpha^2}$ and so we simply write $(\alpha,\beta)$ instead of $(\sqrt[p]{\alpha^2},\alpha, \beta)$, and then modules over $\Hlam$ with central character $(\alpha,\beta)$ correspond to modules over the finite-dimensional algebra $$\Hab:=\frac{\Hlam}{(X-\sqrt[p]{\alpha^2}, Y-\alpha, Z-\beta)}$$

\proposition $\Hlam$ is a finite-dimensional module over its center of generic rank $p^2$.
\begin{proof} The point $(\alpha,\beta)$ is nonsingular point in Spec $\Zlam$ when $\alpha\neq0$, in which case $x=\frac{1}{\alpha }x^{p+1}=\frac{1}{\alpha }(x^2)^{\frac{p+1}{2}}$ and as $x^2$ is linearly dependent on $kG$, then $x$ is too. Thus for a smooth point $(\alpha,\beta)$ the algebra $\Hab$ is of dimension $p^2$.
\end{proof}

While $\Hab$ has dimension $p^2$ for all $\alpha\neq0$, the minimal possible dimension of $\Hzerob$ is $p^2$ and the maximal possible dimension is $2p^2$, and both of these bounds are realized as $\lambda$ ranges through $kG$ and are the only dimensions occurring. In fact, we will show that the set of $\lambda$ such that the quotient algebras $\Hzerob$ of $\Hlam$ have the maximal possible dimension consists of exactly those $\lambda$ which are zero-divisors in $kG$. On the other hand, if $\lambda$ is invertible then dim($\Hzerob)=p^2$.

Consider the three relations defining $\Hzerob$ as a quotient of $\Hlam$:
\begin{align}
&x^2-2\sum_{a=1}^{p-1}\frac{c_a}{a}g^a+2\sum_{a=1}^{p-1}\frac{c_a}{a}=0\\
&x^p=0\\
&y^p-y(x^2-2\sum_{a=1}^{p-1}\frac{c_a}{a}g^a)^{\frac{p-1}{2}}=\beta
\end{align}
As usual, $\sum\limits_{a=1}^{p-1}c_ag^a=\lambda$. Because of the first relation, the third relation becomes a polynomial in $k[y]$. $x^2$ is dependent on $kG$ but $x$ is not. Thus the dimension of $\Hzerob$ will be $2p^2$ if the relations (3.6) and (3.7) do not imply a nontrivial relation among the elements of a basis for $kG$. Writing $x^2=2\sum\limits_{a=1}^{p-1}\frac{c_a}{a}(g^a-1)$, and observing that $0=x^{p+1}=(x^2)^{\frac{p+1}{2}}$, it follows that the dimension of $\Hzerob$ is $2p^2$ if and only if the following equation holds in $kG$:

\begin{align}
&\left(\sum_{a=1}^{p-1}\frac{c_a}{a}(g^a-1)\right)^{\frac{p+1}{2}}=0
\end{align}

Recall that the set of zero-divisors in $kG$ are the elements belonging to the augmentation ideal which has a basis $(g-1), (g^2-1),...,(g^{p-1}-1)$. If $\mu$ is a zero-divisor then, since $(a+b)^p=a^p+b^p$ in characteristic $p$ whenever $a$ and $b$ commute, $\mu^p=0$. Let $\del_g$ be the operator of differentiation with respect to $g$, and $g$ the operator of multiplication by $g$, and consider the action of $g\del_g$ on $kG$: $$g\del_g(\sum_{a=0}^{p-1}m_ag^a)=\sum_{a=1}^{p-1}am_ag^a.$$

\proposition Let $$\mu:=\sum_{a=1}^{p-1}m_a(g^a-1)$$ be an arbitrary zero-divisor of $kG$. Then $\mu^{p-1}\neq0$ if and only if $g\del_g(\mu)$ is invertible. If $g\del_g(\mu)$ is a zero-divisor then $\mu^{\frac{p+1}{2}}=0$.

\begin{proof} Suppose that $\lambda:=g\del_g(\mu)$ is invertible. Then $$\lambda^{-1}g\del_g(\mu^n)=n\mu^{n-1}$$ 
and thus $$(\lambda^{-1}g\del_g)^{p-2}(\mu^{p-1})=(p-1)!\mu\neq0\iff\mu\neq0$$ So $\mu^{p-1}$ could not be $0$.

Now suppose $\lambda=g\del_g(\mu)$ is a zero-divisor. Assume $p>3$. It is easy to verify the statement for $p=3$ by hand. First we claim that $\mu^{p-1}=0$. We know $$g\del_g(\mu^{p-1})=-\lambda\mu^{p-2}$$ must square to $0$ as $2p-4>p\implies\mu^{2p-4}=0$. Observe that any element $\zeta\in kG$ such that $\zeta^2=0$ generates an ideal that is a one-dimensional representation of $G$. We know $G$ has a unique one-dimensional representation, namely that spanned by the sum of all the elements of $G$: $1+g+...+g^{p-1}$. Thus $\zeta^2=0$ implies $\zeta=C\cdot(1+g+...+g^{p-1})$ for some scalar $C\in k$. Therefore $g\del_g(\mu^{p-1})$ must be a scalar multiple of $1+g+g^2+...+g^{p-1}$, which is impossible unless that scalar is zero, since the result of $g\del_g$ applied to anything cannot contain a constant term. So $$g\del_g(\mu^{p-1})=0$$
from which it follows that $\mu^{p-1}$ is a constant, and since $\mu$ is a zero-divisor, this implies $\mu^{p-1}=0$. Since $\mu^{p-1}=0$, $(\mu^{\frac{p-1}{2}})^2=0$ and thus $$\mu^{\frac{p-1}{2}}=C\cdot(1+g+g^2+...+g^{p-1})$$ for some constant $C$ (apparently given by a homogenous polynomial $f$ of degree $\frac{p-1}{2}$ in $k[x_1,...,x_{p-2}]$ evaluated at $c_1,...,c_{p-2}$ with $c_i=\frac{m_a}{a}$). Since $g$ acts trivially on the line spanned by $(1+g+...+g^{p-1})$ and the sum of the coefficients of $\mu$ is equal to $0$, we have:
$$\mu^{\frac{p+1}{2}}=\mu\cdot\mu^{\frac{p-1}{2}}=C\mu\cdot(1+g+g^2+...+g^{p-1})=0$$

\end{proof}

\corollary If $\lambda$ is a zero-divisor, then the central quotient algebras $\Hzerob$ of $\Hlam$, that is those defined over the singular line of Spec Z($\Hlam$), are of the maximal possible dimension, namely $2p^2$. If $\lambda$ is invertible then the algebras $\Hzerob$ are of dimension $p^2$.

\begin{proof}
Let $\mu=\sum_a\frac{c_a}{a}(g^a-1)$. Then $g\del_g(\mu)=\sum\limits_{a=1}^{p-1}c_ag^a=\lambda$, and $2\mu=x^2$ in $\Hzerob$.  By Proposition 5.11 and the defining relations of $\Hzerob$,

\begin{align*}
&\lambda\textrm{ is a zero-divisor }\iff\mu^{\frac{p+1}{2}}=0\\
&\qquad\iff\textrm{ the equations }(x^2)^{\frac{p+1}{2}}=0,\;(x^2)^{\frac{p+3}{2}}=0,\;\textrm{etc., do not impose relations in $kG$}\\
&\lambda\textrm{ is invertible }\iff\mu^{p-1}\neq0\\
&\qquad\iff\textrm{ the equations }(x^2)^{\frac{p+1}{2}}=0,\;(x^2)^{\frac{p+3}{2}}=0,\;\textrm{etc., impose $\frac{p-1}{2}$ relations in $kG$}
\end{align*}
\end{proof}

\vspace{5mm}

Let's return to the description of simple modules over $\Hlam$. We could construct a left module $\Sab$ over $\Hab$ by taking $$\Sab:=\Hab/\Hab(x-\alpha)+\Hab(g-1).$$ which would have $k$-basis given by $y^i$, $i=0,...,p-1$, and then investigate for which $(\alpha,\beta)$ it is a simple module. However, it's easier to get the simple modules as quotients of Verma modules by maximal ideals.

\subsubsection{Simple modules as quotients of Verma modules}
Define Verma modules for $\Hlam$ as $$\Da:=\Hlam/\left(\Hlam(x-\alpha)+\Hlam(g-1)\right)$$
as $\alpha$ ranges through $k$. As a left $\Hlam$-module, $\Da$ is isomorphic to $k[y]$. When $t=0$, $\Da$ contains an infinite descending chain of submodules generated by $(C-\beta)^n$ for each $n\in\bN$, $C=y^p-y\frac{\delta^p(g)}{\delta(g)}$ the central element. An easy criterion for when a submodule is proper does not make reference to the special cases $t=0$ or $t=1$:

\lemma Let $f(y)\in\Da$ generate a submodule of $\Da.$ Then $f(y)$ generates a proper submodule if and only if $g\cdot f(y)=f(y)$ and $x\cdot f(y)=c f(y)$ for some $c\in k$.
\begin{proof}
We may assume the degree of $f(y)$ in $y$ is minimal over elements of $\Hlam f(y)\subset\Da$. Multiplication by $y$ always raises the degree. For any $r\in\R$, $rf=fr+$ lower order terms, so the ``only if" is clear. Equally obvious: suppose $r\in\R$ acts by a scalar on $f(y)$. Writing $h\in\Hlam$ in its normal form $h=\sum y^i r_i$, $r\in\R$, the degree of $hf(y)$ is at least that of $f(y)$ unless $hf(y)=0$ since $r_i$ multiplies $f(y)$ by a scalar and $y^i$ raises the degree.
\end{proof}

\proposition Suppose $\lambda\neq0$ does not contain a constant term. There is a two-parameter family of distinct simple $\Hlam$-modules $\Sab$, $\alpha,\beta\in k$. If $\lambda$ is invertible, all simple $\Hlam$-modules are $p$-dimensional, while if $\lambda$ is a zero-divisor then the simple $\Hlam$-modules over the smooth locus are $p$-dimensional while the simple modules over the singular locus are $1$-dimensional.
\begin{proof}
Suppose $\alpha\neq0$. Then $\delta(g)=\alpha\neq0$ in $\Da$. Suppose $f(y)=y^n+a_{n-1}y^{n-1}+...\in\Da$, $0<n<p$. Then compute the action of $g$ using the multiplication in $\Hlam$ to move $kG[x]$ terms to the right then reduce modulo the ideals defining $\Da$: \begin{align*}
g\cdot f(y)&=\left(\sum_{i=1}^n{n\choose i}y^{n-i}\delta^i(g)\right)+a_{n-1}\left(\sum_{i=1}^{n-1}{n-1\choose i}y^{n-1-i}\delta^i(g)\right)+...\\
&=y^n+a_{n-1}y^{n-1}+{n\choose 1}\alpha y^{n-1}+a_{n-2}y^{n-2}+....\\
\end{align*}
Since $\alpha\neq0$ and ${n\choose1}=n\neq0$ as $n<p$, $g$ cannot fix $f(y)$ and so, by the lemma, $f(y)$ cannot generate a proper ideal.

When $\alpha=0$, $g\cdot (y-\beta)=y-\beta$ and $x\cdot (y-\beta)=c_{\lambda}$ where $c_{\lambda}$ is the sum of the coefficients of $\lambda$. Thus $y-\beta$ generates a proper submodule if and only if $c_{\lambda}=0$ mod $p$. 

If $\alpha=0$ and $\lambda$ is not a zero-divisor then $\delta(x)=c_{\lambda}\neq0$, and by the same argument given in the $\alpha\neq0$ case but with $x$ in place of $g$, $x$ cannot act by a scalar on any polynomial of degree less than $p$, and so there can be no proper submodule generated in degree less than $p$. 

Since the shift by any $\beta\in k$ of the central element $Y=y^p-yx^{p-1}-...$ commutes with $x$ and $g$, the image $Y-\beta$ in $\Da$ generates a maximal submodule of $\Da$ when $\alpha\neq0$ or when $\alpha=0$ and $\lambda$ is invertible.
\end{proof}

The simple modules are then of the form $$\Sab=\Da/(Y-\beta)\Da$$ where $Y$ is the central element $y^p-yx^{p-1}-...$, except over the singular locus when $\lambda$ is a zero-divisor, where the simple modules have the form $$\Mzerob=\Delta_0/(y-\beta)\Delta_0$$

The representation theory of $\Hlam$ has consequences for the structure of $\Hlam$: existence of normal elements, radical, and so forth. \\

An element $v$ in an algebra $H$ is called normal if $vH=Hv$, that is, the left and right ideals generated by $v$ coincide. The description of the Verma modules implies:
\corollary The only normal elements of $\Hlam$ are the elements of the center $\Zlam$.
\begin{proof}
Suppose $v\in\Hlam$ is normal, and look at its image $\bar{v}\in\Da$. Since $v\Hlam=\Hlam v$, $\bar{v}$ generates a proper submodule of $\Da$. But all submodules of $\Da$, $\alpha\neq0$, are generated by powers of $Y$ and their scalar translates. So $v$ would have to be of the form $v=f(Y)+F$ with for some $f\in k[Y]$ and $F\in \Hlam(x-\alpha)+\Hlam(g-1)$. Since this must hold for all $\alpha$, $F$ must be in $\Hlam(g-1)$, $F=h(g-1)$ for some $h\in\Hlam$. Since the difference of normal elements is normal, $F$ must be normal. However, $$h(g-1)y=hy(g-1)+hgx=yh(g-1)+[h,y](g-1)+hgx$$ so  $\Hlam h(g-1)=h(g-1)\Hlam$ implies $[h,y](g-1)+hgx\in\Hlam h(g-1)$ and in particular $[h,y](g-1)+hgx\in\Hlam(g-1)$, so that we must have $$h\in\Hlam(g-1).$$ Writing $h=h'(g-1)$ and repeating the argument, we obtain $h'\in\Hlam(g-1)$, and so on; after the $p$th iteration, we have $$F\in\Hlam(g-1)^p=0.$$
\end{proof}

Let $\Rad\Hlam$ denote the Jacobson radical of $\Hlam$. \proposition$\Rad\Hlam=0$.

\begin{proof}
Take $z$ in $\Rad\Hlam$ and write $$z=y^nf_n(x,g)+y^{n-1}f_{n-1}(x,g)+...+yf_1(x,g)+f_0(x,g).$$
Suppose $n<p$. Since the Jacobson radical of any ring is a two-sided ideal, $gzg^{-1}\in\Rad\Hlam$ and thus $\Rad\Hlam$ contains the following nonzero element of smaller degree in $y$ than $z$:
\begin{align*}
gzg^{-1}-z=&(x+y)^nf_n(x,g)+(x+y)^{n-1}f_{n-1}(x,g)+...+(x+y)f_1(x,g)+f_0(x,g)\\
&-y^nf_n(x,g)-y^{n-1}f_{n-1}(x,g)-...-yf_1(x,g)-f_0(x,g)\\
=&y^{n-1}nxf_n(x,g)+\textrm{lower order terms in }y
\end{align*}
Therefore if $\Rad\Hlam$ contains an element of degree less than $p$ in $y$, then it contains an element of degree $0$, that is, it contains some $r\in R$.

So suppose $r\in R$ and $r\in\Rad\Hlam$. Then $r=0$ in any simple quotient $S$ of $\Hlam$. The simple quotients $S$ of $\Hlam$ are subquotients of the Verma modules $\Da=\Hlam/(\Hlam(g-1)+\Hlam(x-\alpha))$ as $\alpha$ runs through $k=\bar{F}_p$. Thus $g-1$ must divide $r$ in order for $r$ to be $0$ in $S$ for all simple $S$, as $x-\alpha$ can't divide $r$ for all $\alpha$ since $k$ is infinite. Write $r=r_1(g-1)$ for some $r_1\in R$. For generic $S$, a basis is given by $1,y,...,y^{p-1}$. We must have:
$$0=ry=yr+\delta(r)=0+\delta(r_1(g-1))=\delta(r_1)(g-1)+r_1\delta(g-1)=0+r_1xg$$
and so we must have $r_1xg=0$ in any simple $S$, but if $S$ is a quotient of $\Da$ then $$r_1xg=\alpha r_1$$ and thus $r_1=0$ in any simple $S$; it follows that $$r_1=r_2(g-1)$$ for some $r_2\in R$. Repeating this argument, it holds that in order for $r$ to annihilate each the basis elements $1,y,y^2,...,y^d$ in every simple $S$ of dimension $p$, it is necessary that $$r=r_{d+1}(g-1)^{d+1}$$ for some $r_{d+1}\in R$. But then in order for $ry^{p-1}$ to be $0$ (in addition to $ry^{d}=0$ for every $d<p-1$) in every simple $p$-dimensional module $S$, $(g-1)^p=0$ must divide $r$, and therefore $r=0$.

Now suppose $n=p$.  Then $$gzg^{-1}-z=\begin{cases} \textrm{a polynomial of degree }<p\textrm{ in }y&\textrm{ if }z\notin \Zlam\\0&\textrm{ if }z\in \Zlam\end{cases}$$
The first case cannot happen since we've just shown $\Rad\Hlam$ contains no element of degree $<p$ in $y$. The second case cannot happen since if $z=y^p+...\in \Zlam$, then $z$ acts by $\beta$ in the simple representation $\Sab$, so $z$ doesn't act by $0$ in every simple representation and so $z\notin J(\Hlam)$.
If $n>p$, we apply the same arguments to reduce to the cases already considered.
\end{proof}

An element $e$ in a ring $H$ is called idempotent if $e^2=e$. The existence of idempotents or lack there-of is important for understanding the structure of a ring and its modules.

\lemma Let $\R=k[x]\otimes kG$. Then $\R$ contains no idempotents besides $0$ and $1$.
\begin{proof}
Suppose $f(x,g)=x^n\lambda_n+x^{n-1}\lambda_{n-1}+...+x\lambda_1+\lambda_0$ is idempotent, $f(x,g)\neq0,1$. If $f(x,g)=f(g)\in kG$ then $f$ cannot be idempotent as the $p$th power of any element of $kG$ belongs to $k$, so $f\in k$ and so $f=0$ or $1$. Also, for any $f$ as above, $(\lambda_0)^2=\lambda_0\implies \lambda_0=0$ or $1$. Let $d$ be the smallest nonzero integer such that $\lambda_d\neq0$, i.e. $d$ is the lowest power of $x$ that appears in $f$. Then $f(x,g)^2=...+2\lambda_0\lambda_dx^d+\lambda_0=...+\lambda_dx^d+\lambda_0$, which is impossible unless $\lambda_d=0$, contradicting the minimality of $d$.
\end{proof} 
\corollary $\Hzero=k[x,y]\rtimes G$ contains no idempotents besides $0$ and $1$.
\begin{proof}
By the same argument, with respect to the minimality of the power of $y$ occurring in a given element.
\end{proof}

\proposition $\Hlam$ contains no idempotents besides $0$ and $1$.
\begin{proof}
Suppose $e\in\Hlam$ is idempotent. Let $\bar{e}$ denote the image of $e$ in $\gr\Hlam=\Hzero$. Then $$\bar{e}^2=\bar{e^2}=\bar{e}$$ so $\bar{e}$ is idempotent in $\Hzero$, hence $\bar{e}=0$ or $1$. Hence the degree of $e$ is zero, implying $e\in kG$. But $kG$ has no nontrivial idempotents.
\end{proof}
\corollary No simple $\Hlam$-module has a projective cover.
\begin{proof}
 \cite{rings}, Proposition 17.19.
\end{proof}

\subsection{Extensions between simple modules, and $\Hlam$-modules from the viewpoint of noncommutative algebraic geometry} We restrict our attention to the $t=0$ case in this section.
It happens that each simple module has nonzero Ext groups with itself but that unless $\lambda\in \Rad kG\backslash \Rad^2 kG$, all Ext groups between distinct simple modules vanish. No matter the choice of $\lambda$, simple modules from different Verma modules have no nonzero Ext groups with each other. Within a fixed Verma, the exceptional case of nonzero Ext groups occurring between distinct simples when $\lambda\in \Rad kG\backslash \Rad^2kG$ happens over the singular locus. Specifically it only occurs between simple modules $S_{0,\beta}$ and $S_{0,\beta'}$ over the singular locus when the difference between $\beta$ and $\beta'$ equals the square root of the sum of the coefficients of $\frac{\lambda}{g-1}$.

\proposition Let $S$ be a simple $\Hlam$-module of dimension $p$. Then
$$\Ext_{\Hlam}^i(S,S)=
\begin{cases} S & i=0\\
S^{\oplus3} & i=1\\ S^{\oplus4} & i\geq2\end{cases}
$$

\begin{proof}
Without loss of generality, assume $\beta=0$ so that $C=y^p-y\frac{\delta^p(g)}{\delta(g)}$ acts on $S$ by $0$ (recall that the behavior of the simple modules over $\Spec \Zlam$ depends only on the parameter coming from the central elements in $k[x]\otimes kG$ together with whether or not $\lambda$ is invertible). Let $\alpha$ be such that $x$ acts on $S$ by $\alpha$. So we have $S=S_{\alpha,0}$. 

Set $H:=\Hlam$. The following complex is a free resolution of $S$ by $H$-modules, where elements of the modules are thought of as row-vectors and the maps are right multiplication by the matrices:
\begin{align*}
... \stackrel{B}{\lra}H^{\oplus4}\stackrel{A}{\lra} H^{\oplus4}&\stackrel{B}{\lra}H^{\oplus4}\stackrel{A}{\lra} 
H^{\oplus4}\xrightarrow{\begin{pmatrix}\sum\limits_ag^a&\alpha-x&C&0\\0&g-1&0&-C\\0&0&g-1&\alpha-x\\0&0&0&\sum\limits_ag^a\end{pmatrix}}\\
&\lra H^{\oplus4}\xrightarrow{\begin{pmatrix} g-1&x-\alpha&-C&0\\0&\sum\limits_ag^a&0&C\\0&0&\sum\limits_ag^a&x-\alpha\\0&0&0&g-1\end{pmatrix}}\\\lra &H^{\oplus4}\xrightarrow{\begin{pmatrix} C & 0& \alpha-x\\ 0 &- C & g-1\\ g-1 & \alpha-x & 0\\0&\sum\limits_ag^a&0\end{pmatrix}} H^{\oplus3} \xrightarrow{\begin{pmatrix} x-\alpha\\g-1\\C\end{pmatrix}}H\lra 0
\end{align*}

$A$ and $B$ denote respectively the second-to-last and last matrices written explicitly above (going from right to left and bottom to top up the complex). In $\sum\limits_ag^a$, the sum is over $a=0,1,2,...,p-1$. Note that all these matrices have entries in the maximal commutative subalgebra $(k[x]\otimes kG)[C]$. Also, the annihilator of $\sum\limits_ag^a$ is $H(g-1)$ while the annihilator of $g-1$ is $H\sum\limits_ag^a$, which is why these terms keep appearing in the maps above. The maps $A$ and $B$ keep repeating in alternation indefinitely, and the resolution never terminates. 

Taking $\Hom_H(-,S)$ of the resolution, we get the complex
$$...\lla S^{\oplus4}\stackrel{0}{\lla} S^{\oplus4}\stackrel{0}{\lla} S^{\oplus4}\stackrel{0}{\lla} S^{\oplus3}\stackrel{0}{\lla}S\lla0$$
Every map is the zero map because $S$ is a quotient of $H$ which kills the left ideals $H(x-\alpha)$, $H(g-1)$, and $HC$; since $\sum\limits_ag^a$ evaluated at $g=1$ is $p$, H($\sum\limits_ag^a$) maps to $0$ in $S$ as well; and thus precomposition with each of the maps in the resolution of $S$ sends an arbitrary element of $H$ into the kernel of projection onto $S$. Therefore the $\Ext$ groups are just given by the modules appearing in the complex.
\end{proof}
In the situation that $\lambda\in\Rad kG$ and $S$ is a simple module over a point in the singular locus of $\Spec\Zlam$, the ranks of the $\Ext$ groups depend on whether $g-1$ divides $\frac{\lambda}{g-1}$, i.e. whether $\lambda\in\Rad^2 kG$.
\proposition Let $H=\Hlam$ with $\lambda$ a zero-divisor in $kG$. Observe that $g-1$ divides $\lambda$. Consider $S$, a simple $H$-module of dimension $1$. Suppose $(g-1)^2$ does not divide $\lambda$. Then
$$\Ext_H^i(S,S)=\begin{cases} S & i=0,1\\S^{\oplus2} & i=2\\ S^{\oplus3} & i\geq3\end{cases}$$
On the other hand, if $(g-1)^2$ divides $\lambda$ then 
$$\Ext_H^i(S,S)=\begin{cases} S& i=0\\S^{\oplus2} & i=1\\ S^{\oplus 3} & i=2\\ S^{\oplus4}&i\geq3\end{cases}$$
\begin{proof}
The assumption that $S$ is one-dimensional means that $S=\Mzerob$ for some $\beta\in k$, i.e. $x$ acts on $S$ by $0$, $y$ by $\beta$, and $g$ by $1$. The following complex is a free resolution of $S$ over $H$:
\begin{align*}
...\stackrel{B}{\lra}H^{\oplus4}\stackrel{A}{\lra}H^{\oplus4}\stackrel{B}{\lra}&H^{\oplus4}\stackrel{A}{\lra}H^{\oplus(4)}\xrightarrow{\begin{pmatrix}\sum\limits_ag^a&x&-(y-\beta)&\frac{\lambda}{g-1}\\
0&g-1&-g&y-\beta\\0&0&g-1&x\\0&0&0&\sum\limits_ag^a\end{pmatrix}}\\
&\lra H^{\oplus4}\xrightarrow{\begin{pmatrix}g-1&-x&y-\beta&0\\
0&\sum\limits_ag^a&g(g-1)^{p-2}&-(y-\beta) \\0&0&\sum\limits_ag^a&-x\\0&0&0&g-1\end{pmatrix}}\\
&\lra H^{\oplus4}\xrightarrow{\begin{pmatrix}y-\beta&\frac{\lambda}{g-1}&-x\\g&y-\beta&1-g\\1-g&x&0\\0 & \sum\limits_ag^a&0\end{pmatrix}}H^{\oplus3}\xrightarrow{\begin{pmatrix}x\\g-1\\y-\beta\end{pmatrix}}H\lra0
\end{align*}

To see that it is exact, we make a few elementary observations: first, that the annihilator of $\lambda\in kG$ a zero-divisor is the left ideal $H(\sum\limits_ag^a)$ (i.e., exactness at the $4\times3$ step is mostly clear except for where there are terms introduced to kill $\lambda=xy-yx$ as it arises, where we must make sure we included the the whole annihilator of $\lambda$). Secondly, useful identities from $kG$: \begin{align*}&(g-1)^{p-1}=\sum\limits_ag^a\\
&g\sum\limits_ag^a=\sum\limits_ag^a\\&(\sum\limits_ag^a)^2=0\\&g(1-g)^{p-2}=\sum_aag^a=g\del_g(\sum_ag^a).\end{align*} 


Applying $\Hom_{H}(-,S)$ to the resolution of $S$ gives the complex

\begin{align*}
0\to S\stackrel{0}{\lra}S^{\oplus3}\xrightarrow{\begin{pmatrix}0&1&0&0\\ \frac{\lambda}{g-1}&0&0&0\\0&0&0&0\end{pmatrix}}S^{\oplus4}\stackrel{0}{\lra}S^{\oplus4}\xrightarrow{\begin{pmatrix}0&0&0&0\\0&0&0&0\\0&0&0&0\\\frac{\lambda}{g-1}&0&0&0\end{pmatrix}}S^{\oplus4}\stackrel{0}{\lra}S^{\oplus4}\stackrel{\tilde{B}}{\lra}...
\end{align*}
Here $\tilde{B}$ denotes the matrix with all $0$'s except for $\frac{\lambda}{g-1}$ in the lower left corner. The maps $0$ and $\tilde{B}$ continue in alternation indefinitely to the right. Now there are two cases: 

\textbf{Case I.} If $\lambda$ is divisible by $g-1$ but not by $(g-1)^2$ then $\frac{\lambda}{g-1}$ is invertible and its role in the two matrices above is multiplication by the nonzero scalar determined by the action of $\frac{\lambda}{g-1}$ on $S$. Then at $S^{\oplus3}$ the kernel is $1$-dimensional, so that $\Ext^1(S)=S$. The image of the $3\times 4$ matrix is $S^{\oplus2}$ in the first two components of $S^{\oplus4}$, so that $\Ext^2(S)=S\oplus S$. At every $S^{\oplus4}$ from there onwards, either the kernel is $4$-dim and the image $1$-dim, or the image is $0$-dim and the kernel $3$-dim.

\textbf{Case II.} If $\lambda$ is divisible by $(g-1)^2$ then $\frac{\lambda}{g-1}$ is divisible by $g-1$, hence must act by $0$ on $S$, and so every map above is the zero map except the $3\times 4$ matrix, which has $2$-dimensional kernel and $1$-dimensional image.

\end{proof}

As a consequence of Propositions 5.21 and 5.22,

\corollary The algebras $\Hlam$ have infinite global dimension.

Each Verma module $\Delta_{\alpha}$, possesses a line's worth of simple quotient modules $\Sab$ arising as the quotients of $\Delta_{\alpha}$ by $C-\beta$ for each $\beta\in k$. For distinct simple modules $\Sab$ and $\Sabprime$ from the same Verma $\Delta_{\alpha}$ belonging to the Azumaya locus, all $\Ext$ groups vanish.

\proposition Suppose the simple quotients of $\Delta_{\alpha}$ are $p$-dimensional. Then $$\Ext_{\Hlam}^i(\Sab,\Sabprime)=0$$ for all $i\geq0$, $\beta\neq\beta'\in k$.

If the simple quotients of $\Delta_0$ are $1$-dimensional, i.e. if $\lambda\in \Rad kG$, then it is also true that   $$\Ext_{\Hlam}^i(\Mzerob,\Mzerobprime)=0$$ for all $i\geq0$, $\beta\neq\beta'\in k$, except in the following special case: if $\lambda\notin(\Rad kG)^2$, so that the sum $c$ of the coefficients of $\frac{\lambda}{g-1}$ is nonzero, and if $(\beta'-\beta)^2=c$, then $$\Ext_{\Hlam}^i(\Mzerob,\Mzerobprime)=\begin{cases} 0 & i=0\\ \Mzerobprime & i=1,2\\ 0 & i\geq3\end{cases}$$
\begin{proof}
Consider first the case where the simple modules are $p$-dimensional. Let $\alpha$ be any element of $k$ if $\lambda$ is invertible, or let $\alpha\in k^{\times}$ if $\lambda\in\Rad kG$. Suppose $\beta\neq\beta'\in k$ and let $\Sab$, $\Sabprime$ be the corresponding $p$-dimensional simple quotients of $\Delta_{\alpha}$. In the free resolution of $\Sab$ from Proposition 5.19, replace $C$ with $C-\beta$ and then apply $\Hom(-,\Sabprime)$ to obtain the complex:
$$0\lra\Sabprime\xrightarrow{\begin{pmatrix}0,&0,&C-\beta\end{pmatrix}}\Sabprime^{\oplus3}\xrightarrow{\begin{pmatrix} C-\beta&0&0&0\\0&\beta-C&0&0\\0&0&0&0\end{pmatrix}}\Sabprime^{\oplus4}\xrightarrow{\tilde{A}}\Sabprime^{\oplus4}\xrightarrow{-\tilde{A}}\Sabprime^{\oplus4}\xrightarrow{\tilde{A}}...$$
where $$\tilde{A}:=\begin{pmatrix}0&0&0&0\\0&0&0&0\\\beta-C&0&0&0\\0&C-\beta&0&0\end{pmatrix}$$ and the resolution goes on forever to the right, alternating between $\tilde{A}$ and $-\tilde{A}$. As $\beta\neq\beta'$ the first map is injective and its image the copy of $\Sabprime$ in the third component. The kernel of the second map (the $4\times 3$ matrix) is $\Sabprime$ in the third component and its image is $\Sabprime\oplus\Sabprime$ in the first two components. Thus we have $\Ext^0(\Sab,\Sabprime)=0$ and $\Ext^1(\Sab,\Sabprime)=0$. The kernel of $\tilde{A}$ and $-\tilde{A}$ is $\Sabprime\oplus\Sabprime$ in the first two components, which coincides with the image of both maps. Thus all Ext groups vanish.

Now suppose $\lambda\in \Rad \Hlam$ and $\alpha=0$. Take $\beta\neq\beta'$ and let $\Mzerob$ and $\Mzerobprime$ be the corresponding $1$-dimensional simple quotients of $\Delta_0$. Apply $\Hom(-,\Mzerobprime)$ to the free resolution of $\Mzerob$ in Proposition 5.20 to obtain the complex
$$0\lra\Mzerobprime\xrightarrow{(0,0,\beta'-\beta)}\Mzerobprime^{\oplus3}\xrightarrow{\begin{pmatrix} \beta'-\beta&1&0&0\\\frac{\lambda}{g-1}&\beta'-\beta&0&0\\0&0&0&0\end{pmatrix}}\Mzerobprime^{\oplus4}\xrightarrow{\tilde{A}}\Mzerobprime^{\oplus4}\xrightarrow{\tilde{B}}\Mzerobprime^{\oplus4}\xrightarrow{\tilde{A}}...$$
where $$\tilde{A}=\begin{pmatrix}0&0&0&0\\0&0&0&0\\\beta'-\beta&0&0&0\\0&\beta-\beta'&0&0\end{pmatrix}$$ and $$\tilde{B}=\begin{pmatrix}0&0&0&0\\0&0&0&0\\\beta-\beta'&0&0&0\\\frac{\lambda}{g-1}&\beta'-\beta&0&0\end{pmatrix}$$

Call the $1\times 3$ matrix $\phi_0$, the $3\times 4$ matrix $\phi_1$, and denote by $\phi_2$ the first instance of $\tilde{A}$. Since $\beta'\neq\beta$, $\phi_0$ is injective with image the copy of $\Mzerobprime$ in the third component. Thus $\Ext_{\Hlam}^0(\Mzerob,\Mzerobprime)=0$. The kernel and image of $\tilde{A}$ and $\tilde{B}$ all coincide and are equal to $\Mzerobprime\oplus \Mzerobprime$ placed in the first two components. Thus $\Ext_{\Hlam}^n(\Mzerob,\Mzerobprime)=0$ for all $n\geq3$. The kernel and image of $\phi_1$ depend on the sum $c$ of the coefficients of $\frac{\lambda}{g-1}$.

\textbf{Case I.}  Suppose $\frac{\lambda}{g-1}$ either belongs to $\Rad\Hlam$ or if it is invertible, that the sum $c$ of its coefficients is different from $(\beta'-\beta)^2$. Then Ker $\phi_1=\{ 0,0,m\}, m\in \Mzerobprime=$ Im $\phi_0$, and Im $\phi_1=\{(m_1, m_2,0,0)\}, m_1, m_2\in \Mzerobprime=$ Ker $\phi_2$. Therefore $$\Ext_\Hlam^i(\Mzerob,\Mzerobprime)=0\textrm{ for }i=1,2.$$

\textbf{Case II.} If $c=(\beta'-\beta)^2$ then the line $\{(a(\beta'-\beta),a,0,0)\}, a\in \Mzerobprime$ also belongs to Ker $\phi_1$ and thus $$\Ext_{\Hlam}^1(\Mzerob,\Mzerobprime)=\textrm{ Ker }\phi_1/\textrm{Im }\phi_0=\Mzerobprime.$$ Also, in this case Im $\phi_1=\{((\beta'-\beta)a,a,0,0)\} a\in \Mzerobprime$ is only a line so that $$\Ext_{\Hlam}^2(\Mzerob,\Mzerobprime)=\textrm{ Ker }\phi_2/\textrm{Im }\phi_1=\Mzerobprime.$$

\end{proof}

Rather than compute the Ext groups between simple modules on different fibers, we apply some general results on modules over Ore extensions. Work by Smith and Zhang examines the module structure over an Ore extension $H$ of a commutative ring $R$ \cite{SZ}. The motivating idea is that the inclusion $R\into H$ induces a surjection $\Spec H\onto \Spec R$. $\Spec R$ is a genuine variety while  ``$\Spec H$", defined to be $\Mod H$, is a noncommutative space. The simple modules over $H$ can all be found in the fibers of this map. In the case of $\R=k[x]\otimes kG$ and $\mathsf{H}=\Hlam=\R[y;\delta=\lambda\del_x+xg\del_g]$, the fiber module $F_{\alpha}$ as defined in section 3 of Smith and Zhang \cite{SZ} coincides with the Verma module $\Delta_a$. As a consequence of Lemma 8.2 of \cite{SZ}, we may state:
\proposition Let $\Delta_{\alpha}$ and $\Delta_{\alpha'}$ be distinct Verma modules. Then for all simple modules $\Sab$, $S_{\alpha',\beta'}$ obtained as quotients or subquotients of those respective Vermas, and for all $i\geq0$, $$\Ext_{\Hlam}^i(\Sab, S_{\alpha',\beta'})=0$$

Combining Propositions 5.21, 5.22, and 5.25 and following \cite{SZ}, section 8, yields that the noncommutative space $\Mod\Hlam$ is the disjoint union of its fibers, where the fibers are the Verma modules:

\proposition $\operatorname{Mod}\Hlam$ is the disjoint union of the Verma modules $\Delta_{\alpha}$, $\alpha\in k$.\remark It was necessary to prove Propositions 5.21 and 5.22 to conclude that some nonzero Ext groups exist between modules in the same fiber, as we cannot apply Propositions 8.4 and 8.5 of \cite{SZ} since $k=\bar{\mathbb{F}_q}$ is not an uncountable field.


\subsection{Simple $\Hlam$-modules in the $t=1$ case} 
The description of simple left $\Hlam$-modules follows more or less the same pattern when $\lambda$ contains a constant term ($t=1$) as when $\lambda$ does not ($t=0$). The simple modules are the quotients of left Verma modules $\Da=\Hlam/(\Hlam(x-\alpha)+\Hlam(g-1))\cong k[y]$ by maximal submodules (which are necessarily cyclic). As in the $t=0$ case, there is a two-parameter family of simple modules over $\Hlam$.

The description of simple modules for $\Hlam$ when $t=1$ follows from Lemma 5.13.
\proposition Let $\lambda\in kG$ be a polynomial with constant term equal to $1$ and consider the $\Hlam$-modules $\Da=\Hlam/(\Hlam(x-\alpha)+\Hlam(g-1)$. For $\alpha\neq0$, the simple modules are $p^2$-dimensional. When $\alpha=0$ there are two cases:
\begin{itemize}
\item If $\lambda$ is a zero-divisor then $y-\beta$, $\beta\in k$, generates a maximal submodule of $\Delta_0$ and so the simple modules obtained from $\Delta_0$ are $1$-dimensional.
\item If $\lambda$ is invertible then $Y-\beta:=y^p-y\frac{\delta^p(g)}{\delta(g)}-\beta$ generates a maximal submodule and so the simple modules obtained from $\Delta_0$ are $p$-dimensional.
\end{itemize}
\begin{proof}
If $\alpha\neq0$ the subset of $g$-fixed elements of $\Da$ is $k[Y]$. 
$$x\cdot Y=Yx+[x,Y]=Y\alpha+[gyg^{-1}-y,Y]=Y\alpha-g(yx^p)g^{-1}+yx^p=Y\alpha-\alpha^{p+1}$$
and so if $\alpha\neq0$ the submodule generated by $Y$ contains $1$; likewise for $f(Y)$, deg $f<p$, by downwards induction. The submodules generated by the images of the central elements $y^{p^2}-y^p\frac{\delta^{p^2}(g)}{\delta(g)}$ are therefore the maximal ones.

When $\alpha=0$ then $g\cdot (y-\beta)=y-\beta$, $x\cdot (y-\beta)=c_{\lambda}$ where $c_{\lambda}$ is the sum of coefficients of $\lambda$; thus $y$ generates a proper submodule if and only if $\lambda$ is a zero-divisor. Otherwise, by the same argument as in the $t=0$ case, $f(y)$ does not generate a proper submodule for any $f$ of degree less than $p$. However, since \begin{align*} x\cdot (Y-\beta)&=Y\alpha-\alpha^{p+1}-\beta\alpha=0\\g\cdot (Y-\beta)&=(Y-\beta)g=Y-\beta,\end{align*} $Y-\beta$ generates a submodule, and it is maximal.
\end{proof}

\subsubsection{Ext groups between simples when $t=1$}
Propositions 5.21 and 5.22 carry over to the $t=1$ case and give the groups $\Ext^i (S)$ when $S$, a simple module, belongs to the smooth locus or when $S$ belongs to the singular locus and is $1$-dimensional, respectively. That takes care of the case $t=1$, $\lambda$ a zero-divisor. The only other case to consider when $t=1$ is $\lambda$ invertible, when $S$ belongs to the singular locus and is $p$-dimensional, smaller than the dimension of a simple module over the smooth locus. In this case, $C$ as it appears in the resolution in 5.21 goes by the name of $Y$ and doesn't commute with $x$ (rather, $[x,Y]=-x^{p+1}$); however, the resolution in 5.21 works for $S_{0,0}$ over $\Hlam$ if $C$ is replaced with $Y-x^p$, and it appears the resulting Ext groups are the same as in 5.21.

It follows that in the $t=1$ case as well as the $t=0$ case, $\Hlam$ has infinite global dimension and $\Mod\Hlam$ is the disjoint union of its fibers $\Da$.

\subsection{Other representations of $\Hlam$ worth noting: the differential operator representation, the polynomial representation, and the Weyl representation.}
There are several natural infinite-dimensional representations of $\Hlam$, beyond the most obvious one of $\Hlam$ as a module over itself. Because $\Hlam=\R[y;\delta]$, $\Hlam$ acts naturally on $\R=k[x]\otimes kG$ where $\R$ acts via left multiplication and $y$ acts via $\delta$. Furthermore, depending on whether $\lambda\notin \Rad kG$ or $\lambda\in \Rad kG$, then $\Hlam$ either has a representation on the Weyl algebra $\Weyl\cong\frac{k\la x,y\ra}{xy-yx-1}$ or on the polynomial algebra $k[x,y]$, respectively.
\subsubsection{The differential operator representation}
Consider the subring $R$ as a representation $H=R[y;\delta]$ on which elements $r\in R$ act as multiplication by $r$ while $y$ acts via $-\delta$:
$$y\cdot r=-\delta(r)$$
The minus sign is necessary since we write our Ore extension as a left rather than right extension. We will call $R$ with this $H$-action the differential operator representation of $H$. 
\lemma $R$ is not a faithful representation of $\Hlam$.
\begin{proof}
$\Hlam$ being finite over its center implies that $\Hlam$ has a central element of the form $Z=y^{p^n}-y^{m}F(x,g)$ for some $n$,$m$. Then $F(x,g)$ belongs to the kernel of $\delta$ since it commutes with $y$. Hence for any $r\in R$, 
$$Z\cdot r=Zr\cdot 1=rZ\cdot 1=-r(\delta^{p^n}(1)-\delta^{p^m}(F(x,g))=0$$
So the two-sided ideal $(Z)\subset\Hlam$ annihilates $R$.
\end{proof} 

The differential operator representation $R$ contains infinitely many infinite-dimensional subrepresentations in the ideals $R(x^{np})$, and also in the ideals $R(x^2-2g\del_g^{p-2}(\lambda))^n$ in the case $t=0$: since $\delta(x^p)=0$, and $\delta(x^2-2g\del_g^{p-2}(\lambda))=0$ when $t=0$, it holds for any $r$ in $R$ that \begin{align*}&y\cdot rx^{np}=-\delta(r)x^{np}\\&y\cdot r(x^2-2g\del_g^{p-2}(\lambda))^n=-\delta(r)(x^2-2g\del_g^{p-2}(\lambda))^n\end{align*}

\subsubsection{The Weyl representation}
For any $\lambda\notin\Rad kG$, the Weyl algebra is a module for $\Hlam.$ Let $c_\lambda\in k^{\times}$ be the sum of the coefficients of $\lambda;$ then $$A_{1,c_\lambda}=\frac{k\la x,y\ra}{(xy-yx-c_\lambda)}\cong\frac{k\la x,y\ra}{(xy-yx-1)}$$ is an $\Hlam$ module. The action is given by $y$ acting by left multiplication and $x$ and $g$ acting as they do on $\Hlam$ but then with the resulting polynomials in $x$ and $g$ evaluated at $g=1$:
\begin{align*}
&y\cdot y^nx^m=y^{n+1}x^m\\
&r\cdot y^nx^m=\sum_{i=0}^n{n\choose i}y^{n-i}x^m\delta^i(r)|_{g=1}\qquad\textrm{for all }r\in k[x]\otimes kG\\
&\qquad\qquad\qquad\qquad\qquad\qquad\qquad\qquad\qquad(\delta=\lambda\del_x+xg\del_g)
\end{align*}
Moreover, this module can be realized as the quotient of $\Hlam$ by the left ideal generated by $kG$: as a left $\Hlam$-module, $$A_{1,c_{\lambda}}\cong \Hlam/\Hlam(g-1).$$
Note that the isomorphism is an isomorphism of left $\Hlam$-modules only.

How does the Weyl algebra representation $A_{1,c_{\lambda}}$ relate to the Verma modules and simple modules we have examined earlier? By the Third Isomorphism Theorem, all of the Verma modules are subquotients of the Weyl representation: $$\Da=\Hlam/\left(\Hlam(g-1)+\Hlam(x-\alpha)\right)\cong A_{1,c_{\lambda}}/A_{1,c_{\lambda}}(x-\alpha)$$ Consequently every simple representation also arises as a subquotient of $A_{1,c_{\lambda}}$. 

The Ext groups of $A_{1,c_{\lambda}}$ with $\Hlam$ and with itself may be computed using the following free resolution, where the maps above the arrows denote right multiplication: $$...\xrightarrow{g-1}\Hlam\xrightarrow{\sum g^i}\Hlam\xrightarrow{g-1}\Hlam\xrightarrow{\sum g^i}\Hlam\onto A_{1,c_{\lambda}}\lra0$$ 
Applying $\Hom(-,\Hlam)$ we obtain:
$$0\lra\Hlam\xrightarrow{\sum g^i}\Hlam\xrightarrow{g-1}\Hlam\xrightarrow{\sum g^i}\Hlam\xrightarrow{g-1}...$$ 
from which we see that \begin{align*}&\Ext_{\Hlam}^0(A_{1,c_{\lambda}},\Hlam)=\Hlam(g-1)\\
&\Ext_{\Hlam}^i(A_{1,c_{\lambda}},\Hlam)=0\qquad\qquad\textrm{for all }i>0.\end{align*}
On the other hand, applying $\Hom(-,A_{1,c_{\lambda}})$ we obtain a complex in which every map is $0$ since right multiplication by $g-1$ or $\sum g^i$ annihilates $A_{1,c_{\lambda}}$. Therefore $$Ext_{\Hlam}^i(A_{1,c_{\lambda}})=A_{1,c_{\lambda}}\;\;\textrm{for all }i\geq0.$$
\corollary $A_{1,c_{\lambda}}$ has infinite projective dimension as an $\Hlam$-module.
\subsubsection{The polynomial representation} If $\lambda\in\Rad kG$ then the sum of the coefficients of $\lambda$ is zero mod $p$. The polynomial algebra $k[x,y]$ then becomes a representation of $\Hlam$ via the same action as in the Weyl representation above.
Moreover, $$k[x,y]\cong \Hlam/\Hlam(g-1)$$ and every Verma is a quotient of $k[x,y]$ as an $\Hlam$-module when $\lambda\in\Rad kG$, just as happened with the Weyl representation for $\lambda\notin\Rad kG$. Likewise the projective dimension of $k[x,y]$ as an $\Hlam$-module is infinite.

\section{The algebras $\Hlam$ for $E=G^r$}

Let $E=(\bb{Z}/p)^r$ be an elementary abelian $p$-group and realize $E$ as the unipotent subgroup of $\SL(2,\bb{F}_q)$, $q=p^r$: $$E=\big\{\begin{pmatrix} 1& a\\0&1\end{pmatrix} \;|\;a\in\bb{F}_q\big\}.$$ Let $\xi_1,\xi_2,...,\xi_r$ be a basis for $\bb{F}_q$ as an $\bb{F}_p$-vector space, so that $$g_i:=\begin{pmatrix}1&\xi_i\\0&1\end{pmatrix}$$ is the generator of the $i$th copy of $\bb{Z}/p$ and together $g_1,...,g_r$ generate $E$. We will always take $\xi_1=1$ so that $g_1=\begin{pmatrix}1&1\\0&1\end{pmatrix}$, $\xi_2=\xi$ where $\xi$ generates $\bb{F}_q^{\times}$ and  $\xi_i=\xi^{i-1}$.  Let $k=\bar{\bb{F}}_q$. The group algebra $kE$ is a truncated polynomial ring in the $r$ variables given by the generators $g_i$:
$$kE\cong k[g_1,g_2,...,g_r]/(g_i^p-1)_{i=1}^r$$
Let $$\gdelg=g_1\del_{g_1}+\xi g_2\del_{g_2}+\xi^2g_3\del_{g_3}...+\xi^{r-1}g_r\del_{g_r}.$$ By $\del_{g_i}$ we mean partial differentiation with respect to the variable $g_i$, thinking of the elements of $kE$ as polynomials in $r$ variables of degree at most $p-1$ in each variable. For any $\lambda\in kE$, the algebra $$\Hlam:=\frac{k\la x,y\ra \rtimes E}{(xy-yx-\lambda)}$$ admits a presentation as a differential operator Ore extension over its commutative subring $$\R:=k[x]\otimes kE$$ as follows: $$\Hlam=\R[y;\delta],\qquad \delta:=\lambda\del_x+x\gdelg.$$ This means that for any $r\in\msf{R},$ $$ry=yr+\delta(r).$$ As with the algebras $\Hlam$ for $G$ the cyclic group of order $p$ studied so far in this paper, the family of algebras $\Hlam$ for $E$ consists of PBW deformations of the skew group ring $\Hzero=k[x,y]\rtimes E$ analogous to certain symplectic reflection algebras, namely those deformations of $\bb{C}[x,y]\rtimes G$ where $G$ is a finite subgroup of $\SL(2,\bb{C})$ and where $[x,y]$ takes a value $\mathsf{c}$ in the center of the group algebra of $G$. 

As we have mentioned earlier in this paper, the family of symplectic reflection algebras associated to a vector space $V$ and finite group $G$ (when such a family exists) has two branches: the side of the family where $\mathsf{c}$ contains a constant term (``$t=1$") and the side of the family where it does not (``$t=0$"). The two sides of this family behave differently in their algebraic structure and representation theory. In our work so far we have seen that the same dichotomy holds for the algebras $\Hlam$ and that the action of the operator $\gdelg$ is key: the two cases $t=0$ and $t=1$ correspond to whether or not $\lambda$ lies in the image of $\gdelg$, that is, whether or not $(g\del_g)^{q-1}(\lambda)=\lambda$, and the behavior of the operator $\gdelg$ is the underlying reason for the different algebra structures that appear in the two cases. If we write $\bar{a}=(a_1,...,a_r)\in\E$ and $$\lambda=\sum_{\bar{a}}c_{\bar{a}}g_1^{a_1}g_2^{a_2}\cdot\cdot\cdot g_r^{a_r}$$ then the condition $(\gdelg)^{q-1}(\lambda)=\lambda$ is equivalent to requiring that we have $c_{\bar{a}}=0$ whenever $a_1+\xi_2a_2+...+\xi_ra_r=0\;\textrm{ mod } p$; since the $\xi_i$ are linearly independent over $\bb{F}_p$ this is no more than the requirement that the coefficient of $1$ is $0$. Thus we expect that the algebras $\Hlam$ constructed from $\E$ will likewise fall into two families, determined by whether $c_0=0$ or $c_0\neq0$.

We will exclusively work with the $t=0$ case from now on. Thus $\lambda\in kE$ is a nonzero polynomial in $g_i$, $i=1,..,r$, of degree at most $p-1$ in each variable and without a constant term. 

\subsection{Some remarks about the center of $\Hlam$.}

Let $\Zlam$ denote the center of $\Hlam$. We expect that $\msf{Z}_{\lambda}$ is isomorphic to the quotient of a polynomial ring on three variables by one relation (between the first and second variables), and that $\Hlam$ is of generic rank $qp$ as a module over $\msf{Z}_{\lambda}$. We first identify those generators of $\msf{Z}_{\lambda}$ belonging to the commutative subalgebra $\R=k[x]\otimes kE$. To check an element $r\in \R$ commutes with $y$ it is sufficient to check that $\delta(r)=0$, since $ry=yr+\delta(r)$ for all $r\in \R$ by the definition of a differential operator extension.
\lemma $x^p\in\msf{Z}_{\lambda}$.
\begin{proof}
$\delta(x^p)=px^{p-1}\delta(x)=0$.
\end{proof}
Since $\lambda$ does not have a constant term, the operator $\gdelg$ has an inverse on $\lambda$: $(\gdelg)^{q-1}(\lambda)=\lambda$ and we write $(\gdelg)^{-1}(\lambda)$ for $(\gdelg)^{q-2}(\lambda)$. This allows us to identify a second, quadratic central element of $\Hlam$ contained in the subalgebra $\R$:
\lemma $x^2-2(\gdelg)^{-1}(\lambda)\in\msf{Z}_{\lambda}$.
\begin{proof}
$\delta(x^2-2(\gdelg)^{-1}(\lambda))=2x\lambda-2x\gdelg((\gdelg)^{-1}(\lambda))=2x\lambda-2x\lambda=0$.
\end{proof}
Let $$A=x^2-2(\gdelg)^{-1}(\lambda)+2\sum_{\bar{a}}(a_1+\xi_2a_2+...+\xi_na_n)^{-1}c_{\bar{a}}$$$$B=x^p$$ Between $A$ and $B$ there is the obvious relation: $$A^p=B^2$$

Recall that $\Hlam$ is filtered by degree in $x$ and $y$ with both variables assigned degree $1$ and $E$ put in degree $0$. Then the associated graded algebra with respect to this filtration is $$\gr\Hlam=\Hzero=k[x,y]\rtimes E.$$ The center $\msf{Z}_0$ of $\Hzero$ is the subalgebra of $E$-invariants in $k[x,y]$: $$\msf{Z}_0=k[x,y]^E=k[x,y^q-yx^{q-1}]$$ There is a natural homomorphism from $\Zlam$ to $\Zzero$ but it is not necessarily onto. We wish to find a central element which deforms the second generator $y^q-yx^{q-1}$ of $\Zzero$, that is, a central element whose top degree term is $y^q-yx^{q-1}$. A computation shows:

\lemma If $Z\in\Hlam$ commutes with $y$ and $Z$ is of the form $y^q-yx^{q-1}+\textit{lower order terms}$ then $Z$ is of the form $$Z=y^q-y(x^{q-1}+x^{q-3}\gamma_3+x^{q-5}\gamma_5+...)$$ where  $\gamma_{i}\in kE$ satisfy $\gdelg(\gamma_{i})=(i-2)\gamma_{i-2}$.

\lemma If $C\in\Hlam$ commutes with $y$ and $g_1$ and $C=y^q-yx^{q-1}+\textit{lower order terms}$ then this uniquely determines $C$: $$C=y^q-y\frac{\delta^q(g_1)}{\delta(g_1)}$$

\begin{proof}
By the previous lemma, $C=y^q-yF$, $F\in k[x]\otimes kG$. Then $$g_1C=g_1(y^q-yF)=y^qg_1+\delta^q(g_1)-yFg_1-\delta(g_1)F=(y^q-yF)g_1=Cg_1$$ so $\delta^q(g_1)-\delta(g_1)F=0.$ Note that $\delta(g_1)=xg_1$ is not a zero divisor and it divides $\delta^q(g_1)$; thus we can solve for $F$.
\end{proof}

Returning to the the degree $q$ generator of $\Zzero$, it may be expressed as: $$y^q-yx^{q-1}=\prod_{a\in\bb{F}_q}(y-ax)$$
A natural question is whether the expression $\prod_{a\in\bb{F}_q}(y-ax)$ considered as an element of $\Hlam$ gives a central element of $\Hlam$ if the factors are put in the right order. In general this fails but in one special case the answer is yes: when $\lambda$ is an element of the group. 

\subsection{The center $\msf{Z_g}$ of $\msf{H_g}$}
Consider the case that $[x,y]=g$ where $g\in E$ is an honest element of the group different from the identity. All of the $q-1$ algebras obtained this way are isomorphic, as the following lemma shows.

\lemma For any $\alpha\in\bb{F}_q^{\times}$, set $g_{\alpha}=\begin{pmatrix} 1&\alpha\\0&1\end{pmatrix}$ and $$\msf{H_{g_{\alpha}}}:=\frac{k\la x,y\ra\rtimes E}{(xy-yx-g_{\alpha})}.$$ Then $$\Hg\cong\msf{H_{g_{\alpha}}}.$$ 

\begin{proof} Define an algebra homomorphism $$\phi_{\alpha}:\Hg\to\msf{H_{g_{\alpha}}}$$ such that \begin{align*} \phi_{\alpha}\begin{pmatrix}1&\beta\\0&1\end{pmatrix}&=\begin{pmatrix}1&\alpha\beta\\0&1\end{pmatrix}\mbox{ for any }\beta\in\bb{F}_q\\ \phi_{\alpha}(x)&=\sqrt{\alpha}x\\ \phi_{\alpha}(y)&=\frac{1}{\sqrt{\alpha}}y\end{align*}
Then $\phi_{\alpha}$ takes generators to generators and preserves all defining relations.
\end{proof}

\proposition Take $q=p$ and $\lambda=g:=\begin{pmatrix}1&1\\0&1\end{pmatrix}$. Then the third generator of $\msf{Z_{g}}$ is  $$C:=\prod_{a=1}^p(y-ax)$$
{\textit{Note that the order in which the product is taken matters; let indices going from bottom to top correspond to multiplicands going from left to right.}
\begin{proof}
Using the observation that $[y-ax,y]=-ag$ together with the Leibniz rule, we compute the commutator of $C$ and $y$:
\begin{align*}
[C,y]&=\sum_{a\in\bb{F}_p}\left(-a\prod_{b=1}^{j-1}(y-bx)\prod_{c=j+1}^p(y-(c-1)x)g\right)\\
&=\sum_{a\in\bb{F}_p}\left(-a\prod_{b=1}^{p-1}(y-bx)g\right)\\&=\left(\sum_{a\in\bb{F}_p}a\right)\prod_{b=1}^{p-1}(y-bx)g\\&=0
\end{align*}
since $\sum_{a\in\bb{F}_p}a=0$.

As for the commutator with $g$, use that $[g,y]=xg$, first move the $g$ past the terms to its right, then move the $x$ past the resulting terms using that the lower-order cost of moving $x$ past each term in parenthesis is to replace that term by a $g$ :
\begin{align*}
[g,C]&=\sum_{a\in\bb{F}_p}\left(\prod_{b=1}^{a-1}(y-bx)x\prod_{c=a}^{p-1}(y-cx)g\right)\\&=p\prod_{a=1}^{p-1}(y-ax)xg+\sum_{a=1}^{p-1}\sum_{b=a}^{p-1}\left(\prod_{c=1}^{b-1}(y-cx)g\prod_{d=b+1}^{p-1}(y-dx)g\right)\\&=0+\left(\sum_{a=1}^{p-1}a\right)\prod_{b=1}^{p-2}(y-bx)g^2\\&=0
\end{align*}
\end{proof}

\corollary Take $q=p$ and $\lambda=g_n:=\begin{pmatrix}1&n\\0&1\end{pmatrix}$ for any $1<n<p$. Then the third generator of $\msf{Z_{g_n}}$ is  $$C:=\prod_{a=1}^p(y-nax)$$\\

An analogous statement holds for $r>1.$ From now on, we take $r>1, \;q=p^r$, and $E=(\bb{Z}/p)^r$; $1, \xi,\xi^2,...,\xi^{r-1}$ is an $\bb{F}_p$ basis for $\bb{F}_q$ and we set $$\Xi:=\left<\xi,\xi^2,...,\xi^{r-1}\right>_{\bb{F}_p}$$ to be the $\bb{F}_p$-linear span of the $\xi_i$, $i>1$. Form $\Hg:=\Hlam$ with $\lambda:=g_1=[x,y]$ where $g_1=\begin{pmatrix}1&1\\0&1\end{pmatrix}$. 
\proposition  The element $$C=\prod_{\zeta\in\Xi}\left(\prod_{a=1}^p(y-(\zeta+a)x)\right)$$ in $\Hg$ commutes with $y$ and $g_1$. 

\textit{N.b. the first product does not require any particular order but the innermost product over $a$ requires the order specified, up to a cyclic permutation.}
\begin{proof}
For a given $\zeta$, the terms involving $\zeta$ appear in succession with $a$'s increasing by $1$ from left to right. When we apply the Leibniz rule to $C$ to compute $[C,y]$, it has the effect of picking out one term $(y-(\zeta+a)x)$ at a time and replacing it with $g_1$ which then moves past the remaining terms to its right decreasing the number $a$ in each of them by $1$ (since the terms involve subtracting $ax$). This does not affect any terms with $\zeta'\in\Xi$ which appeared to the left of $\zeta$; moreover, the effect on all terms $\zeta''\in\Xi$ which appear to the right of $\zeta$ is the same for each for value of $a=1,...,p$. There are $p$ terms with $\zeta$, and (as we saw in the proof of the cyclic case) deleting one term and replacing it with $g_1$ always gives $$(y-(\zeta+1)x)(y-(\zeta+2)x\cdot\cdot\cdot(y-(\zeta+(p-1)x))$$ once $g_1$ is moved all the way to the right of the terms with $\zeta$. Thus there are $p$ identical summands for each $\zeta\in\Xi$, and so the sum is $0$.

Similarly, applying the Leibniz rule to compute $[g_1,C]$ has the effect of picking out one term $(y-(\zeta+a)x)$ at a time and replacing it with $xg_1$. By the same argument, this reduces to the cyclic case; for a fixed $\zeta\in\Xi$ we get $$\left(\prod_{\zeta'<\zeta}\prod_a(y-(\zeta'+a)x)\right)\left(\sum_{b\in\bb{F}_p}b\right)\left(\prod_{a=1}^{p-2}(y-(\zeta+a)x)\right)xg_1\left(\prod_{a\zeta''>\zeta}\prod_b(y-(\zeta''+a)x)\right)=0 $$ because $\sum_{b\in\bb{F}_p}b=0$. 
\end{proof}

\corollary When $\lambda=g_1$, $$C:=\prod_{\zeta\in\Xi}\prod_{b=1}^p(y-(\zeta+a)x)=y^q-y\frac{\delta^q(g_1)}{\delta(g_1)}$$ and $C$ either is the third generator of the center $\Zg$, or there is no such third generator of degree $q$.

\begin{proof}
We have seen that $C$ commutes with $y$ and $g_1$; also, the top degree term of $C$ is $y^q-yx^{q-1}$. By Lemma 6.4 then, $$C=y^q-y\frac{\delta^q(g_1)}{\delta(g_1)}.$$ On the other hand, suppose a generator $\tilde{C}$ of the center $\Zg$ exists with $\tilde{C}$ of degree $q$. Then $\tilde{C}$ maps to $y^q-yx^{q-1}$ in the associated graded map $\gr:\Hg\to\Hzero$ (the map chops off all lower order terms). Since $\tilde{C}$ is central, it commutes with $y$ and $g_1$. It follows from Lemma 6.4 that $C=\tilde{C}$.
\end{proof}

Our aim now is to show that $C$ is indeed central, and so we must show that $C$ commutes with the other generators $g_i$ of $E$. 

\example Let $p=3$ and $r=2$, and let $$C=\prod_{a=1}^3\prod_{b=1}^3(y-(a\xi+b)x)$$
where $1,\xi$ are an $\bb{F}_3$-basis for $\bb{F}_9$ over $\bb{F}_3$. Then $C\in \msf{Z}_g$.
\begin{proof}
Fix $a$. Then a direct calculation shows: $$\prod_{b=1}^3(y-(a\xi+b)x)=y^3-yx^2-yg-\prod_{b=1}^3(a\xi+b)x^3$$ Set $Z=y^3-yx^2-yg$ and note that $Z=y^3-yA$ where $A=x^2-2\gdelg(g)$ is the quadratic central element of $\Hg$. Since $x^3\in\msf{Z}_g$ we see that each of these $3$ subproducts of $C$ commute with each other: $$C=\left(Z-\prod_{b=1}^3(\xi+b)x^3\right)\left(Z-\prod_{b=1}^3(2\xi+b)x^3\right)Z$$ 
The effect of multiplication by $g_2$ is to cyclically permute these three factors. Since they commute, this preserves $C$:
$$g_2C=\prod_{a\in\bb{F}_p}\prod_{b=1}^p(y-((a-1)\xi+b)x)g_2=Z\left(Z-\prod_{b=1}^3(\xi+b)x^3\right)\left(Z-\prod_{b=1}^3(2\xi+b)x^3\right)g_2=Cg_2$$
\end{proof}

The example gives hope that the $p^{r-1}$ products $$Y_{\zeta}:=\prod_{a=1}^p(y-(\zeta+a)x)$$ will commute with one another for any $p$. Set $$Y:=\prod_{a=1}^p\left(y-ax\right)$$
Since $[x,y]=g_1$ this is simply the same element which was central in the $\bb{F}_p$ case: $$Y=y^p-\frac{\delta^p(g)}{\delta(g)}=y^p-y\left(x^2-2g\right)^{\frac{p-1}{2}}.$$
On the one hand, the effect of $g_i$ on $C$ is to permute the $Y_{\zeta}$; for any $g_{\zeta}=\begin{pmatrix}1&\zeta\\0&1\end{pmatrix}\in E$, $$g_{\zeta}\cdot Y=Y_{\zeta}.$$ On the other hand,  $$g_{\zeta}\cdot Y=g_{\zeta}\cdot\left(y^p-y\frac{\delta^p(g)}{\delta(g)}\right)=(y+\zeta x)^p-(y+\zeta x)\left(x^2-2g\right)^{\frac{p-1}{2}}$$ and thus $$Y_{\zeta}=\prod_{a=1}^p\left(y-(a-\zeta)x\right)=(y+\zeta x)^p-(y+\zeta x)A^{\frac{p-1}{2}}$$ 

\proposition $$Y_{\zeta}=Y-\prod_{a\in\bb{F}_p}\left(a-\zeta\right)x^p=Y+(\zeta^p-\zeta)x^p$$
\begin{proof}
Let $A_{k,j}$ be the coefficient of $x^jg^k$ in $\delta^{n-2k}(g)$ in the algebra $\Hg$ when $E=G=<g>$. We claim that the $\bb{F}_q$-coefficient of $y^{n-2k-j}x^jg_1^k$ in $(y+\zeta x)^n$ is given by ${n\choose 2k+j}T_{k,j}$ where $T_{k,j}$ is defined by the linear recursive formula

 $$T_{0,0}=1$$
 
 $$T_{k,j}=(\zeta+k)T_{k,j-1}+(j+1)T_{k-1,j+1}$$

$T_{k,j}$ is assumed to be $0$ if $k$ or $j$ is negative. The coefficients $T_{k,j}$ form a square array that's a variation on the Andr\'{e} numbers $A_{k,j}$, the coefficients of Andr\'{e} polynomials (cf. Section 4.3). The formula is easily checked by writing $(\zeta x+y)^n=(\zeta x+y)(\zeta x+y)^{n-1}$ and expanding out $(\zeta x+y)^{n-1}$ according to the formula by induction. Thus $Y_{\zeta}$ will have the desired form if and only if:

 \begin{enumerate}
 
 \item $T_{0,p}$, the coefficient of $x^p$, is $\zeta^p$\\
 
  \item for each $k\geq1$, $T_{k,p-2k}=\zeta A_{k,p-2j}$ where $A_{k,j}$ is the $k,j$th Andr\'{e} number (mod $p$).\\

\end{enumerate}
The requirement (1) is obvious: for any $j\geq0$, it follows immediately from the recursive formula that $$T_{0,j}=\zeta^j$$
We prove (2) by induction. First, $k=1$: by induction on $j$, $$T_{1,j}={j+2\choose2}\zeta^{j+1}+{j+2\choose3}\zeta^j+...+{j+2\choose j+1}\zeta^2+\zeta$$
Therefore, $$T_{1,p-2}={p\choose2}\zeta^{p-1}+...+{p\choose p-1}\zeta^2+\zeta=\zeta\qquad\mbox{(mod } p)$$
Since $A_{1,p-2}=1$ mod $p$ this establishes (2) for $k=1$.
For the induction step, assume: $$T_{k,p-2k}=\zeta A_{k,p-2k}.$$ Then consider the element in the sequence whose incoming arrows come from $T_{k,p-2k}$ and $T_{k+1,p-2k-2}$, namely $T_{k+1,p-2k-1}$: 
$$T_{k+1,p-2k-1}=(\zeta+k+1)T_{k+1,p-2k-2}-2kT_{k,p-2k}=(\zeta+k+1)T_{k+1,p-2k-2}-2k\zeta T_{k,p-2k}$$

The next ingredient is a quadratic recursive formula for $T_{k,j}$ which can be proven inductively by writing $(\zeta x+y)^n=(\zeta x+y)^{n-1}(\zeta x+y)$ and expanding $(\zeta x+y)^{n-1}$ according to the quadratic recursive formula below by induction:

$$T_{k,j}=\zeta T_{k,j-1}+\zeta\sum_{l=0}^{k-1}\sum_{m=0}^j{2k+j-1\choose 2l+m}T_{l,m}A_{k-1-l,j-m}$$

This formula takes a simple form when $2k+j-1=p$ since all but one of the binomial coefficients in the sum vanishes:

$$T_{k,p-2k+1}=\zeta T_{k,p-2k-2}+\zeta A_{k-1,p-2k+1}$$

(since $T_{0,0}=1$) and so replacing $k$ with $k+1$, we have:

$$T_{k+1,p-2k-1}=\zeta T_{k+1,p-2k-2}+\zeta A_{k,p-2k}$$

Setting the two expressions for $T_{k+1,p-2k-1}$ equal to each other, we have:

$$(\zeta+k+1)T_{k+1,p-2k-2}-2k\zeta T_{k,p-2k}=\zeta T_{k+1,p-2k-2}+\zeta A_{k,p-2k}$$

Solving for $T_{k+1,p-2k-2}$, we get:

$$T_{k+1,p-2k-2}=\zeta\frac{2k+1}{k+1}A_{k,p-2k}=\zeta A_{k+1,p-2k-2}$$
 
 The last equality uses the fact mentioned in Section 4.3 on the center of $\Hg$ for $E=G=<g>$ that $A_{k,p-2k}=\frac{(2k-1)!!}{k!}$ mod $p$ for each $A_{k,p-2k}$ along the $k$th knight's-move antidiagonal of the sequence. That numerical formula can be proven on its own by comparing the two recursive definitions of $A_{k,j}$ and using induction on $k$, although it also follows, after reindexing, from the more general formula that appeared in the proof of Theorem 4.2, with $\lambda=g$: $$\gdelg(C_{p,k})=-(k+1)\lambda C_{p,k+2}.$$
\end{proof}

\theorem The element $$C=\prod_{\zeta\in\Xi}\prod_{a=1}^p(y-(\zeta+a)x)=\prod_{\zeta\in\Xi}Y_{\zeta}$$ belongs to the center $\Zg$ of $\Hg$. Thus $$\Zg\cong \frac{k[A,B,C]}{(A^p-B^2)}$$
\begin{proof}
Proposition 2.8 says that $C$ commutes with $y$ and $g_1$, while Proposition 2.11 implies that $C$ commutes with $g_i=\begin{pmatrix}1&\xi^{i-1}\\0&1\end{pmatrix}$: $$g_i C=g_i\prod_{\zeta\in\Xi}Y_{\zeta}=\left(\prod_{\zeta\in\Xi}Y_{\xi^{i-1}+\zeta}\right)g_i=\left(\prod_{\zeta\in\Xi}Y_{\zeta}\right)g_i=Cg_i$$ Since $y$ and the $g_i$ generate $\Hg$, that proves the first statement. For the second, Lemmas 2.1 and 2.2 give the generators $A$ and $B$ and their relation; it is easy to argue $A$,$B$, and $C$ generate the entire center. \end{proof}
 
Applying the isomorphism of Lemma 6.5 gives an analogous formula for the other algebras where the bracket is group-valued. The effect of the isomorphism is to permute the linear factors. Let $\alpha$ be any nonzero element of $\bb{F}_q$ and form $\msf{H_{g_{\alpha}}}$ where $[x,y]=g_{\alpha}$. The following corollary explicitly desribes the center $\msf{Z_{g_{\alpha}}}$ of $\msf{H_{g_{\alpha}}}$.
 
\corollary $$\msf{Z_{g_{\alpha}}}\cong\frac{k[X,Y,Z]}{(X^p-Y^2)}$$ where  generators for $\msf{Z_{g_{\alpha}}}$ may be taken as follows: \begin{align*}
A&=x^2-2\gdelg^{-1}(g_{\alpha})\\
B&=x^p\\
C&=\prod_{\zeta\in\Xi}\left(\prod_{a=1}^p(y-\alpha(\zeta+a)x)\right)
\end{align*}\\

\subsection {Simple representations over $\Hg$ when $E=(\bb{Z}/p)^r$, $r>1$}

Simple modules over $\Hg$ correspond to central characters and thus to points in $\Spec \Zg$. The most interesting question concerns what happens over the singular locus of the center. $\Spec\Zg$ is singular along the line $A=B=0$. First, a Verma module over $0$ can be defined as $$\Delta_0:=\Hg/\left(\Hg x+\sum_{i=1}^r\Hg(g_i-1)\right)$$ with the expectation that its simple quotients will coincide with the simple $\Hg$-modules on which $x^p$ acts by $0$. These are exactly the modules corresponding to the singular line in $\Spec\Zg$. $\Delta_0$ is infinite-dimensional and has $k$-basis $1,y,y^2,....$. 

\proposition For each $\beta\in k$, $Y-\beta=y^p-y\left(x^2-2g\right)^{\frac{p-1}{2}}-\beta$ generates a maximal left ideal of $\Delta_0$. Therefore the simple modules over the singular locus are $p$-dimensional.

\begin{proof}
First, the action of $y$ can only raise the degree of anything in $\Delta_0$. Second, any $g\in E$ acts trivially on $Y-\beta$:  $g(Y-\beta)=(Y-\beta)g+cx^pg$ for some $c\in k$, and $g=1$ and $x^p=0$ in $\Delta_0$, so $gY=Y$ in $\Delta_0$. Moreover, $x$ acts by $0$ on $Y-\beta$ since $x=g_1yg_1^{-1}-y$ commutes with $Y-\beta$. By Lemma 5.13, $Y-\beta$ then generates a proper submodule.

As for maximality: it holds for $G$ cyclic that $x,y,$ and $g_1$ do not preserve any larger submodule of $\Delta_0$ so a fortiori it holds for $E$.

\end{proof}
Over the smooth locus of the center, simple modules are $q$-dimensional. They are quotients of Verma modules $\Da$ whose simple quotients give the simple modules on which $x^p$ acts by $\alpha^p\in k^{\times}$: $$\Da:=\Hg/\left(\Hg(x-\alpha)+\sum_{i=1}^r(g_i-1)\right)$$

\proposition $C-\beta=\prod Y_{\eta}-\beta$ generates a maximal submodule of $\Da$. Thus the dimension of the simple modules over the smooth locus is $q$.
\begin{proof}
By repeating the arguments in the case $E=G$ it is clear that a maximal submodule $M\subset\Da$ must be cyclically generated by an element of degree a power of $p$; moreover for $g_1$ to fix the generator it must be of the form $$f=(Y-\beta_1)(Y-\beta_2)\cdot\cdot\cdot(Y-\beta_d)-\beta$$ 
Thus the degree of $f$ is $dp$ and it may be supposed this is the minimal degree of any element of $M$. Then if $\eta\notin\bb{F}_p$, $\eta^p-\eta\neq0$ and 
$$g_{\eta}\cdot f=(Y-\beta_1-(\eta^p-\eta)\alpha^p)(Y-\beta_2-(\eta^p-\eta)\alpha^p)\cdot\cdot\cdot(Y-\beta_d-(\eta^p-\eta)\alpha^p)$$ ($x^p$ may be replaced with its image $\alpha^p$ in $\Da$ because it is central) and either $g_{\eta}\cdot f=f$ in $\Da$ or $f-g_{\eta}\cdot f\in M$ is a polynomial of smaller degree, contradicting minimality. But the former is only possible if for every $i$, $\beta_i=\beta_j-(\eta^p-\eta)\alpha^p$ for some $j$, and as this is true for every $\eta$, $d=p^{r-1}$ and the $\beta_i$ must range over all $\eta$ in the $\bb{F}_p$-span of $\xi^i$, $i=1,...,r-1$. It follows that $$f=\prod_{\eta} Y_{\eta}-\beta=C-\beta.$$


\end{proof}

To summarize, simple modules over $\Hg$ are parametrized by pairs $\{(\alpha,\beta)\}\in k^2$. The simple modules over the smooth locus of the center correspond to pairs with $\alpha\neq0$; in this case, the simple modules are $q$-dimensional. When $\alpha=0$ the corresponding simple modules lie above the singular locus of the center and have dimension $p$. Thus when $q=p^r$ with $r>1$ the Azumaya locus and the smooth locus coincide for $\Hg$, unlike what happens for the cyclic group $G$.

\subsection{The center of $\Hlam$ when $\lambda\in kE$}

The general solution is less involved than the special case of $\lambda=g$, building off the results of Section 4; the formulation of the big central element, while precise, is less explicit than the factorization into $q$ linear factors given for the big generator of $\Zg$:

\theorem Let $E=\E\subset\SL(2,\bb{F}_q)$, $q=p^r$, $g_i$ the generator of the $i$th copy of $G$, $\lambda=\sum\limits_{i=1}^r\sum\limits_{j=0}^{p-1} c_{i,j}g_i^j\in kE$, and $$\Hlam=\frac{k\la x,y\ra\rtimes E}{xy-yx-\lambda}.$$  Scale $\lambda$ so that $c_{0,0}$, the coefficient of $1$ in $\lambda$, is either $0$ or $1$, and let $\Zlam$ denote the center of $\Hlam$. Then:

\begin{itemize}
\item (\textit{$t=0$ case}). When $c_{0,0}=0$ then $\Zlam$ is generated by three elements which have degrees $2,\;p,\mbox{ and }q$:
\begin{align*}
&x^2-2D^{q-1}(\lambda),\\
&x^p,\\
&y^q-y\frac{\delta^q(g_1)}{\delta(g_1)}
\end{align*}
Therefore $$\Zlam\cong \frac{k[A,B,C]}{(A^p-C^2)}$$
\item(\textit{$t=1$ case}). When $c_{0,0}=1$ then $\Zlam$ is generated by two elements which have degrees $p$ and $pq$:
\begin{align*}
&x^p,\\
&\prod_{a\in \bb{F}_p}(C-a x^q)
\end{align*}
where $C$ is the element that was central in the $t=0$ case. Therefore $$\Zlam\cong k[B,D]$$ 
\end{itemize}
\begin{proof}
This follows quite easily from the combinatorics of $\delta$ used in the proof of Theorem 4.2 (the case of $E=G=\bb{Z}/p$. As at the beginning of this section, take $\xi$ any generator of $\bb{F}_q^{\times}$, $\xi_i=\xi^{i-1}$ a basis for $\bb{F}_q$ over $\bb{F}_p$, and $g_i=\begin{pmatrix}1&\xi_i\\0&1\end{pmatrix}$. Set $$D:=\sum_{i=1}^r\xi_ig_i\del_i.$$ Then the $r$ trees for $\delta^n(g_i)$ all look almost identical to that for $\delta^n(g)$ in the cyclic case (see the diagrams preceding Lemma 4.1) when polynomials in powers $D^m(\lambda)$ are translated into partitions and $g_i$ is not written in the tree, except that each partition is ``homogenized" by some power of $\xi_i$, so that if $D^j(\lambda)$ has degree $j+1$ and $\xi_i$ has degree $1$ then every partition in $\delta^n(g_i)$ has degree $n$ as a polynomial in $\xi_i,\lambda,D(\lambda),...,D^{q-2}(\lambda)$ for $n$ up to $q$. As for $\delta^n(\lambda)$, its tree will be identical to what it was in the cyclic case, so long as the sizes of parts in the partitions are not reduced mod $p-1$ but are allowed to grow.



Following the reasoning of the proof of Theorem 4.2, the coefficients $C_{p^j,k}$ at all but the rightmost node $C_{p^j,p^j}$ in $\delta^{p^j}(g_i)$, the $p^j$th row of the tree for $\delta^n(g_i)$, which survive mod $p$ belong to those partitions whose parts sum to the maximal possible number: $p^j-1$. That means all $C_{p^j,k}$, $k<p$ are divisible by the same power of $\xi_i$, namely $\xi_i$, and so as in the proof of Theorem 4.2, \\$D(C_{p^j,k})=-(k+1)C_{p^j,k+2}$ for $k$ up to $p^{j-4}$. It is easy to see that $C_{p^j,p^j-2}=\xi_iD^{p-2}(\lambda)$ and $C_{p^j,p^j}=\xi_i^{p^j}$ and thus
 $$\delta\left(\frac{\delta^{p^j}(g_i)}{\delta(g_i)}\right)=\left(D^{p^j-1}(\lambda)-\xi_i^{p^j-1}\lambda\right)x^{p^j}$$

It follows that in the $t=0$ case, 
$$\delta\left(\frac{\delta^q(g_i)}{\delta(g_i)}\right)=\left(D^{q-1}(\lambda)-\lambda\right)x^q=0$$

while in the $t=1$ case, 

$$\delta\left(\frac{\delta^q(g_i)}{\delta(g_i)}\right)=\left(D^{q-1}(\lambda)-\lambda\right)x^q=-x^q$$

And so in the $t=0$ case, $C=y^q-y\frac{\delta^q(g_1)}{\delta(g_1)}$ is central. In the $t=1$ case, $C,\;x^q,\;\mbox{and }y$ generate a subalgebra isomorphic to $k[\tilde{x},\tilde{y}]\rtimes \la \tilde{g}\ra$ with $\tilde{y}=C$, $\tilde{x}=x^q$, and $\tilde{g}=y$, and thus $D=\prod\limits_{a\in\bb{F}_p}(C-a x^q)$ is central.
\end{proof}

\subsection{Verma modules and simple modules for $\Hlam$, $\lambda\in kE$}
Form the left Verma modules $$\Da=\Hlam\left/\left(\Hlam(x-\alpha)+\sum_{i=1}^r\Hlam(g_i-1)\right)\right.$$ for each $\alpha\in k$. The spectrum of the center of $\Hlam$ is singular over $x=0$ when $\lambda$ has no constant term, while in the $t=1$ case $\Zlam$ is smooth everywhere since then $\Zlam$ is a polynomial ring. Nonetheless even in the $t=1$ case there is a difference between the sizes of simple quotients of $\Delta_0$ and those of $\Da$ for $\alpha\neq0$. 

Recall that $\Da$ is the fiber over $x=\alpha$ of the projection $\Mod\Hlam\onto\Spec\R$, $\R=k[x]\otimes kE$. The $0$-fiber $\Delta_0$ coincides with the singular locus when $t=0$. $\Da$ is cyclic and isomorphic to $k[y]$ as a $k$-vector space. We look for maximal submodules of $\Da$. Clearly the big central element generates a submodule, and this bounds the degrees of simples by the degree of the big central element: $q$ when $t=0$ and $pq$ when $t=1$. These upper bounds are realized over the generic fiber $\Da$.

\theorem Consider $\Hlam$ for $E$. In both the $t=0$ and $t=1$ cases, simple $\Hlam$-modules over the $0$-fiber are $p$-dimensional if $\lambda\notin\Rad kE$, and $1$-dimensional if $\lambda\in\Rad kE$. Simple modules over the generic fiber (that is, simple quotients of $\Da$ for $\alpha\neq0$) have dimension $q$ when $t=0$ and dimension $pq$ when $t=1$.

\begin{proof}

Consider $\Delta_0$ when $\lambda\in\Rad kE$. Then for any $\beta\in k$, $g_i$ fixes $y-\beta$ for each $i$ and $x$ annihilates $y-\beta$ so $y-\beta$ generates a maximal submodule of $\Delta_0$; this is true in both the $t=0$ and $t=1$ cases.

By the arguments we've made previously in this type of theorem, a maximal submodule needs to be generated by an element of degree a power of $p$. Set $$Y=y^p-y\frac{\delta^p(g_1)}{\delta(g_1)}$$
By Section 6.4, $$[y,Y]=y(D^{p-1}(\lambda)-\lambda)x^p$$
and thus \begin{align*}
x\cdot Y=Yx+[x,Y]&=Y\alpha+[g_1yg_1^{-1}-y,Y]\\
&=Y\alpha-\alpha^{p+1}(D^{p-1}(\lambda)-\lambda)|_{g_i=1}
\end{align*}
so that $x\cdot Y$ is a scalar multiple of $Y$ if $\alpha=0$. Furthermore, observe that 
$$\delta^p(g_i)-\xi_i^pg_ix^p=\xi_i\left(\delta^p(g_1)-x^pg_1\right)g_1^{-1}g_i$$
(this follows from thinking about the recursion trees for $\delta^n(g_i)$ and $\delta^n(g_1)$).
Therefore:
$$[g_i,Y]=\delta^p(g_i)-\delta(g_i)\frac{\delta^p(g_1)}{\delta(g_1)}=\left(\xi_i^p-\xi_i\right)x^pg_i$$ 
It follows that each $g_i$ fixes $Y$ in $\Delta_0$. Then by Lemma 5.13, $Y$, and more generally $Y-\beta$ for any $\beta\in k$, generates a proper submodule if $\alpha=0$.  This proves the statement in both the $t=1$ and $t=0$ cases for the dimensions of simple modules over the $0$-fiber.

Suppose $\alpha\neq0$. From the computation of the action of $g_i$ above, it is evident that $g_i$ does not fix $Y$ when $\alpha\neq0$, since $\xi_i^p-\xi_i\neq0$. By the same argument given in Proposition 6.15, $g_i$ does not fix any polynomial in $Y$ of degree less than $p^r$. Same thing if $Y$ is replaced by $y^p-y\frac{\delta^p(g_i)}{\delta(g_i)}$. On the other hand, we know that in the $t=0$ case, the image of $C=y^q-y\frac{\delta^q(g_1)}{\delta(g_1)}$ in $\Da$ generates a submodule. So it must be maximal, and the simples over the smooth locus are $q$-dimensional in the $t=0$ case. In the $t=1$ case, $g_1\cdot C=C-yx^q=C-\alpha^qy$, so $C$ does not generate a proper submodule. Likewise, argue by downwards induction on degree in $y$ that $g_1$ will not fix $f(C)$ if $f\in k[C]$ is of degree less than $p$. 

\end{proof}

\bibliographystyle{plain}
\bibliography{andrepaper}
\end{document}